\documentclass[onefignum,onetabnum]{siamart220329}
\usepackage{graphics,graphicx,epstopdf,url}
\usepackage{arydshln}
\usepackage{color}
\usepackage{algorithmic}
\usepackage{url,cases,supertabular}
\usepackage{enumerate,makecell}
\usepackage{graphicx,epstopdf}
\usepackage{color,verbatim}
\usepackage{mathrsfs,subeqnarray,subfigure}
\usepackage{listings,array}
\usepackage{booktabs}
\usepackage{dashrule}
\usepackage{arydshln}
\usepackage{graphicx}
\usepackage{amsfonts}
\usepackage{graphicx}
\usepackage{subfigure}
\usepackage{lineno}
\usepackage{multirow,bm}
\usepackage{color}
\usepackage{algorithm}
\usepackage{algorithmic}
\usepackage{longtable}
\usepackage{booktabs}
\usepackage{rotating}
\usepackage{enumitem}
\usepackage{lineno}

\usepackage{url,cases,supertabular}
\usepackage{amssymb,mathtools}
\usepackage{enumerate,makecell}
\usepackage{amsfonts}
\usepackage{graphicx,epstopdf}
\usepackage{color,verbatim}
\usepackage{mathrsfs,subeqnarray,subfigure}
\usepackage{listings,array}
\usepackage{booktabs}
\usepackage{dashrule}
\usepackage{arydshln}
\usepackage{graphicx}
\usepackage{mathtools}
\usepackage{longtable}
\usepackage{booktabs}
\usepackage{rotating}
\usepackage{enumitem}
\usepackage{lineno}
\usepackage{url}
\usepackage{pdflscape}
\usepackage{makecell}
\usepackage{lipsum}
\usepackage{clrscode}
\usepackage{appendix}
\usepackage{booktabs}
\usepackage{url}
\usepackage{multirow}
\usepackage{textcomp}
\usepackage{amsmath}
\usepackage{amsfonts}
\usepackage{amssymb}
\usepackage{mathrsfs}
\usepackage{hyperref}
\usepackage{graphicx,graphics,subfigure}
\usepackage{hyperref,url}
\usepackage{epsf,epstopdf}
\usepackage{algorithm}
\usepackage{algorithmic}

\numberwithin{equation}{section}

\newtheorem{Assumption}{Assumption}[section]

\newcommand{\ba}{\begin{array}}
\newcommand{\ea}{\end{array}}

\newcommand{\bit}{\begin{itemize}}
\newcommand{\eit}{\end{itemize}}
\newcommand{\be}{\begin{equation}}
\newcommand{\ee}{\end{equation}}
\newcommand{\bee}{\begin{equation*}}
\newcommand{\eee}{\end{equation*}}
\newcommand{\bea}{\begin{eqnarray}}
\newcommand{\eea}{\end{eqnarray}}

\newcommand{\st}{\mathrm{s.t.}}

\newcommand{\argmin}{\mathop{\mathrm{arg\,min}}}

\newcommand{\Rmn}[1]{\uppercase\expandafter{\romannumeral#1}}
\renewcommand{\theenumi}{\roman{enumi}}

\numberwithin{equation}{section}

\newcommand{\Mcal}{\mathcal{M}}

\newcommand{\R}{\mathbb{R}}

\newcommand{\prox}{\mathrm{prox}}
\numberwithin{theorem}{section}
\newcommand{\iprod}[2]{\left \langle #1, #2 \right \rangle }

\usepackage{lipsum}
\usepackage{clrscode}
\usepackage{appendix}
\usepackage{longtable}
\usepackage{booktabs}
\usepackage{url}
\usepackage{multirow}
\usepackage{textcomp}
\usepackage{amssymb}
\usepackage{mathrsfs}
\usepackage{hyperref}
\usepackage{graphicx,graphics,subfigure}
\usepackage{hyperref,url}
\usepackage{epsf,epstopdf}
\usepackage{algorithm}
\usepackage{algorithmic}
\usepackage{dsfont}

\usepackage{indentfirst}
\usepackage{pdfpages}
\usepackage{pdflscape}
\usepackage{longtable}
\usepackage{booktabs}
\ifpdf
  \DeclareGraphicsExtensions{.eps,.pdf,.png,.jpg}
\else
  \DeclareGraphicsExtensions{.eps}
\fi
\newtheorem{rmk}{remark}[section]

\headers{ALPDSN for multi-block composite optimization}{Zhanwang Deng, Kangkang Deng, Jiang Hu, Zaiwen Wen}

\usepackage{amsopn}

\ifpdf
\hypersetup{
  pdftitle={},
  pdfauthor={}
}
\fi

\begin{document}

\title{
An Augmented Lagrangian Primal-Dual Semismooth Newton method for multi-block composite  optimization\thanks{{Z. Wen was supported in part by the NSFC grant 11831002.}
}}

\author{Zhanwang Deng\thanks{Academy for Advanced Interdisciplinary Studies, Peking University, Beijing 100871, CHINA (\email{dzw\_opt2022@stu.pku.edu.cn}) } \and
Kangkang Deng\thanks{Department of Mathematics,  National University of Defense Technology, Changsha, 410073,
CHINA (\email{freedeng1208@gmail.com}).} \and
Jiang Hu\thanks{Corresponding author. Massachusetts General Hospital and Harvard Medical School, Harvard University, Boston, MA
02114, US (\email{hujiangopt@gmail.com}).} \and
Zaiwen Wen\thanks{Beijing International Center for Mathematical Research, Center for Machine Learning Research and College of Engineering, Peking University, Beijing 100871, CHINA (\email{wenzw@pku.edu.cn})}}
\maketitle

\begin{abstract}
In this paper, we develop a novel primal-dual semismooth Newton method for solving linearly constrained multi-block convex composite optimization problems. First, a differentiable augmented Lagrangian (AL) function is constructed by utilizing the Moreau envelopes of the nonsmooth functions. It enables us to derive an equivalent saddle point problem and establish the strong AL duality under the Slater’s condition. Consequently, a semismooth system of nonlinear equations is formulated to characterize the optimality of the original problem instead of the inclusion-form KKT conditions. We then develop a semismooth Newton method, called ALPDSN, which uses purely second-order steps and a nonmonotone line search based globalization strategy. Through a connection to the inexact first-order steps when the regularization parameter is sufficiently large, the global convergence of ALPDSN is established. Under the regularity conditions, partial smoothness, the local error bound, and the strict complementarity, we show that both the primal and the dual iteration sequences possess a superlinear convergence rate and provide concrete  examples where these regularity conditions are met. Numerical results on the image restoration with two regularization terms and the corrected tensor nuclear norm problem are presented to demonstrate the high efficiency and robustness of our ALPDSN.
\end{abstract}


\section{Introduction}
The goal of this paper is to design an efficient primal-dual second-order method for solving the following linearly constrained  multi-block convex composite  optimization problem:
\begin{equation}\label{prob}
\min_{\bm{x}_1,\cdots,\bm{x}_n} \sum_{i \in \mathcal{I}_1} h_i(\bm{x}_i) + \sum_{i \in \mathcal{I}_2} h_i(\bm{x}_i), \; \mathrm{s.t.} \; \sum_{i=1}^{n-1} \mathcal{A}_i \bm{x}_i + \bm{x}_n  = \bm{b},
\end{equation}
where $\mathcal{I}_1 = [1,\cdots,n_1], \mathcal{I}_2 = [n_1+1,\cdots,n]$,
 $0 \leq n_1 \leq n - 1$ is an integer,
$h_i:\mathcal{X}_i \rightarrow (-\infty, +\infty)$ is a convex and twice continuously differentiable function with Lipschitz continuous gradient for any $i \in \mathcal{I}_1$, $h_{i}:\mathcal{X}_i \rightarrow (-\infty, +\infty]$ is a proper closed convex function for any $i \in \mathcal{I}_2$,
 $\mathcal{A}_i: \mathcal{X}_i\rightarrow \mathcal{Y}$ is a linear map for any $i = 1,\cdots,n-1$, $\mathcal{A}_n=\mathcal{I}$ is the identity mapping, $\bm{b} \in \mathcal{Y}$, $\mathcal{X}_i$ and $\mathcal{Y}$ are real finite-dimensional Euclidean spaces and each is equipped with an Euclidean inner product $\iprod{\cdot}{\cdot}$ and its induced norm $\|\cdot\|$. If $n_1 = 0$, all functions $h_i$ are proper closed convex. Note that the partitioning of the indices $[1,\cdots,n]$ varies across different problems. The proper closed convex functions $h_i$ in the objective function give us the flexibility to handle nonsmooth terms or even extra constraints, such as the $\ell_1$ regularization, nuclear norm or indicator function over nonnegative constraints, etc.  While the objective function in \eqref{prob} is separable, the linear constraint couples all variables together.

\subsection{Particular cases}
The form of \eqref{prob}  encompasses a broad range of practical scenarios and real-world applications.
It indeed covers the cases where $\mathcal{A}_n \neq \mathcal{I}$
by adding an additional block $\bm{x}_{n+1}$ and setting an indicator function $h_{n+1}(\bm{x}_{n+1}) = \delta_{\{0\}}(\bm{x}_{n+1})$. One of the  three-term composite optimization problems is given by
\begin{equation}\label{pro:three}
    \min_{x} \;\; f(\mathcal{A}x) + g(\mathcal{B}x) + h(x),
\end{equation}
where $f$ is continuously differentiable,  $g$ and $h$ are two regularization terms, $\mathcal{A}$ and $\mathcal{B}$ are two linear maps. Typical applications include the TV-$\ell_1$ regularization problem.
The dual problem of \eqref{pro:three} can be formulated as:
\begin{equation}\label{pro:three:dual}
    \min_{\lambda,\mu,w} \;\; f^*(\lambda) + g^*(\mu) + h^*(w), \;\; \mathrm{s.t.} \;\; \mathcal{A}^*\lambda + \mathcal{B}^*\mu + w = 0,
\end{equation}
where $\mathcal{A}^*$ and $\mathcal{B}^*$ are the adjoint operators of $\mathcal{A}$ and $\mathcal{B}$, respectively, and $f^*$ and $h^*$ are the Fenchel conjugate functions of $f$ and $h$, respectively.

\subsection{Literature review}

\subsubsection{First-order algorithms}
One of the popular methods for solving multi-block problem \eqref{prob} is the alternating direction method of multipliers (ADMM) and its variants due to its extreme simplicity and efficiency.
Various applications for ADMM in statistical
learning can be found in \cite{boyd2011distributed}.
 An inexact symmetric Gauss-Seidel based majorized semi-proximal ADMM is considered in \cite{li2016schur} for solving the convex minimization problem whose
objective function is a sum of a multi-block quadratic function and a non-smooth function
involving only the first block. Leveraging this technique, an inexact multi-block ADMM-type first-order
method for solving a class of high-dimensional convex composite conic optimization problems to moderate accuracy is presented in \cite{chen2017efficient}. Motivated by the theoretical advances in \cite{chen2017efficient},
  a semi-proximal augmented Lagrangian-based decomposition method for a class of convex composite quadratic conic programs with a primal block-angular structure is proposed in \cite{lam2021semi}.

 Primal-dual algorithms are also a group of
 mainstream methods for solving multi-block problems by decoupling the variables.
 One of the early primal-dual algorithms is the primal-dual hybrid gradient (PDHG) method \cite{zhu2008efficient}. Distinguishing itself from PDHG,   the Chambolle–Pock method \cite{chambolle2011first}  incorporates an additional extrapolation step. Additionally, the proximal alternating predictor–corrector (PAPC) method in \cite{drori2015simple} also serves as a primal-dual algorithm designed for problems involving two-block functions.
When there are three functions, the Condat–Vu approach   \cite{condat2013primal}  is a generalization of the Chambolle–Pock method while the primal–dual fixed-point algorithm (PDFP) \cite{chen2016primal} is a generalization of PAPC.
The PD3O method  \cite{yan2018new} is a generalization of both Chambolle-Pock and PAPC, and it has the advantages of both Condat-Vu and PDFP. The primal-dual algorithms for two and three functions are reviewed in \cite{komodakis2015playing} with a  focus on the equivalent saddle point reformulation.

\subsubsection{Second-order algorithms}

In the context of the two-block composite optimization problem,
 the authors in \cite{zhao2010newton} propose an AL semismooth Newton (SSN) method, called SDPNAL,  to solve semidefinite programming.
 It employs an inexact augmented Lagrangian method (ALM) for its dual problem and utilizes the SSN method to minimize each AL function. Based on this algorithmic framework, the efficiency of  augmented Lagrangian SSN method has also
been verified in sparse optimization \cite{li2018efficiently}, linear programming \cite{li2020asymptotically} and quadratic semidefinite programming \cite{li2018qsdpnal}.

The SSN methods based on the residual mappings of two-block composite optimization problems induced by the first-order methods, e.g., the proximal gradient method \cite{xiao2018regularized,milzarek2014semismooth} and the Douglas–Rachford splitting \cite{li2018semismooth}, have also been shown to be  competitive in solving the structured composite optimization problems.
Furthermore, a globalized stochastic semismooth Newton method for
solving stochastic optimization problems involving smooth nonconvex and nonsmooth convex terms
in the objective function is presented in \cite{milzarek2019stochastic}. Several studies have used this framework to address various sparse regularization terms \cite{chen2020proximal,chen2020alternating}. The key is the derivation of the Jacobian for the proximal operator associated with the regularization term and the effective integration of sparsity into the computation of the semismooth Newton step.

\subsection{Contributions}

Our main contribution is the development of  a versatile semismooth Newton method for solving  general multi-block convex composite optimization problems.
 A few novel perspectives are summarized as follows.

(1) \textbf{An easy-to-follow paradigm of semismooth system.}  We construct a differentiable AL function  for problem \eqref{prob} possibly with nonsmooth regularization terms (including indicator functions). It is a summation of the smooth part of the objective function and the Moreau envelopes of the nonsmooth part. This results in an equivalent AL based saddle point problem such that the strong duality holds under Slater's condition. Consequently, a semismooth system of nonlinear equations is derived to characterize the optimality of the original multi-block problem. This reformulation is an easily accessible routine, but approaches in    \cite{xiao2018regularized,li2018semismooth} heavily depend on the construction of first-order type methods using the specific problem structures.

(2) \textbf{A tractable pure semismooth Newton strategy.}
Since the generalized Jacobian of the above nonlinear system can be assembled in a standard way and is positive semidefinite, we propose a semismooth Newton method, called ALPDSN, with a nonmonotone line search procedure.  We further exploit the structures of the semismooth Newton system and develop a computationally tractable way to solve it efficiently. Note that the convergence of the approaches in  \cite{xiao2018regularized,milzarek2014semismooth,li2018semismooth,milzarek2019stochastic} must be ensured through additional first-order steps. By checking whether the residual has sufficient decrease or the regularization parameter is large enough, ALPDSN is able to only execute second-order steps.
Compared to the outer and inner loops in the ALM-based methods \cite{li2018highly}, our method updates both primal and dual variables simultaneously in a single semismooth Newton step.

(3) \textbf{A rigorous global and local convergence analysis.}
By recognizing the $\alpha$-averaged property of the residual mapping from problem \eqref{prob}, we demonstrate that ALPDSN converges globally  through diminishing residual norms. The central aspects of our proof revolve around the boundedness of residuals and the convergence of the inexact first-order iterative procedure. On the other hand, the nonlinear system is essentially locally smooth under the assumption of partial smoothness of $\{h_i\}_{i \in \mathcal{I}_2}$. It means that ALPDSN transits to a regularized Newton step   near a solution with strict complementarity. Hence, the superlinear convergence is established if the local error bound condition holds for the residual mapping without requiring the nonsingularity of the generalized Jacobian.

(4) \textbf{A promising performance on complicated applications.}
 The concepts of ALPDSN including the regularity conditions in the local  analysis are well demonstrated in two practical and interesting applications, image restoration involving two regularization terms and the corrected tensor nuclear norm (CTNN) problem.  Although there are four groups of variables in them, we only need to solve a semismooth Newton system with two blocks. Numerical comparisons with state-of-the-art first-order methods show the computational efficiency of ALPDSN across both high and low accuracy scenarios. 
\subsection{Organization}
The rest of this paper is organized as follows. In Section \ref{sec2}, we present a nonlinear system from the AL saddle point, where the strong duality holds under Slater’s condition. In Section \ref{sec3}, we propose a semismooth Newton method to solve problem \eqref{prob} based on this nonlinear equation.  A theoretical analysis of the global and local convergence is established in Section \ref{sec4}.  The derivations of ALPDSN on two practical applications are demonstrated  in Section \ref{sec5}. Numerical
experiments are presented in Section \ref{sec6} and we conclude this paper in Section \ref{sec7}.

\subsection{Notations} For a linear operator $\mathcal{A}$, its adjoint operator is denoted as $\mathcal{A}^*$.   The Fenchel conjugate function of a proper convex function $g$ is $g^*(\bm{z}) = \sup_{\bm{x}}\{\left<\bm{x},\bm{z}\right> - g(\bm{x})  \}$ and the subdifferential  is $ \partial g(\bm{x}): = \{\bm{z}:~ g(\bm{y}) - g(\bm{x}) \geq \left<\bm{z}, \bm{y} - \bm{x}  \right>,~\forall \bm{y}  \}. $
 For any proper convex function $g$ and a constant $t>0$, the proximal operator of $g$ is defined by $
    \prox_{tg}(\bm{x}) = \arg\min_{\bm{y}}\{g(\bm{y}) + \frac{1}{2t}\|  \bm{y} - \bm{x}\|^2  \}
$.
For a convex set $\mathcal{C}$, we use $\delta_{\mathcal{C}}$ and ${\rm ri}(\mathcal{C})$ to denote its indicator function and relative interior, respectively.
Other notations will be defined when they appear.

\section{A nonlinear system from the AL saddle point} \label{sec2}

Throughout this paper, we make the following assumption.
\begin{Assumption}\label{assum}
 Problem \eqref{prob} has an optimal solution $\bm{x}^* = (\bm{x}_1^*,\cdots,\bm{x}_n^*)$. Furthermore, the Slater's condition holds, i.e., there exists $\bm{x}_i \in  {\rm ri}({\rm dom}(h_i)) (i=1,\cdots,n) $ such that $\sum_{i=1}^{n-1} \mathcal{A}_i \bm{x}_i + \bm{x}_n = \bm{b} $.
\end{Assumption}

Under Assumption \ref{assum}, there exist  Lagrangian multipliers $\bm{z}^* = \{\bm{z}_i^*\}_{i\in\mathcal{I}_2}$, where $\bm{z}_i^*\in \mathcal{X}_i$ for $i\in \mathcal{I}_2$, such that the optimal solution $\bm{x}^* $ of problem \eqref{prob} satisfies the following Karush–Kuhn–Tucker (KKT) optimality conditions \cite{karush1939minima,kuhn2014nonlinear}:
\begin{equation}\label{kkt}
    \left\{ \begin{aligned}
    0 & =  \sum_{i=1}^{n-1} \mathcal{A}_i \bm{x}_i^* + \bm{x}_n^* - \bm{b}, \bm{z}_i^*  \in \partial h_i(\bm{x}_i^*), \,i \in \mathcal{I}_2 \backslash n ,\quad  \bm{z}_n^*  \in - \partial h_n(\bm{x}_n^*),\\
    0 & =  \mathcal{A}_i^*(\bm{z}_n^*) + \nabla h_i(\bm{x}_i^*), i \in \mathcal{I}_1, \quad 0  =    \mathcal{A}_i^*(\bm{z}_n^*) + \bm{z}_i^*, i \in \mathcal{I}_2 \backslash n,
    \end{aligned} \right.
\end{equation}
where $\partial h_i(\bm{x}_i^*)$ denotes the subdifferential of $h_i$ at $\bm{x}_i^*$.

\subsection{Saddle point problem}

In this subsection, we present the AL function of problem \eqref{prob} and derive a compact form of this function by eliminating variables. 
  We first introduce the auxiliary variables $\bm{y}_i = \bm{x}_i, i \in \mathcal{I}_2 \backslash n$ to decouple the variables for
the constraint $\sum_{i=1}^{n-1} \mathcal{A}_i \bm{x}_i +\bm{x}_n = \bm{b}$ from the possibly nonsmooth term $h_i, i\in \mathcal{I}_2\backslash n$.
Therefore, one can equivalently recast \eqref{prob}  as
\begin{equation}\label{prob:dual}
\begin{aligned}
   \min & \sum_{i \in \mathcal{I}_1}h_i(\bm{x}_i) +  \sum_{i \in \mathcal{I}_2} h_i(\bm{y}_i) + h_n(\bm{x}_{n}), \; \\
   \mathrm{s.t.} & \; \sum_{i=1}^{n-1} \mathcal{A}_i \bm{x}_i + \bm{x}_n = \bm{b}, \; \bm{x}_i = \bm{y}_i, i \in \mathcal{I}_2 \backslash n.
\end{aligned}
\end{equation}
Denote $\bm{x}: =(\bm{x}_1,\cdots,\bm{x}_{n-1}),\bm{y}: =(\bm{y}_{n_1 + 1},\cdots,\bm{y}_{n-1})$ and $\bm{z}: =(\bm{z}_{n_1 + 1},\cdots,\bm{z}_{n})$, then the AL function of \eqref{prob:dual} is given as follows:
\begin{equation}\label{func:al}
\begin{aligned}
\mathcal{L}_{\sigma}&(\bm{x},\bm{x}_n,\bm{y};\bm{z}):=\sum_{i \in \mathcal{I}_1} h_i(\bm{x}_i)+  \sum_{i \in \mathcal{I}_2\backslash n} h_i(\bm{y}_i) + h_n(\bm{x}_n)  - \frac{1}{2\sigma} \sum_{i \in \mathcal{I}_2} \|\bm{z}_i\|^2 \\
  & + \frac{\sigma}{2}\| \sum_{i=1}^{n-1} \mathcal{A}_i \bm{x}_i +\bm{x}_n - \bm{b} + \bm{z}_n/\sigma\|^2   + \frac{\sigma}{2}\sum_{i \in \mathcal{I}_2 \backslash n} \|\bm{x}_i - \bm{y}_i - \bm{z}_i/\sigma\|^2,
\end{aligned}
\end{equation}
where $\sigma>0$ is a penalty parameter.  It is worth mentioning that problem \eqref{prob:dual} can be presented as the following minimax problem:
\begin{equation}\label{class:al}
    \min_{\bm{x},\bm{x}_n, \bm{y}} \max_{\bm{z}} \mathcal{L}_{\sigma}(\bm{x},\bm{x}_n,\bm{y};\bm{z}).
\end{equation}

While the AL function \eqref{func:al} possesses desirable properties, its non-differentiability poses a limitation on the development of second-order algorithms. To alleviate this issue, we introduce an AL function defined by:
\begin{equation}\label{func:saddle}
\Phi(\bm{x};\bm{z}):  = \min_{\bm{y},\bm{x}_n} \mathcal{L}_{\sigma}(\bm{x},\bm{x}_n,\bm{y};\bm{z}).
\end{equation}
Problem \eqref{func:saddle} only contains the primal variable $\bm{x}$  while the others are dual variables $\bm{z}$, that is associated with the nonsmooth term $h_i,i\in \mathcal{I}_2$.   Since the minimization problem in \eqref{func:saddle}  has closed-form solutions:
\begin{align}
\bm{y}_i &= \prox_{h_i/\sigma}(\bm{x}_i - \bm{z}_i /\sigma),\; \forall i \in \mathcal{I}_2 \backslash n, \\ \label{eq:v}
\bm{x}_n &= \prox_{h_n/\sigma}\left(\bm{b} - \bm{z}_n/\sigma - \sum_{i=1}^{n-1}\mathcal{A}_i\bm{x}_i\right).
\end{align}
Plugging \eqref{eq:v} into \eqref{func:al} yields the explicit expression of $\Phi$:
\begin{equation}\label{eq:alfunc}
\begin{aligned}
&\Phi(\bm{x};\bm{z}) = h_n(\prox_{h_n/\sigma}(\bm{b} - \bm{z}_n/\sigma-\sum_{i=1}^{n-1}\mathcal{A}_i\bm{x}_i) )  - \underbrace{\frac{1}{2\sigma}  \sum_{i \in \mathcal{I}_2} \|\bm{z}_i\|^2}_{\rm multipliers}\\
 &  +\sum_{i \in \mathcal{I}_2\backslash n } \underbrace{  \left( h_i(\prox_{h_i/\sigma}(\bm{x}_i - \bm{z}_i/\sigma)) + \frac{\sigma}{2} \|\bm{x}_i - \bm{z}_i/\sigma - \prox_{h_i/\sigma}(\bm{x}_i - \bm{z}_i/\sigma) \|^2   \right) }_{{\rm Moreau~envelope ~of~}h_i}\\
      & + \sum_{i \in \mathcal{I}_1} h(\bm{x}_i) +  \underbrace{ \frac{\sigma}{2}\left\| \sum_{i=1}^{n-1} \mathcal{A}_i \bm{x}_i + \prox_{h_n/\sigma}(\bm{b} - \bm{z}_n/\sigma-\sum_{i=1}^{n-1}\mathcal{A}_i\bm{x}_i) - \bm{b} + \bm{z}_n/\sigma\right\|^2}_{\sum_{i=1}^{n-1} \mathcal{A}_i \bm{x}_i + \bm{x}_n= b}.
\end{aligned}
\end{equation}
   To further simplify the expression,
 the Moreau envelope function for $h_i$ is defined by:  $e_{\sigma}h_i(\bm{x}): = \min_{\bm{y}} h_i(\bm{y}) + \frac{\sigma}{2} \| \bm{y} - \bm{x} \|^2. $
 If no confusion can arise, we also write   $\bm{w} = (\bm{x},\bm{z})$ for simplicity. One can see that $\Phi$ can be rewritten as
\begin{equation*}
 \Phi(\bm{w}) = \sum_{i\in \mathcal{I}_1}h_i(\bm{x}_i) + \sum_{i \in \mathcal{I}_2\backslash n} e_{\sigma}h_i(\bm{x}_i - \frac{\bm{z}_i}{ \sigma}) + e_{\sigma}h_n \left(\bm{b} - \frac{\bm{z}_n}{\sigma} - \sum_{i=1}^{n-1}\mathcal{A}_i\bm{x}_i\right)
  - \frac{1}{2\sigma} \sum_{i \in \mathcal{I}_2} \|\bm{z}_i\|^2.
\end{equation*}
Then the corresponding saddle point problem is given by:
  \begin{equation}\label{prob:minimax}
 \min_{\bm{x}}\max_{\bm{z}} \Phi(\bm{x};\bm{z}).
\end{equation}

\subsection{Strong duality}\label{subsec:strong}
In this subsection, we further analyze the relationship between the saddle point problem \eqref{prob:minimax} and problem \eqref{prob}  and give the strong duality of $\Phi(\bm{w})$ based on the AL function. It follows from the Moreau envelope theorem \cite{beck2017first} that $e_{\sigma}h_i$ is continuously differentiable, which implies that $\Phi$ is continuously differentiable. 

\begin{lemma} \label{Strong:duality}
1) The gradient of the  function $\Phi$ in \eqref{eq:alfunc} is given as follows:
\begin{equation} \label{eq:grad-F}
    \begin{aligned}
    \nabla_{\bm{x}_i} \Phi(\bm{w})  =& \nabla h_i(\bm{x}_i) -  \mathcal{A}_i^*\prox_{\sigma h_n^*}\left( \sigma \bm{b} - \bm{z}_n - \sigma \sum_{i=1}^{n-1}\mathcal{A}_i\bm{x}_i \right) ,\,~~ i \in \mathcal{I}_1, \\
    \nabla_{\bm{x}_i} \Phi(\bm{w})  =&   \prox_{\sigma h_i^*}(\sigma \bm{x}_i - \bm{z}_i) - \mathcal{A}_i^*\prox_{\sigma h_n^*}\left( \sigma \bm{b} - \bm{z}_n - \sigma \sum_{i=1}^{n-1}\mathcal{A}_i\bm{x}_i \right),~~\,i \in \mathcal{I}_2 \backslash n,\\
      \nabla_{\bm{z}_i} \Phi(\bm{w})  =&  -\bm{z}_i/\sigma  - \prox_{\sigma h_i^*}(\sigma \bm{x}_i - \bm{z}_i)/\sigma,~~\, i \in \mathcal{I}_2 \backslash n,  \\
     \nabla_{\bm{z}_n} \Phi(\bm{w}) =& -\bm{z}_n/\sigma + \prox_{\sigma h_n^*}(\sigma \bm{b} + \bm{z}_n - \sigma \sum_{i=1}^{n-1}\mathcal{A}_i \bm{x}_i)/\sigma.
    \end{aligned}
\end{equation}
2) Suppose that Assumption \ref{assum} holds. Given $\sigma > 0$, the strong duality holds for $\Phi(\bm{x},\bm{z})$, namely,
\begin{equation}\label{lemma:strong}
    \min_{\bm{x}}\max_{\bm{z}} \Phi(\bm{x},\bm{z}) = \max_{\bm{z}} \min_{\bm{x}}\Phi(\bm{x},\bm{z}),
\end{equation}
where both sides of \eqref{lemma:strong} are equivalent to problem \eqref{prob}.
\end{lemma}
\begin{proof}

1) The gradient expression \eqref{eq:grad-F} follows from the explicit gradient expression of the Moreau envelope and the
Moreau decomposition \cite{beck2017first} for any proper convex function  $h$ with $t>0$:
$\bm{x}  = \prox_{t h}(\bm{x})  +  t \prox_{t^{-1} h^*}( \bm{x}/t).$
2) To prove the strong duality, we first consider the following equivalent formulation of \eqref{prob:dual}:

\begin{equation}\label{thm:dual}
\begin{aligned}
   \min_{\bm{x},\bm{y}} & \sum_{i \in \mathcal{I}_1 \cup\, n}h_i(\bm{x}_i) + \sum_{i \in \mathcal{I}_2 \backslash n} h_i(\bm{y}_i) + \frac{\sigma}{2}\| \sum_{i=1}^{n-1} \mathcal{A}_i \bm{x}_i + \bm{x}_n - \bm{b}\|^2 + \frac{\sigma}{2}\sum_{i \in \mathcal{I}_2 \backslash n} \|\bm{x}_i -\bm{y}_i\|^2, \; \\
   \mathrm{s.t.} & \; \sum_{i=1}^{n-1} \mathcal{A}_i \bm{x}_i + \bm{x}_n = \bm{b}, \; \bm{x}_i = \bm{y}_i, i \in \mathcal{I}_2 \backslash n,
\end{aligned}
\end{equation}
whose Lagrangian function is $\mathcal{L}_{\sigma}(\bm{x},\bm{x}_n,\bm{y};\bm{z})$. According to Assumption \ref{assum} and the equivalence of \eqref{prob} and \eqref{prob:dual}, there exists  $\bm{x}_i \in {\rm ri}({\rm dom}(h_i)),\, i \in \mathcal{I}_1 \cup n,\, \bm{y}_i \in {\rm ri}({\rm dom}(h_i)),\, i \in \mathcal{I}_2 \backslash n$ and $\bm{x}_i = \bm{y}_i, i \in \mathcal{I}_2 \backslash n$, such that $\sum_{i=1}^{n-1} \mathcal{A}_i \bm{x}_i + \bm{x}_n = \bm{b}, \; \bm{x}_i = \bm{y}_i, i \in \mathcal{I}_2 \backslash n$. It further indicates that Slater’s condition holds for problem \eqref{thm:dual}.
Thus, the strong duality implies
\begin{equation} \label{thm:eqn:s}
       \min_{\bm{x},\bm{x}_n, \bm{y}} \max_{\bm{z}} \mathcal{L}_{\sigma}(\bm{x},\bm{x}_n,\bm{y};\bm{z}) =  \max_{\bm{z}} \min_{\bm{x},\bm{x}_n, \bm{y}} \mathcal{L}_{\sigma}(\bm{x},\bm{x}_n,\bm{y};\bm{z}),
\end{equation}
where both sides of \eqref{thm:eqn:s} are equivalent to the original minimization problem \eqref{prob}.

By introducing auxiliary variables $\bar{\bm{y}}_i = \bm{y}_i - \bm{x}_i$ and $\bar{\bm{x}}_n =  \sum_{i=1}^{n-1} \mathcal{A}_i \bm{x}_i + \bm{x}_n-\bm{b} $ and defining
\begin{equation*}
    \phi( \bm{x}, \bar{\bm{y}},\bar{\bm{x}}_n  ) :=  \sum_{i=1}^{n-1} h_i(\bm{x}_i + \bar{\bm{y}}_i) +  h_n(\bm{b} - \sum_{i=1}^{n-1} \mathcal{A}_i \bm{x}_i + \bar{\bm{x}}_n) +\frac{\sigma}{2}(\|\bar{\bm{y}}\|^2+\|\bar{\bm{x}}_n\|^2),
\end{equation*}
one can see that for any given $\bm{x}$:
\begin{subequations}
    \begin{align}
        &\max_{\bm{z}} \Phi(\bm{x},\bm{z})
         =   \max_{\bm{z}} \min_{\bm{y},\bm{x}_n} \mathcal{L}_{\sigma }(\bm{x},\bm{x}_n,\bm{y};\bm{z})
 =  \max_{\bm{z}}\min_{\bar{\bm{y}},\bar{\bm{x}}_n}  \left\{  h_n(\bm{b}-\sum_{i=1}^{n-1} \mathcal{A}_i \bm{x}_i +\bar{\bm{x}}_n) \right. \notag \\
&+ \sum_{i=1}^{n-1} h_i(\bm{x}_i + \bar{\bm{y}}_i ) + \frac{\sigma}{2}\sum_{i \in \mathcal{I}_2 \backslash n} \| \bm{z}_i/\sigma -\bar{\bm{y}}_i \|^2
    \left. + \frac{\sigma}{2}\|  \bm{z}_n/\sigma -\bar{\bm{x}}_n\|^2 - \frac{1}{2\sigma} \sum_{i \in \mathcal{I}_2} \|\bm{z}_i\|^2 \right\} \notag \\
 &=  \max_{\bm{z}}
 \min_{\bar{\bm{y}},\bar{\bm{x}}_n} \left\{   \phi( \bm{x}, \bar{\bm{y}},\bar{\bm{x}}_n  )- \sum_{i \in \mathcal{I}_2 \backslash n}\left<\bar{\bm{y}}_i,\bm{z}_i\right> - \left<\bar{\bm{x}}_n,\bm{z}_n\right>  \right\}\notag \\
& \leq  \min_{\bar{\bm{y}},\bar{\bm{x}}_n} \max_{\bm{z}} \left\{   \phi( \bm{x}, \bar{\bm{y}},\bar{\bm{x}}_n  )- \sum_{i \in \mathcal{I}_2 \backslash n}\left<\bar{\bm{y}}_i,\bm{z}_i\right> - \left<\bar{\bm{x}}_n,\bm{z}_n\right>  \right\} \notag  = \phi(\bm{x},0,0). \notag
    \end{align}
\end{subequations}
As a consequence, we have
\begin{equation} \label{eqn:strong2}
    \min_{\bm{x}} \max_{\bm{z}} \Phi(\bm{x},\bm{z}) \le \min_{\bm{x}} \phi(\bm{x},0,0) ,
\end{equation}
where the right hand side is the same as \eqref{prob}.
According to the min-max inequality, we also obtain
\begin{equation} \label{eqn:min-max}
\max_{\bm{z}} \min_{\bm{x}}\Phi(\bm{x},\bm{z})
   \le  \min_{\bm{x}}\max_{\bm{z}} \Phi(\bm{x},\bm{z}).
\end{equation}
 Since the left hand side of \eqref{eqn:min-max} is equal to \eqref{prob},
\eqref{lemma:strong} holds according to inequality \eqref{eqn:strong2}. Then both sides of \eqref{lemma:strong} are equivalent to \eqref{prob}.  This proof is completed.
\end{proof}

\subsection{Analysis of nonlinear system}
 We focus on the nonlinear operator:
\begin{equation} \label{eq:F}
     F(\bm{w}) =
     \begin{pmatrix}
          \nabla_{\bm{x}} \Phi(\bm{w})^\top,
         - \nabla_{\bm{z}} \Phi(\bm{w})^\top
     \end{pmatrix}^\top.
\end{equation}
Prior to analyzing the properties of $F$, we give the definition of the generalized Jacobian and the semismoothness for a locally Lipschitz continuous mapping \cite{mifflin1977semismooth}.
\begin{definition} \label{def:Jacobian}
    Let $F$ be a locally Lipschitz continuous mapping. Denote $D_F$ by the set of differentiable points of $F$. The B-subdifferential of $F$ at $\bm{x}$ is defined by
\[
\partial_B F(\bm{w}) := \left\{\lim_{k \rightarrow \infty} J(\bm{w}^k)\, |\,  \bm{w}^k \in D_F, \bm{w}^k \rightarrow \bm{w}\right\},
\]
where $J(\bm{w})$ represents the generalized Jacobian of $F$ at $\bm{w} \in D_F$. The set $\partial F(\bm{w})$ = $co(\partial_B F(\bm{w}))$ is called Clarke’s generalized Jacobian, where $co$ denotes the convex hull.
\end{definition}

\begin{definition}
    Let $F$ be a locally Lipschitz continuous mapping. We say $F$ is semismooth (or strongly semismooth) at $\bm{w}$ if $F$ is directionally differentiable at $\bm{w}$ and for any $\bm{d}$ and $J \in \partial F(\bm{w}+\bm{d})$, it holds that
$
\| F(\bm{w}+\bm{d}) -  F(\bm{w}) - J\bm{d} \| = o(\|\bm{d}\|) ~({\rm or ~}O(\|\bm{d}\|^2)), \;\; \bm{d} \rightarrow 0. $
We say $F$ is semismooth (or strongly semismooth) if $F$ is semismooth (or strongly semismooth) for any $\bm{w}$.
\end{definition}

We also give the definitions of nonexpansive mapping and $\alpha$-averaged mapping.
\begin{definition}
    A mapping $F: \R^n \rightarrow \R^n$ is said to be nonexpansive if $F$ is Lipschitz continuous with modulus $L \leq 1$. A mapping $T$ is called $\alpha$-averaged $(\alpha \in (0,1])$ if there exists a nonexpansive operator $R$ such that $T = (1 - \alpha) I + \alpha R$ with $I$ being the identity mapping. Generally, for a monotone residual mapping $F$, the corresponding fixed point mapping $T = I - F$ may not be averaged. If in addition, $F$ is  $L$-Lipschitz continuous, there exists a $\beta$ such that $T = I - \beta F$ is $\alpha$-averaged for some $\alpha$.
\end{definition}

Overall, we have the following lemma.
 \begin{lemma} \label{lem:property-F}
 Assume that the functions $\{h_i\}_{i \in \mathcal{I}_1}$ are twice continuously differentiable and with Lipschitz continuous gradients. The following assertions hold.
 \begin{enumerate}
     \item[{\rm (i)}]  The nonlinear operator $F(\bm{w})$ is monotone, i.e.,
     \begin{equation}
         \left<\bm{w}_1 - \bm{w}_2,  F(\bm{w}_1) -  F(\bm{w}_2)\right> \geq 0, \quad \forall~ \bm{w}_1,\bm{w}_2.
     \end{equation}
     In addition, $F$ is Lipschitz continuous, i.e., there exists a $L> 0$ such that
     \be
     \| F(\bm{w}_1) - F(\bm{w}_2)\| \leq L\|\bm{w}_1 - \bm{w}_2\|, \quad \forall~ \bm{w}_{1}, \bm{w}_2.
     \ee
     Then each element of $\partial F(\bm{w})$ is positive semidefinite for any $\bm{w}$.
     \item[{\rm (ii)}] If the proximal operators ${\rm prox}_{h_i}$ are (strongly) semismooth, then the map $F$, defined in \eqref{eq:F}, is (strongly) semismooth.
     \item[{\rm (iii)}] $\bm{w}^*$ is the  stationary point of \eqref{prob:dual} if and only if $\bm{w}^*$ satisfies $F(\bm{w}^*) = 0$.
     \item[{\rm (iv)}] For $\alpha \in (0, 1]$, the mapping $T = I - 2\alpha F/L$ is $\alpha$-averaged.
 \end{enumerate}
 \end{lemma}
 \begin{proof}
(i). Since $\Psi$ is convex-concave, the monotonicity of $F$ can be proved by following \cite{nemirovski2004prox}.
The Lipschitz continuity of $F$ comes from the 1-Lipschitz continuity of ${\rm prox}_{h_i/\sigma}$ and ${\rm prox}_{\sigma h_i^*}$ for any $i=n_1 + 1, \ldots, n$, Lipschitz continuity of $\nabla h_i$ for any $i=1, \ldots, n_1$, and the linearity of $\mathcal{A}$. Then it follows from the monotonicity of $F$ \cite[Proposition 2.3]{jiang1995local} that each element of $\partial F(\bm{w})$ is positive semidefinite for any $\bm{w}$.

(ii).  By definition any twice differentiable
functions are all strongly semismooth everywhere \cite{facchinei2003finite}. Since the composition of any two (strongly) semismooth mappings is also (strongly) semismooth \cite{hintermuller2010semismooth},  $F$ is (strongly) semismooth with respect to $J$ if the proximal operator ${\rm prox}_{h_i}$'s ($i \in \mathcal{I}_2$) are (strongly) semismooth.

 (iii). $F(\bm{w}^*) = 0$ is equivalent to $\nabla \Phi(\bm{w}) = 0$, which implies that $\bm{w}^*$ satisfies the first-order stationary condition. According to Lemma  \ref{Strong:duality}, this proof is completed. We should note that since \eqref{prob} is convex, $\bm{x}^*$ is the solution.

 (iv). Since the monotonicity and Lipschitz continuity of $F$, the mapping $I - 2 F/L$ is nonexpansive. Note that $ T = I - \frac{2 \alpha}{L} F = (1-\alpha)I + \alpha(I - \frac{2}{L}F).  $
Thus we have the $\alpha$-averaged property of $T$. Therefore, the proof is completed.
 \end{proof}

Throughout the paper, we assume that the proximal operators ${\rm prox}_{h_i}$ are semismooth. This is the case where $h_i$ is the $\ell_1$/$\ell_{\infty}$ norm function, the projection over a polyhedral set \cite[Example 12.31]{rockafellar2009variational}, or the projection over symmetric cones \cite{sun2002semismooth}.

\section{A semismooth Newton method} \label{sec3}
As shown in Lemma \ref{lem:property-F}(iii), one can solve the system $F(\bm{w}) = 0$ to obtain the solution of \eqref{prob}. Since $F$ is semismooth, we propose a semismooth Newton method for solving the equation. One typical benefit of second-order methods is the superlinear or quadratic local convergence rate.

 \subsection{Jacobian of $F$}
 We first compute the generalized Jacobian of $F(\bm{w})$. Note that for a convex function $h$, its proximal operator $\prox_{th}$ is Lipschitz continuous. Then, by Definition \ref{def:Jacobian}, the following sets:
 \begin{equation}
 \begin{aligned}
 \hat{\mathcal{D}}_{h_i}  &:= \nabla^2 h_i(\bm{x}_i), i\in \mathcal{I}_1,\;     \mathcal{D}_{h_i}  := \partial  \prox_{\sigma h_i^*}(\sigma \bm{x}_i - \bm{z}_i), i \in \mathcal{I}_2 \backslash n, \\
 \mathcal{D}_{h_n} &:=   \partial \prox_{\sigma h_n^*}(\bm{b}+\bm{z}_n/\sigma - \sum_{i=1}^{n-1}\mathcal{A}_i\bm{x}_i)
\end{aligned}
 \end{equation}
 are the Clarke's generalized Jacobian of $\prox_{\sigma h_i^*}(\sigma \bm{x}_i + \bm{z}_i)$, for any $i=1,\cdots,n$.

  We define the following matrix operators:
  \begin{equation} \label{def:H123}
  \mathcal{H}_{\bm{xx}}:=
      \left( \begin{array}{ccc}
           \sigma \hat{\mathcal{D}}_{h_1}  + \sigma \mathcal{A}_1^* \hat{\mathcal{D}}_{h_n}
 \mathcal{A}_1
           & \cdots &  \sigma \mathcal{A}_1^* \hat{\mathcal{D}}_{h_n} \mathcal{A}_{n-1} \\
          \vdots &  \ddots & \vdots \\
           \sigma \mathcal{A}_{n-1}^* \hat{\mathcal{D}}_{h_n} \mathcal{A}_1 &  \cdots &   \sigma \hat{\mathcal{D}}_{h_{n-1}}+ \sigma \mathcal{A}_{n-1}^* \hat{\mathcal{D}}_{h_n} \mathcal{A}_{n-1}
      \end{array} \right),
  \end{equation}
  $$\mathcal{H}_{\bm{xz}} := - \left[\mathrm{blkdiag}\left(\left\{ \hat{\mathcal{D}}_{h_i} \right\}_{i=1}^{n-1}\right),((\mathcal{A}_1^*;\cdots,\mathcal{A}_{n-1}^*) \hat{\mathcal{D}}_{h_n})^{\top} \right],$$
  and
   $$\mathcal{H}_{\bm{zz}} := \mathrm{blkdiag}\left(\left\{\frac{1}{\sigma} I- \frac{1}{\sigma} \hat{\mathcal{D}}_{h_i} \right\}_{i=1}^n\right),$$  where $\hat{\mathcal{D}}_{h_i} \in \mathcal{D}_{h_i}$ for $i \in \mathcal{I}_2$, $\mathrm{blkdiag}$ denotes the block diagonal matrix operator. For any $\bm{w}$, we define the generalized  Jacobian of $F(\bm{w})$ as
 \begin{equation}\label{equ:jaco}
      \hat{\partial} F(\bm{w}) : = \left\{ \left(
    \begin{array}{cc}
    \mathcal{H}_{\bm{xx}}     & \mathcal{H}_{\bm{xz}}  \\
    -\mathcal{H}_{\bm{xz}}^{\top}     & \mathcal{H}_{\bm{zz}}   \\
    \end{array}
    \right):  \hat{\mathcal{D}}_{h_i} \in \mathcal{D}_{h_i} {\rm~for~all~} i \in \mathcal{I}_2 \right\}.
\end{equation}
It follows from \cite[Lemma 1]{chan2008constraint} that $\hat{\partial} F(\bm{w})[\bm{d}] =  \partial F(\bm{w})[\bm{d}], ~\forall \bm{d}$. Based on this equivalence and Lemma \ref{lem:property-F}(i), any element of $\hat{\partial} F(\bm{w})$ is also positive semidefinite. Such positive definiteness is important in the design of semismooth Newton methods.

\subsection{A globalized semismooth Newton method}
Basically, the semismooth Newton method
is designed to be a generalization of the standard Newton method in the sense that the coefficient matrix of the linear system uses an element of the Clarke's generalized Jacobian. One major difference  is that the semismooth Newton direction is not necessarily a descent direction corresponding to certain merit functions.

We now explain the computation of the semismooth Newton step  at the $k$-th iteration. First,
 an element of the Clarke's generalized Jacobian  defined by \eqref{equ:jaco} is taken as $J^k \in \hat{\partial}F(\bm{w}^k)$. Since $J^k$ may only be positive semidefinite, a small correction $\tau_{k,i}I$ is added to $J^k$, where $\tau_{k,i} := \kappa \gamma^i \|F(\bm{w}^k)\|$ for a suitably chosen integer $i$ and constants $\kappa > 0$ and $\gamma > 1$.
Then, we compute the semismooth Newton direction $\bm{d}^{k,i}$ as the solution of the following linear system
\be \label{eq:ssn} (J^k + \tau_{k,i} \mathcal{I}) \bm{d}^{k, i} = -  F(\bm{w}^k). \ee
 The shift term $\tau_{k,i} \mathcal{I}$ guarantees the existence and uniqueness of $\bm{d}^{k,i}$. The corresponding trial semismooth Newton step is then defined by
\be \label{eq:ssn-step}
  \bar{\bm{w}}^{k,i}  = \bm{w}^k + \bm{d}^{k,i}.
\ee

We next present a globalization scheme to ensure convergence using pure semismooth Newton steps. It utilizes both line searches on the shift parameter $\tau_{k,i}$ and the nonmonotone decrease on the residuals $F(\bm{w}^k)$. Specifically, for an integer $\zeta\geq 1$, $\nu \in (0,1)$, and $\beta \in (1/2, 1]$, we increase the integer $i$ starting from $0$ gradually and check the following conditions
\begin{align}
    \|F(\bar{\bm{w}}^{k,i})\|  & \leq  \nu \max_{\max(1, k-\zeta+1) \leq j \leq k}\|F(\bm{w}^j)\|, \label{eq:decrease-1} \\
    \tau_{k,i} & \geq k^{\beta} \label{eq:decrease-2}.
\end{align}
When \eqref{eq:decrease-1} or \eqref{eq:decrease-2} holds, we set $\bm{w}^{k+1} = \bar{\bm{w}}^{k,i}$. The condition \eqref{eq:decrease-1} checks whether a nonmonotone descent property of the residuals holds. It allows a temporary increase of the residuals. The parameter $\zeta$ controls the number of points to be referenced  back. A larger $\zeta$ gives a loose requirement on the acceptance of the semismooth Newton step. If this condition does not hold even for a large  value $i_{\max}$, we set $\tau_{k,i}$ sufficiently large, i.e.,   \eqref{eq:decrease-2} holds, which in conjunction with the Lipschitz continuity of $F$ guarantees the decrease of the residuals $\{F(\bm{w}^k)\}$ through an implicit mechanism. The detailed approach, dubbed as  ALPDSN, is outlined in Algorithm \ref{alg:ssn}.

\begin{algorithm}[htbp]
\caption{A semismooth Newton method  ALPDSN for solving \eqref{prob}.}
\label{alg:ssn}
\begin{algorithmic}[1]
\REQUIRE The constants $\gamma > 1$, $\nu \in (0,1)$, $\beta \in (1/2, 1]$, $\kappa >0$, an integer $\zeta \geq 1$, and an initial point $\bm{w}^0 $, set $k = 0$.
\WHILE {\emph{stopping condition not met}}
\STATE Compute $F(\bm{w}^k)$ and select $J(\bm{w}^k) \in \hat{\partial} F(\bm{w}^k)$.
\STATE Find the smallest integer $i \ge 0$ such that $\bar{\bm{w}}^{k,i}$ in \eqref{eq:ssn-step} satisfies \eqref{eq:decrease-1} or \eqref{eq:decrease-2}.
\STATE Set $\bm{w}^{k+1} = \bar{\bm{w}}^{k,i}$.
\STATE Set $k=k+1$.
\ENDWHILE
\end{algorithmic}
\end{algorithm}

\subsection{Solving the semismooth Newton system}
One of the key parts of Algorithm \ref{alg:ssn}   is the computation of the search direction $\bm{d}^{k,i}$ from  \eqref{eq:ssn}. By omitting the superscripts or subscripts $k$ and $i$ for simplicity, and denoting $\bm{d}_{\bm{x}} = (\bm{d}_{\bm{x}_{1}},\cdots,\bm{d}_{\bm{x}_{n-1}})^{\top} $,
$\bm{d}_{\bm{z}} = (\bm{d}_{\bm{z}_{n_1+1}},\cdots,\bm{d}_{\bm{z}_n})^{\top} $,
the system \eqref{eq:ssn} can be represented by:
\begin{equation}\label{eq:ssn-xz}
\left(
    \begin{array}{cc}
    \mathcal{H}_{\bm{xx}} + \tau \mathcal{I}     & \mathcal{H}_{\bm{xz}}  \\
    -\mathcal{H}_{\bm{xz}}^T     & \mathcal{H}_{\bm{zz}} + \tau \mathcal{I} \\
    \end{array}
    \right) \left(
    \begin{array}{c}
    \bm{d}_{\bm{x}}       \\
    \bm{d}_{\bm{z}}       \\
    \end{array}
    \right) = \left(
    \begin{array}{c}
    -F_{\bm{x}}       \\
    -F_{\bm{z}}    \\
    \end{array}
    \right),
\end{equation}
where $F_{\bm{x}} = (\nabla_{\bm{x}_1}\Phi(\bm{w}),\cdots,\nabla_{\bm{x}_{n-1}}\Phi(\bm{w}))^{\top}$,
$F_{\bm{z}} = (-\nabla_{\bm{z}_{n_1 +1}}\Phi(\bm{w}),\cdots,-\nabla_{\bm{z}_n}\Phi(\bm{w}))^{\top}. $
 We apply a Gaussian elimination to \eqref{eq:ssn-xz}. For a given $\bm{d}_{\bm{x}}$, the direction $\bm{d}_{\bm{z}}$ can be calculated by
 \begin{equation} \label{cor:z}
 \bm{d}_{\bm{z}} = (\mathcal{H}_{\bm{zz}} + \tau \mathcal{I})^{-1}(\mathcal{H}_{\bm{xz}}^{\top}\bm{d}_{\bm{x}}-F_{\bm{z}}).
  \end{equation}
Hence, \eqref{eq:ssn-xz} can be reduced to a  linear system with respect  to $\bm{d}_{\bm{x}}$:
\begin{equation} \label{eqn:simp}
\mathcal{H}_{\bm{r}}\bm{d}_{\bm{x}}
= F_{\bm{r}},
\end{equation}
where  $F_{\bm{r}}:=\mathcal{H}_{\bm{xz}} (\mathcal{H}_{\bm{zz}}+\tau \mathcal{I})^{-1}F_{\bm{z}} - F_{\bm{x}}$ and $\mathcal{H}_{\bm{r}}:= (\mathcal{H}_{\bm{xx}} +
 \mathcal{H}_{\bm{xz}}(\mathcal{H}_{\bm{zz}} + \tau \mathcal{I})^{-1}\mathcal{H}_{\bm{xz}}^{\top} + \tau \mathcal{I}).$ The definition of  $\mathcal{H}_{\bm{xz}}$ in \eqref{def:H123} yields
\begin{equation} \label{eqn:equation}
\begin{aligned}
\mathcal{H}_{\bm{r}}
 &=
      \left( \begin{array}{ccc}
            \overline{\mathcal{D}}_{h_1}  +  \mathcal{A}_1^* \overline{\mathcal{D}}_{h_n}
 \mathcal{A}_1
           & \cdots &   \mathcal{A}_1^* \overline{\mathcal{D}}_{h_n} \mathcal{A}_{n-1} \\
          \vdots &  \ddots & \vdots \\
            \mathcal{A}_{n-1}^* \overline{\mathcal{D}}_{h_n} \mathcal{A}_1 &  \cdots &    \overline{\mathcal{D}}_{h_{n-1}}+  \mathcal{A}_{n-1}^* \overline{\mathcal{D}}_{h_n} \mathcal{A}_{n-1}
      \end{array} \right),
      \end{aligned}
\end{equation}
      where $ \overline{\mathcal{D}}_{h_i} = \sigma \hat{\mathcal{D}}_{h_i} + \tilde{\mathcal{D}}_{h_i}, \tilde{\mathcal{D}}_{h_i} = \hat{\mathcal{D}}_{h_i} (\frac{1}{\sigma}\mathcal{I} - \frac{1}{\sigma} \hat{\mathcal{D}}_{h_i} + \tau \mathcal{I})^{-1} \hat{\mathcal{D}}_{h_i}.  $
     Consequently, we essentially only need to handle $n-1$ variables in \eqref{eqn:simp} for  problem \eqref{prob}.

 For certain problems, the linear system \eqref{eqn:simp} can be solved directly using the sparse Cholesky factorization,  or the Sherman-Morrison-Woodbury formula or some other techniques depending on the specific structures of $\mathcal{H}_{\bm{r}}$. Otherwise,
it is a common  practice to solve \eqref{eqn:simp} using iterative methods such as preconditioned conjugate gradient, QMR \cite{toh2004solving}, etc.  According to \eqref{eqn:equation},  the $i$-th component of the matrix-vector product
$\bm{h}_{\bm{r}} =  \mathcal{H}_{\bm{r}}\bm{d}_{\bm{x}} $ is
$$\bm{h}_{\bm{r}_i} = (\mathcal{H}_{\bm{r}}\bm{d}_{\bm{x}})_i  = \overline{D}_{h_i}\bm{d}_{\bm{x}_i}+\mathcal{A}_{i}^*\overline{\mathcal{D}}_{h_n} \sum_{i=1}^{n-1}\mathcal{A}_i \bm{d}_{\bm{x}_i}. $$
Hence, $\sum_{i=1}^{n-1}\mathcal{A}_i \bm{d}_{\bm{x}_i}$ can be computed firstly and shared among all components.
 A detailed calculation process  is presented in Algorithm \ref{alg:comp}. In addition, if the corresponding solution is sparse or low-rank, then the special structures of $\overline{\mathcal{D}}_{h_i}$ can further be used to improve the computational efficiency.

\begin{algorithm}[t] \label{alg:com}
\caption{The computational process of $\bm{h}_{\bm{r}} = \mathcal{H}_{\bm{r}}\bm{d}_{\bm{x}}$.} \label{alg:comp}
\begin{algorithmic}[1]
\REQUIRE $\mathcal{A}_i,\bm{d}_{\bm{x}_i}, \overline{\mathcal{D}}_{h_i},\, i = 1,\cdots,n-1$ and $\overline{\mathcal{D}}_{h_n}$.
\ENSURE $\bm{h}_{\bm{r}} = \mathcal{H}_{\bm{r}}\bm{d}_{\bm{x}}.$
\STATE Compute $\bm{a}_1 =  \sum_{i=1}^{n-1}\mathcal{A}_i \bm{d}_{\bm{x}_i},$ where $d_{\bm{x}_i}$ denotes the $i$-block of $\bm{d}_{\bm{x}}$.
\STATE Compute $\bm{a}_2 =  \overline{\mathcal{D}}_{h_n} (\bm{a}_1)$.
\STATE Compute $\bm{h}_{\bm{r}_i} = \mathcal{A}_{i}^*(\bm{a}_2) + \overline{D}_{h_i}\bm{d}_{x_i}, i =1,\cdots,n-1$. 
\end{algorithmic}
\end{algorithm}

\section{Convergence analysis} \label{sec4}
In this section, we study both the global  and  local convergence of  ALPDSN (i.e., Algorithm \ref{alg:ssn}).
\subsection{Global convergence}
The global convergence of ALPDSN relies on the conditions \eqref{eq:decrease-1} and \eqref{eq:decrease-2}, which ensure that each iterate either gives a sufficient decrease on $\|F(\bm{w}^k)\|$ or performs like an inexact fixed-point step. As these two conditions may be satisfied alternatively, we need careful investigation on the estimation of the decrease on $\|F(\bm{w}^k)\|$. We have the following convergence theorem of ALPDSN.
\begin{theorem} \label{thm:global-con}
Suppose that Assumption \ref{assum} holds and
 let $\{\bm{w}^k\}$ be a sequence generated by Algorithm \ref{alg:ssn}. Then, it holds
\be \label{eq:con-w} \lim_{k \rightarrow \infty} \; F(\bm{w}^k) = 0. \ee
\end{theorem}
\begin{proof}
As shown in Lemma \ref{lem:property-F}, the residual mapping $F$ defined in \eqref{eq:F} is monotone and $L$-Lipschitz continuous.
Let us define $\tilde{F} = F/L$. Then, the operator $\tilde{T}:= I - \tilde{F}$ is nonexpansive. Without loss of generality, we use $F$ and $T$ to denote $\tilde{F}$ and $\tilde{T}$, respectively. It then follows from Lemma \ref{lem:property-F} (iii) that there exists a $\bm{w}^*$ such that $F(\bm{w}^*) = 0$. Note that the $\lambda$-averaged operator of $T$ is $T_{\lambda} :=(1-\lambda ) I + \lambda T$.  The semismooth Newton update \eqref{eq:ssn-step} can be written by
\[ \bm{w}^k + \bm{d}^{k,i} = \bm{w}^k - (J^k + \tau_{k,i} I)^{-1} F(\bm{w}^k) = T_{\tau_{k,i}^{-1}}(\bm{w}^k) + (\tau_{k,i}^{-1}I - (J^k + \tau_{k,i} I)^{-1}) F(\bm{w}^k). \]
Define $\bm{r}_k = (\tau_{k,i}^{-1}I - (J^k + \tau_{k,i} I)^{-1}) F(\bm{w}^k)$. It holds that
\[ \begin{aligned}
\|\bm{r}_k\| & = \|(J^k + \tau_{k,i} I)^{-1}\left(\tau_{k,i}^{-1} (J^k + \tau_{k,i}I) - I \right) F(\bm{w}^k) \| \\ & \leq \| (J^k + \tau_{k,i} I)^{-1} \tau_{k,i}^{-1} J^k F(\bm{w}^k) \| \leq  \tau_{k,i}^{-2} L \|F(\bm{w}^k)\|.
\end{aligned}\]
Let $\bm{e}_k = \bm{r}_k\tau_{k,i}$. Then, $\|\bm{e}_k\| \leq \tau_{k,i}^{-1} L \|F(\bm{w}^k)\|$.
The non-expansiveness of $T$ yields
\[ \begin{aligned}
& 2\iprod{F(\bm{w}^{k+1}) - F(\bm{w}^{k})}{\bm{w}^{k+1} - \bm{w}^{k}} \\
= & \| T(\bm{w}^{k+1}) - T(\bm{w}^{k}) \|^2 - \|\bm{w}^{k+1} - \bm{w}^{k}\|^2 - \|F(\bm{w}^{k+1}) - F(\bm{w}^{k})\|^2 \\
\leq & - \|F(\bm{w}^{k+1}) - F(\bm{w}^{k})\|^2.
\end{aligned}
 \]
Consequently, if $\tau_{k,i} \geq k^\beta$ and $k \geq 4$, we have
\be \label{eq:residual-increase}
\begin{aligned}
& \| F(\bm{w}^{k+1})\|^2  = \| F(\bm{w}^{k})\|^2 + \|F(\bm{w}^{k+1}) - F(\bm{w}^{k})\|^2 + 2 \iprod{F(\bm{w}^{k+1}) - F(\bm{w}^{k})}{F(\bm{w}^{k})}\\
=& \| F(\bm{w}^{k})\|^2 + \|F(\bm{w}^{k+1}) - F(\bm{w}^{k})\|^2 + 2\tau_k \iprod{F(\bm{w}^{k+1}) - F(\bm{w}^{k})}{ \bm{w}^{k+1} - \bm{w}^{k} -\bm{r}_k}\\
= & \| F(\bm{w}^{k})\|^2 - (\tau_k - 1)\| F(\bm{w}^{k+1}) - F(\bm{w}^{k}) + \frac{\tau_k}{\tau_k - 1}\bm{r}_k \|^2 + \frac{\tau_k^2}{\tau_k - 1}\| \bm{r}_k \|^2 \\
 \leq & \left( 1+ \frac{L^2}{\tau_{k,i}^{3}(1-\tau_{k,i}^{-1})}\right) \|F(\bm{w}^k)\|^2 \leq \left( 1 + \frac{L^2}{k^{2\beta}} \right) \|F(\bm{w}^k)\|^2.
\end{aligned}
\ee
According to the conditions \eqref{eq:decrease-1} and \eqref{eq:decrease-2}, if
\eqref{eq:ssn-step} is performed with $\tau_{k,i} < k^{\beta}$, we have
\be \label{eq:residual-decrease}
\| F(\bm{w}^{k+1})\|^2 \leq \nu^2 \max_{\max(1, k-\zeta + 1) \leq j \leq k}\| F(\bm{w}^{j})\|^2.
\ee
Then we prove the following assertion using mathematical induction: there exists a $c > 0$ such that
\be \label{globalcov:assertion}
\|F(\bm{w}^{k+1})\| \le c\,\Pi_{n=1}^{k+1} (1 + \frac{L^2}{n^{2\beta}}) \|F(\bm{w}^0)\|^2, \quad \forall \; k.
\ee
It is easy to verify for a large enough $c$ that \eqref{globalcov:assertion} holds when $k \le 4$. Suppose for some $K > 4$, it follows that $\|F(\bm{w}^{K})\| \le c\,\Pi_{n=1}^K (1 + \frac{L^2}{n^{2\beta}}) \|F(\bm{w}^0)\|^2$. Then, in the $K$-th iteration, if \eqref{eq:ssn-step} is performed with $\tau_{K,i} \ge K^{\beta}$, we have from \eqref{eq:residual-increase} that $\|F(\bm{w}^{K+1})\| \leq
c\,\Pi_{n=1}^{K+1} (1 + \frac{L^2}{n^{2\beta}}) \|F(\bm{w}^0)\|^2$.
For the remaining case where $\tau_{K,i} \ge K^{\beta}$, by using \eqref{eq:residual-decrease}, we can establish the same upper bound for $\|F(\bm{w}^{K+1})\|$.
Hence, \eqref{globalcov:assertion} holds. By $\Pi_{k=1}^\infty \left(1 + \frac{L^2}{k^{2\beta}} \right) < \exp(\sum_{k=1}^\infty \frac{L^2}{k^{2\beta}} ) < \infty$, we conclude that $\|F(\bm{w}^{k})\|$ is bounded. Let $M > 0$ be the constant such that for all $k$, $\|F(\bm{w}^{k})\| \leq M$. Then, if $\tau_{k,i} \geq k^\beta$, it holds that
\be \label{eq:residual-increase2} \| F(\bm{w}^{k+1})\|^2 \leq \|F(\bm{w}^k)\|^2 + \frac{L^2 M}{k^{2\beta}}. \ee
To prove \eqref{eq:con-w}, let us consider the following three cases, i.e., (i) the total number of steps \eqref{eq:ssn-step} using \eqref{eq:decrease-1} are finite, (ii) the total number of steps \eqref{eq:ssn-step} using \eqref{eq:decrease-2} are finite, (iii) steps \eqref{eq:ssn-step} with both \eqref{eq:decrease-1} and \eqref{eq:decrease-2} are encountered infinitely many times. Specifically, if (i) happens, one can follow \cite[Theorem 1]{davis2016convergence} to show \eqref{eq:con-w}. For the case (ii),  there exists a $K_0 > 0$ such that $\|F(\bm{w}^{k+1})\| \leq \nu \max_{\max(1, k-\zeta +1)\leq j \leq k} \|F(\bm{w}^{j})\|$ for  $k \geq K_0$. It is obvious that \eqref{eq:con-w} holds. In terms of (iii), let $k_1, k_2, \ldots$ be the indices corresponding to the steps \eqref{eq:ssn} with $\tau_{k_j,i} < k_j^\beta$.
By \eqref{eq:residual-increase2}, we have for any $k_n + 1 < j \leq k_{n+1}$,
    \be \label{eq:increment-res}\begin{aligned}
    \|F(\bm{w}^{j})\|^2 & \leq \|F(\bm{w}^{j-1})\|^2 + \frac{L^2 M}{(j-1)^{2\beta}} \leq \|F(\bm{w}^{k_{n}+1})\|^2 + \sum_{k=k_n+1}^{j-1}  \frac{L^2 M}{k^{2\beta}}.
    \end{aligned}
    \ee
    Let $a_n = \max_{\max(1, n- \zeta +1) \leq j \leq n} \|F(\bm{w}^{k_j + 1})\|$ and $\xi_n:= \sum_{k=k_n+1}^{k_{n+1}-1}  \frac{M}{k^{2\beta}}$.
    Then, due to $k_{n+ \zeta} - k_n \geq \zeta$, we have the recursion
    \be \label{eq:est1}  \| F(\bm{w}^{k_{n+1} + 1}) \| \leq \nu a_n + \nu \max_{ \max(1, n- \zeta +1) \leq j \leq n} \xi_j. \ee
    Denote $c_n := a_{n\zeta}$. Inequation \eqref{eq:est1} gives $c_{n+1} \leq \nu c_n + \nu \sum_{j=n \zeta}^{(n+1)\zeta - 1} \max_{ \max(1, n- \zeta +1) \leq j \leq n} \xi_j$.
    It follows from \cite[Lemma 2]{xu2015augmented} that there exist constants $C_1 > 0$ and $C_2 > 0$ such that
    \be \label{eq:an} \sum_{n=0}^{\infty}\|c_{n+1}\|^2 \leq C_1 \sum_{n=0}^{\infty} \left[ \sum_{j=n \zeta}^{(n+1)\zeta-1} \max_{ \max(1, n- \zeta +1) \leq j \leq n} \xi_j \right]^2 + C_2. \ee
    Note that for sufficiently large $n \geq N$,
    \[ \left[ \sum_{j=n \zeta}^{(n+1)\zeta-1} \max_{ \max(1, n- \zeta +1) \leq j \leq n}\xi_j \right]^2 \leq \sum_{j=n \zeta}^{(n+1)\zeta-1} \max_{ \max(1, n- \zeta +1) \leq j \leq n}\xi_j \leq \zeta \sum_{k=k_{(n-1)\zeta }+1}^{k_{(n+1) \zeta}-1}  \frac{M}{k^{2\beta}}. \]
    This implies $ \sum_{n=N}^{\infty} \left[ \sum_{j=n \zeta}^{(n+1)\zeta-1} \max_{ \max(1, n- \zeta +1) \leq j \leq n} \xi_j \right]^2 \leq 2 \zeta \sum_{k=k_{(N-1) \zeta}}^{\infty} \frac{M}{k^{2\beta}} < \infty  $ and $\lim_{n \rightarrow \infty} \|c_n\| = 0$. Furthermore, by the definition of $c_n$, we have $\lim_{j \rightarrow \infty} \|F(\bm{w}^{k_j + 1})\| \rightarrow 0$. Using \eqref{eq:increment-res} and $\beta > \frac{1}{2}$, we conclude that \eqref{eq:con-w} holds.
\end{proof}

\subsection{Local convergence}
The local superlinear convergence of the semismooth Newton method is contingent upon two distinct conditions: the BD regularity (specifically, the nonsingularity of every element within the $B$-Jacobian) or the combination of $F$'s local smoothness and the local error bound condition \cite{hu2022local}. While the BD regularity demands the solution point to be isolated—a rather restrictive condition—the latter does not impose such a limitation and can accommodate nonisolated solution points. Establishing the local smoothness is not straightforward. One needs to incorporate the fundamental concepts in nonsmooth optimization, including the strict complementarity (SC) and the partial smoothness. When the local smoothness holds, our  ALPDSN reduces to a regularized Newton method for solving a smooth nonlinear system locally, which enables us to show the local superlinear convergence by requiring $F$ to satisfy the local error bound condition.

\subsubsection{Local smoothness}
 To derive the local smoothness, we use SC and partial smoothness as sufficient conditions.
We can reformulate \eqref{prob} as
\be \label{eqn:pro}
\min_{\bm{x}_1,\cdots,\bm{x}_n} \phi(\bm{x}):=\sum_{i \in \mathcal{I}_1} h_i(\bm{x}_i) + \sum_{i \in \mathcal{I}_2} h_i(\bm{x}_i) + \delta_{\mathcal{C}}(\bm{x}),
\ee
where $\mathcal{C}:= \{\bm{x}: \sum_{i=1}^{n-1} \mathcal{A}_i \bm{x}_i + \bm{x}_n  = b \}$. We then introduce the concept of SC  \cite[Definition 4.7]{hu2022local} for \eqref{eqn:pro}.

\begin{definition}[Strict Complementarity] \label{def:sc}
For problem \eqref{eqn:pro}, we say SC holds at $\tilde{\bm{x}}$ if $
 0 \in \emph{ri} \left(\partial \phi(\tilde{\bm{x}})\right)$.
\end{definition}

As the summation between nonsmooth functions is involved in $\phi$, we generally do not have the expression of the subgradient of $\phi$ though it is convex. However, if  the Slater's condition in Assumption \ref{assum} holds, by noting its equivalence to the relative interior qualification condition in \cite[Theorem 6.2]{mordukhovich2017geometric},we have
\bee \label{eq:subg-phi} \partial \phi(\tilde{\bm{x}}) = \sum_{i \in \mathcal{I}_1} \nabla h_i(\tilde{\bm{x}}_i) + \sum_{i \in \mathcal{I}_2} \partial h(\tilde{\bm{x}}_i) +  \mathcal{N}_{\mathcal{C}}(\tilde{\bm{x}}),  \eee
where $\mathcal{N}_{\mathcal{C}}(\tilde{\bm{x}})$ represents the normal cone of the convex set $\mathcal{C}$ at $\tilde{\bm{x}}$. Note that
\[ \mathcal{N}_{\mathcal{C}}(\bm{x}) = \left \{\bm{u}: \bm{u}=\sum_{i=1}^{n-1} \mathcal{A}_i^*(\bm{u}_i) + \bm{u}_n , \bm{u}_n \in \mathcal{X}_i,\, i \in \mathcal{I}_2 \backslash n  \right \} \]
and its relative interior is itself due to the linear structure.
Furthermore, according to Definition \ref{def:sc} and \cite[Exercise 2.45]{rockafellar2009variational}, a point $\tilde{\bm{x}}$ satisfies the SC if and only if
\[
\bm{0} \in \sum_{i \in \mathcal{I}_i} \nabla h_i(\tilde{\bm{x}}_i) + \sum_{i \in \mathcal{I}_2}{\rm ri}(\partial h(\tilde{\bm{x}}_i)) + {\rm ri}\left( \mathcal{N}_{\mathcal{C}}(\tilde{\bm{x}}) \right).
\]
Let $(\tilde{\bm{x}}, \tilde{\bm{z}})$ be a solution pair of \eqref{prob}. Recalling the KKT optimality \eqref{kkt}, the SC holds at $(\tilde{\bm{x}}, \tilde{\bm{z}})$ if and only if they satisfy \eqref{kkt} and
\be \label{eq:scz}  \tilde{\bm{z}}_i \in {\rm ri}(\partial h_i(\tilde{\bm{x}}_i)), \; i \in \mathcal{I}_2 \backslash n ,\quad \tilde{\bm{z}}_n \in -{\rm ri}(\partial h_n(\tilde{\bm{x}}_n)) . \ee

Next, let us present the concept of partial smoothness \cite{lewis2002active,drusvyatskiy2014optimality}.
\begin{definition}[$C^p$-partial smoothness]\label{def-psmooth}
Consider a proper, closed, convex function $\phi\!:\R^n\rightarrow \bar{\mathbb{R}}$ and a $C^p$ $(p \ge 2)$ embedded submanifold $\Mcal$ of $\R^n$. The function $\phi$ is said to be $C^p$-partly smooth at $x \in \Mcal$ for $v \in \partial \phi(\bm{x})$ relative to $\Mcal$ if
\begin{itemize}
    \item[(i)] Smoothness: $\phi$ restricted to $\mathcal{M}$ is $C^{p}$-smooth near $x$.
    \item[(ii)] Prox-regularity: $\phi$ is prox-regular at $x$ for $v$.

    \item[(iii)] Sharpness: $\mathrm{par}\, \partial \phi(\bm{x}) = N_{\Mcal} (\bm{x})$, where $\mathrm{par}\, \Omega$ is the subspace parallel to $\Omega$ and $N_{\Mcal} (\bm{x})$ is the normal space of $\Mcal$ at $x$.
    \item[(iv)] Continuity: There exists a neighborhood $V$ of $v$ such that the set-valued mapping $V \cap \partial \phi$ is inner semicontinuous at $x$ relative to $\mathcal{M}.$
\end{itemize}
\end{definition}

As shown in \cite{lewis2002active,hu2022local}, the indicator function of the positive semidefinite cone, the nuclear norm and the piecewise linear functions (e.g., the $\ell_1$ norm and the total variation) are all partly smooth. To derive the local smoothness of $F$ given by \eqref{eq:F}, we give the following assumption.
\begin{Assumption} \label{assum:partsmooth}
    The functions $h_i, i=1, \ldots, n$ are $C^p$-partly smooth.
\end{Assumption}

Then we have the local smoothness of $F$  around a solution of \eqref{prob}.
\begin{lemma} \label{lem:local-smooth}
Suppose that Assumptions \ref{assum} and \ref{assum:partsmooth} hold. For any optimal solution $(\tilde{\bm{x}},\tilde{\bm{z}})$ of \eqref{prob}, if the SC is satisfied, $F$ defined by \eqref{eq:F} is locally $C^{p-1}$-smooth in a neighborhood of $(\tilde{\bm{x}},\tilde{\bm{z}})$.
\end{lemma}
\begin{proof} For a proper, closed and convex function $h$, if it is $C^p$-partly smooth $(p \ge 2)$ at $\tilde{\bm{x}}$ for $\tilde{\bm{v}}\in \mathrm{ri} \big(\partial \phi(\tilde{\bm{x}})\big)$ relative to a $C^p$-manifold $\mathcal{M}$, then, by \cite[Lemma 4.6]{hu2022local}, the proximal mapping $\mathrm{prox}_{t \phi}$ is $C^{p-1}$-smooth near $\tilde{\bm{x}}+t \tilde{\bm{v}}$ for all sufficiently small $t > 0$. Consequently, by the SC \eqref{eq:scz}, we have for proper $\sigma > 0$, the proximal mappings
${\rm prox}_{h/\sigma}(\bm{x}_i + \bm{z}_i/\sigma), i=1, \ldots, n$, are locally smooth around $(\tilde{\bm{x}}, \tilde{\bm{z}})$. By the definition of $F$, we immediately obtain the local $C^{p-1}$-smoothness of $F$.
\end{proof}

\subsubsection{Superlinear convergence}

The local error bound condition is crucial for the analysis of semismooth Newton-type methods in the literature \cite{nocedal1999numerical,fan2004inexact,zhou2005superlinear,yue2019family}. We next give its definition with respect to the system $F(\bm{w}) = 0$.
\begin{definition}
    The local error bound holds for $F$ if there exist $\epsilon > 0$ and $\gamma > 0$ such that for any $\bm{w}$ with ${\rm dist}(\bm{w}, \bm{W}^*) \leq \epsilon$,
    \be \label{eq:eb2} \|F(\bm{w})\| \geq \gamma {\rm dist}(\bm{w}, \bm{W}^*), \ee
    where $\bm{W}^*$ is the solution set of $F(\bm{w}) = 0$ and ${\rm dist}(\bm{w}, \bm{W}^*):=\argmin_{\bm{u}} \|\bm{w} - \bm{u}\|$.
\end{definition}

For a piecewise polyhedral mapping $F$ (i.e., the graph of $F$ is the union of finitely many polyhedral sets), it is shown in \cite[Corollary]{robinson1981some} that $F$ satisfies the local error bound condition. Then, we have the following lemma on the local error bound of $F$ given by \eqref{eq:F}.
\begin{lemma} \label{lem:eb}
    If the functions $h_i, i=1, \ldots, n$ are proper, closed, convex, and piecewise linear-quadratic, then the residual mapping $F$ given by \eqref{eq:F} satisfies the local error bound condition.
\end{lemma}
\begin{proof}
    It follows from \cite[Proposition 12.30]{rockafellar2009variational} that the proximal mappings, ${\rm prox}_{\sigma h_i}$ and ${\rm prox}_{\sigma h_i^*}, i=1,2, \ldots,n$, are piecewise linear. Since the composition between the piecewise linear mapping and the linear mapping is piecewise linear \cite[Proposition 3.55]{rockafellar2009variational}, the residual mapping is piecewise linear. Applying \cite[Corollary]{robinson1981some} gives the local error bound of $F$.
\end{proof}

It is worth noting that the category of piecewise linear-quadratic functions encompasses a variety of instances, e.g., quadratic functions, the indicator function of polyhedral sets, the $\ell_1$-norm, and the total variation. We now show the locally superlinear convergence rate of our semismooth Newton method without assuming the BD regularity of $F$.
\begin{theorem} \label{thm:main}
Suppose that Assumptions \ref{assum} and \ref{assum:partsmooth} with $p \geq 3$ hold and $F$ satisfies the local error bound condition \eqref{eq:eb2}. If $\bm{w}_{\ell}$ is close enough to $\tilde{\bm{w}}\in \bm{W}^*$ where the SC is satisfied,  $\bm{w}^k$ converges to $\tilde{\bm{w}}$ Q-superlinearly. 
\end{theorem}
\begin{proof}
By Lemma \ref{lem:local-smooth}, $F$ is locally $C^2$ smooth around $\tilde{\bm{w}}$. Then, the rest of the proof can be split into the following two steps. Firstly, for a locally $C^2$-smooth $F$, it follows from \cite[Lemma 5.8]{hu2022local}  that there exists a fixed constant $c > 0$ such that $ {\rm dist}(\bm{w}^{k+1}, \bm{W}^*) \leq c {\rm dist}^2(\bm{w}^k, \bm{W}^*). $
Combining this inequality and the error bound condition \eqref{eq:eb2} gives
\[ \|F(\bm{w}^{k+1})\| \leq L{\rm dist}(\bm{w}^{k+1}, \bm{W}^*) \leq  cL  {\rm dist}^2(\bm{w}^k, \bm{W}^*) \leq \frac{cL}{\gamma^2} \| F(\bm{w}^k)\|^2. \]
This means $\|F(\bm{w}^{k+1})\| \leq \nu \|F(\bm{w}^k)\|$ if $\|F(\bm{w}^k)\|$ is sufficiently small (which can be guaranteed if $\bm{w}_{k}$ is close enough to $\tilde{\bm{w}}$). Hence, there is no line search needed in Step 3 of Algorithm \ref{alg:ssn}. Secondly, by similar argument of \cite[Theorem 5.9]{hu2022local}, we can prove that $\{\bm{w}_{k}\}_{k\geq \ell}$ stay in the same neighborhood of $\tilde{\bm{w}}$ with $\bm{w}_{\ell}$ for some $\ell > 0$. Therefore, the superlinear convergence of $\{\bm{w}^k\}$ is established by using a similar argument.
\end{proof}
\begin{rmk}
    In Theorem \ref{thm:main}, we use Assumptions \ref{assum}, \ref{assum:partsmooth} and the local error bound condition. Notably, the partial smoothness and the Slater's condition are commonly encountered in various applications. For a given optimal solution, it becomes straightforward to validate the SC.
    Even though the local error bound condition may appear more restrictive, Lemma \ref{lem:eb} shows that such a condition is satisfied when the functions $\{h_i\}_{i=1}^n$ are piecewise linear-quadratic.
\end{rmk}

\section{Applications} \label{sec5}
In this section, we give the detailed derivations of ALPDSN on two practical applications.
\subsection{Image restoration with two regularizations}
In this subsection, we consider:
\begin{equation} \label{image}
\min_{\bm{u}\in \mathbb{R}^n} \quad \frac{1}{2}\|\bm{A}\bm{u}-\bm{b}\|^2 + \delta_{\mathcal{Q}}(\bm{u}) + \lambda \|\bm{D}\bm{u}\|_1,
\end{equation}
where $\bm{b}\in \mathbb{R}^m$ is the observed term, $\mathcal{Q} = \{ \bm{u}\in \mathbb{R}^n:0 \le \bm{u} \le 255\}   $ is the pixel range,  $\bm{A} \in \mathbb{R}^{m \times n}$ is a given linear operator, $\bm{D}$ is a given image transform operator, such as the TV operator or Wavelet transform, and $\| \cdot \|_1$ represents the $\ell_1$ norm. When $\bm{A}$ is the identity operator, then the corresponding problem is the image denoising problem.   When $\bm{A}$ is the blurring operator, it leads to the image deblurring problem.  When $\bm{A}$ is the Radon or Fourier transform, the corresponding problem is usually the CT/MRI reconstruction problem. As a consequence, problem \eqref{image} contains various kinds of image restoration problems. The dual problem of \eqref{image} is

\begin{equation} \label{dual:image2}
\begin{aligned}
\min_{\bm{y},\bm{\mu},\bm{\nu}} \quad & \frac{1}{2}\|\bm{y}\|^2+\iprod{\bm{y}}{\bm{b}} + \delta^*_{\mathcal{Q} }(\bm{\nu}) + \delta_{\|\cdot\|_{\infty}<\lambda}(\bm{\mu}),\\
\st \quad & \bm{A}^{\top}\bm{y}+\bm{D}^{\top}\bm{s}+\bm{\nu} =0, \quad \bm{s}=\bm{\mu},
\end{aligned}
\end{equation}
where $\|\cdot \|_{\infty}$ denotes the $\ell_{\infty}$ norm.
The corresponding modified AL function of \eqref{dual:image2} is constructed as follows:
\begin{equation}
\begin{aligned}
    \Phi(\bm{y},\bm{s};\bm{q},\bm{u}) & =  \frac{1}{2}\|\bm{y}\|^2+\iprod{\bm{y}}{\bm{b}}  +\frac{1}{2\sigma}\|\mathbb{P}_{\mathcal{Q}} \left( \bm{u} - \sigma(\bm{A}^{\top}\bm{y} + \bm{D}^{\top}\bm{s}) \right)\|^2  \\
   & + \frac{\sigma}{2} \| \psi_{\lambda}\left(\bm{q}/\sigma - \bm{s}    \right)\|^2 - \frac{1}{2\sigma}\left(\|\bm{u}\|^2 + \|\bm{q}\|^2 \right),
\end{aligned}
\end{equation}
where $\bm{y}, \bm{s}$ correspond to the primal variables and $\bm{q}$, $\bm{u}$ are the Lagrangian multipliers, $\psi_{\lambda}(\bm{u})= \text{sign}(\bm{u})\max\left\{|\bm{u}|-\lambda,0\right\}$, $\mathbb{P}_{\mathcal{Q}}$ is the projection operator.
 Then the nonlinear system defined in \eqref{eq:F} has four terms as follows:
\begin{equation*} \label{eq:F-ir}
   F(\bm{w}) = (F_{\bm{y}},F_{\bm{s}},F_{\bm{q}},F_{\bm{u}})=(\nabla_{\bm{y}}\Phi(\bm{w})^\top,\nabla_{\bm{s}}\Phi(\bm{w})^\top,-\nabla_{\bm{q}}\Phi(\bm{w})^\top,-\nabla_{\bm{u}}\Phi(\bm{w})^\top)^{\top}.
\end{equation*}
Based on the deduction of \eqref{eqn:simp}, the simplified symmetry semismooth Newton system is a two-block system where the specifics of each block are
 \begin{equation*} \label{eqn:image-final}
     \begin{aligned}
       (\mathcal{H}_{\bm{r}})_{\bm{yy}} & =  \bm{A}\overline{D}_1 \bm{A}^{\top} + (\tau_1+1) \mathcal{I} , \,   (\mathcal{H}_{\bm{r}})_{\bm{ys}}  = \bm{A}\overline{D}_{1}\bm{D}^{\top},\, (\mathcal{H}_{\bm{r}})_{\bm{ss}} = \bm{D} \overline{D}_1\bm{D}^{\top} + \overline{D}_{2}  + \tau_2 \mathcal{I}, \\
      (F_{\bm{r}})_{\bm{y}} & = -\bm{A}D_{1}(D^{\tau_3}_1)^{-1} F_{\bm{u}} -  F_{\bm{y}},
      (F_{\bm{r}})_{\bm{s}}  =   -\bm{D}D_{1}(D^{\tau_4}_{1})^{-1} F_{\bm{u}}  - D_2  (D^{\tau_4}_{2})^{-1} F_{\bm{q}} - F_{\bm{z}} ,
     \end{aligned}
 \end{equation*}
 where $D_1 \in \partial \mathbb{P}_{\mathcal{Q}}(\bm{u} - \sigma(\bm{A}^{\top}\bm{y} + \bm{D}^{\top}\bm{s})) $ and $D_2 \in \partial  \psi_{\lambda}\left(\bm{q}/\sigma - \bm{s}    \right) $, $\tau_1,\tau_2,\tau_3,\tau_4 $ correspond to the regularization parameter of $\bm{y},\bm{s},\bm{q},\bm{u}$ respectively.
 The definition of $\overline{D}_1, \overline{D}_2, D_1^{\tau_3}$ and $D_1^{\tau_4}$ are the same as in \eqref{eqn:equation}. In our implementation, we use the QMR method to solve the corresponding equation to obtain $\bm{d}_{\bm{y}}$, $\bm{d}_{\bm{s}}$ and then obtain $\bm{d}_{\bm{q}}$, $\bm{d}_{\bm{u}}$ by using \eqref{cor:z}.

 We note that $F$ defined by \eqref{eq:F-ir} satisfies the local error bound condition because all functions involved in \eqref{dual:image2} are piecewise linear-quadratic. In addition, one can verify that $\delta^*_{\mathcal{Q}}$ and $\delta_{\|\cdot \|_{\infty} \leq \lambda}$ are partly smooth due to their piecewise linear characteristics \cite[Example 3.4]{lewis2002active}. It is observed from our numerical experiment in Section \ref{sec:ir} that the converged solution satisfies the SC and ALPDSN has superlinear convergence.

\subsection{The CTNN  Minimization  Model}

\subsubsection{Tensor Notations}

Throughout this paper,  tensors are denoted by Euler script letters, e.g. $\mathcal{X} \in \mathbb{R}^{n_1 \times n_2 \times n_3}$. We use the MATLAB notation $\mathcal{X}(i,:,:), \mathcal{X}(:,i,:)$ and $\mathcal{X}(:,:,i)$ to express its $i$-th  horizontal, lateral, and frontal slices. $\widehat{\mathcal{X}}:=\text{fft}(\mathcal{X},[\,],3)$ to denote the fast Fourier transform (FFT) of all tubes along the third dimension of  $\mathcal{X}$. We can also get $\mathcal{X}$ from $\widehat{\mathcal{X}}$ by using the inverse FFT operation along the third dimension, i.e., $\mathcal{X}=\text{ifft}(\widehat{\mathcal{X}},[\,],3).$ Let $\overline{\mathcal{X}}$ denote the block diagonal matrix of the tensor $\widehat{\mathcal{X}}$, where the $i$-th diagonal block of $\overline{\mathcal{X}}$ is the $i$-th frontal slice $\widehat{\mathcal{X}}^{(i)} \in \mathbb{R}^{n_1 \times n_2}$ of $\widehat{\mathcal{X}}$, i.e., $\overline{\mathcal{X}}:= \text{blkdiag}(\widehat{\mathcal{X}})$.
In order to define the $t$-product \cite{kilmer2011factorization}, we define a block circular matrix from the frontal slices $X^{(i)}$ of $\mathcal{X}$:
\begin{equation}
 \text{bcirc}(\mathcal{X}):= \begin{bmatrix}
     X^{(1)} & X^{(n_3)}& \cdots & X^{(2)} \\
     X^{(2)} & X^{(1)} & \cdots & X^{(3)} \\
     \vdots & \vdots & \ddots & \vdots \\
     X^{(n_3)}& X^{(n_3-1)}& \cdots & X^{(1)}
 \end{bmatrix}.
\end{equation}

\begin{definition}
The $t$-product of $\mathcal{X} \in \mathbb{C}^{n_1 \times n_2 \times n_3} $ and $\mathcal{Y} \in \mathbb{C}^{n_2 \times n_4 \times n_3}$ is a tensor $\mathcal{Z} \in \mathbb{C}^{n_1 \times n_4 \times n_3}$ given by $\mathcal{Z} = \mathcal{X} * \mathcal{Y} :=  {\rm fold}({\rm bcirc}(\mathcal{X}) \cdot {\rm unfold}(\mathcal{Y}))$,  where ${\rm unfold}(\mathcal{X})$ takes $\mathcal{X}$ into a block $n_1n_3 \times n_2$ matrix:
\[
{\rm unfold}(\mathcal{X}):= \begin{bmatrix}
 (X^{(1)})^{\top}, (X^{(2)})^{\top}, \cdots, (X^{(n_3)})^{\top}
\end{bmatrix}^{\top},
\]
and ${\rm fold}$ is the inverse of ${\rm unfold:~}{\rm fold}({\rm unfold}(\mathcal{X})) = \mathcal{X}$.
Moreover, $\mathcal{X}* \mathcal{Y}=\mathcal{Z}$ is equivalent to $\overline{\mathcal{X}} \, \overline{\mathcal{Y}} = \overline{\mathcal{Z}}$.
\end{definition}

The tensor nuclear norm is defined as follows. One can check \cite{zhao2022robust,lu2019tensor} and the references therein for more details.
\begin{definition}
   The tensor nuclear norm of a tensor $\mathcal{X} \in \mathbb{R}^{n_1 \times n_2 \times n_3}$ is defined as $\|\mathcal{X}\|_{\rm TNN} = \frac{1}{n_3} \|\overline{\mathcal{X}}\|_*,$ where $\|\overline{\mathcal{X}}\|_*$ denotes the nuclear norm of $\overline{\mathcal{X}}.$
\end{definition}
\subsubsection{The CTNN model}
In this subsection, we consider the following dual problem of the CTNN model \cite{zhang2019corrected}:

\begin{equation} \label{model:dual:CTNN}
\begin{aligned}
\min_{\bm{v},\mathcal{W},\mathcal{Z}}\,\,& \, \frac{m}{2}\|\bm{v}\|^2 - \iprod{\bm{v}}{\bm{y}} + \delta^*_{\mathfrak{U}}(-\mathcal{W}) + \delta_{\mathfrak{X}}(\mathcal{Z})\\
\st\,\, &\,  \mathcal{D}^*_{\Omega}(\bm{v})  - \mathcal{Z} + \mathcal{W} = - \mu F(\mathcal{X}_m)  ,
\end{aligned}
\end{equation}
where $\mathcal{D}_{\Omega}: \mathbb{R}^{n_1 \times n_2 \times n_3} \rightarrow \mathbb{R}^{|\Omega|}$ is the sampling operator, $\Omega$ is the given sample region, $\bm{y}$ and $F(\mathcal{X}_m)$ is the given observed term, $ \mathfrak{U}:= \left\{\mathcal{X} \in \mathbb{R}^{n_1 \times n_2 \times n_3}| \|\mathcal{X}\|_{\infty} \le c \right\} $. $\mathfrak{X}:= \{\mathcal{X} \in \mathbb{R}^{n_1 \times n_2 \times n_3}| \|\mathcal{X}\|_{\rm op}\le \mu \} $, where $\|\cdot\|_{\rm op}$ is the  dual norm of $\| \cdot\|_{\rm TNN}$ \cite{lu2019tensor}.

 When introducing a slack variable $\mathcal{Z} = \mathcal{Y}$, the AL function defined in \eqref{eq:alfunc} can be constructed as follows:
\begin{equation} \label{CTNN:ALM3}
\begin{aligned}
\Phi(\bm{v},&\mathcal{Y};\mathcal{P},\mathcal{X}) := \frac{m}{2}\|\bm{v}\|^2 -  \iprod{\bm{v}}{\bm{y}} + \frac{\sigma}{2}\| \mathcal{U} * \hat{\mathcal{S}}_t * \mathcal{V}^{\top} \|_F^2\\
&+\frac{1}{2\sigma} \|\prox_{\sigma \delta_{\mathfrak{U}} }\left(\mathcal{X}  +\sigma( F(\mathcal{X}_m) +\mathcal{D}_{\Omega}^*(\bm{v}) -\mathcal{Y}) \right) \|_F^2  -  \frac{1}{2\sigma}\left( \|\mathcal{X}\|_F^2+\|\mathcal{P}\|_F^2  \right ),
\end{aligned}
 \end{equation}
where $\| \cdot \|_F$ denotes the tensor Frobenius norm \cite{lu2019tensor}, $\bm{v}, \mathcal{Y}$ correspond to the primal variables and $\mathcal{P}, \mathcal{X}$ are the Lagrangian multipliers, $\hat{\mathcal{S}}_t = \prox_{\mu\|\cdot\|_1}(\mathcal{S}) ,\,
\mathcal{U} * \mathcal{S} * \mathcal{V}^{\top}$  is the tensor singular value decomposition \cite{kilmer2011factorization} of $\mathcal{Y} - \mathcal{P}/\sigma$, where $\mathcal{V}^{\top}$ is the  conjugate transpose of $\mathcal{V}$.Then
the nonlinear system defined in \eqref{eq:F} has four parts as follows:
\[
F(\bm{w})= (F_{\bm{v}} , F_{\mathcal{Y}}  ,F_{\mathcal{P}},F_{\mathcal{X}} ) =(\nabla_{\bm{v}}\Phi(\bm{w})^{\top},\nabla_{\mathcal{Y}}\Phi(\bm{w})^{\top},-\nabla_{\mathcal{P}}\Phi(\bm{w})^{\top},-\nabla_{\mathcal{X}}\Phi(\bm{w})^{\top})^{\top}.
\]
According to the deduction of  \eqref{eqn:simp}, the simplified semismooth Newton system is a two-block system where the specifics of each block are:

 \begin{equation}
     \begin{aligned}
      & (\mathcal{H}_{\bm{r}})_{\bm{vv}}  =    \mathcal{D}_{\Omega} \overline{D}_1 \mathcal{D}_{\Omega}^* + (\tau_1+1) \mathcal{I}, ~~   (\mathcal{H}_{\bm{r}})_{\bm{v}\mathcal{Y}} = - \mathcal{D}_{\Omega} \overline{D}_{1} ,~~ (\mathcal{H}_{\bm{r}})_{\mathcal{Y}\mathcal{Y}} = \overline{D}_1 + \overline{D}_{2}  + \tau_2 \mathcal{I}, \\
      & (F_{\bm{r}})_{\bm{v}}  = \mathcal{D}_{\Omega}D_{1}(D^{\tau_3}_1)^{-1} F_{\mathcal{X}} -  F_{\bm{v}},\,
      (F_{\bm{r}})_{\mathcal{Y}}  =   -D_{1}(D^{\tau_4}_{1})^{-1} F_{\mathcal{X}}  - D_2  (D^{\tau_4}_{2})^{-1} F_{\mathcal{P}} - F_{\mathcal{Y}} ,
     \end{aligned}
 \end{equation}
 where $D_1 \in \partial \prox_{\sigma \delta_{\mathfrak{U}}}\left(\mathcal{X}  +\sigma( F(\mathcal{X}_m) +\mathcal{D}_{\Omega}^*(\bm{u}) -\mathcal{Y}) \right) , D_2 \in \partial  \left( \mathcal{U} * \hat{\mathcal{S}}_t * \mathcal{V}^{\top} \right)$,  $\tau_1,\tau_2, \tau_3$ and $\tau_4 $ correspond to the regularization parameters of $\bm{v},\mathcal{Y},\mathcal{P},\mathcal{X}$, respectively. The definitions of $\overline{D}_1, \overline{D}_2, D_1^{\tau_3}$ and $D_2^{\tau_4}$ are the same as in \eqref{eqn:equation}.

In our implementation, we use the QMR method to solve the corresponding equation  to obtain $\bm{d}_{\bm{v}}$  and $\bm{d}_{\mathcal{Y}}$, then we obtain $\bm{d}_{\mathcal{X}}$ and $\bm{d}_{\mathcal{P}}$ according to  \eqref{cor:z}.Next, we give a description of the operation of $D_2$. To the best of our knowledge, this is probably the first time that the generalized Jacobian of $\mathcal{U} * \hat{\mathcal{S}}_t * \mathcal{V}^T$ with respect to $\mathcal{X}$ has been defined.
Since the FFT of a tensor is of block diagonal form \cite{zhao2022robust,lu2019tensor}, we can extend the generalized Jacobian of singular value thresholding of matrices into tensors. More specifically,  let us denote $\mathcal{V} = [\mathcal{V}_1,\, \mathcal{V}_2],\, \mathcal{V}_1 \in \mathbb{R}^{n_2 \times n_1 \times n_3} , \mathcal{V}_2 \in \mathbb{R}^{n_2 \times (n_2 -n_1) \times n_3}$.  Then applying $D_2$ to a tensor $\mathcal{G} \in \mathbb{R}^{n_1 \times n_2 \times n_3}$ yields
\begin{equation} \label{tensor:D2}
\begin{aligned}
    D_2(\mathcal{G}) &= \mathcal{U}*\ \left[    \frac{\Omega_{\sigma,\sigma}^{\mu} + \Omega_{\sigma,-\sigma}^{\mu}}{2} \odot \mathcal{G}_1 +   \frac{\Omega_{\sigma,\sigma}^{\mu} - \Omega_{\sigma,-\sigma}^{\mu}}{2} \odot \mathcal{G}_1^{\top}, (\Omega_{\sigma,0}^{\mu} \odot\left( \mathcal{G}_2 \right)) \right] *\mathcal{V}^{\top},
\end{aligned}
\end{equation}
where $\sigma =[\sigma^{(1)},\cdots,\sigma^{(n_3)}]$  is the tensor singular value of $\mathcal{X}$,  $ \mathcal{G}_1 = \mathcal{U}^{\top} * \mathcal{G} * \mathcal{V}_1 \in \mathbb{R}^{n_1 \times n_1 \times n_3}, \mathcal{G}_2 = \mathcal{U}^{\top} * \mathcal{G} *\mathcal{V}_2 \in \mathbb{R}^{n_1 \times  (n_2 - n_1) \times n_3 }$ and for the $k$-th frontal slices, $(\Omega_{\sigma,\sigma}^{\mu})^{(k)}$ is defined by:
\begin{equation}
    (\Omega_{\sigma,\sigma}^{\mu})_{ij}^{(k)} := \begin{cases}
  \partial_B \prox_{\mu \| \cdot \|_1}(\sigma_i^{(k)}), & \mbox{if } \sigma_i^{(k)} = \sigma_j^{(k)}, \\
  \left\{ \frac{\prox_{\mu\|\cdot\|_1}(\sigma_i^{(k)}) -\prox_{\mu\|\cdot\|_1}(\sigma_j^{(k)}) }{\sigma_i^{(k)}-\sigma_j^{(k)}} \right\}, & \mbox{otherwise}.
\end{cases}
\end{equation}
 $\Omega_{\sigma,-\sigma}$ and $\Omega_{\sigma,0}$ are defined analogously.  Through simple mathematical derivations, $\overline{D}_2(\mathcal{G})$ and $(D_2^{\tau_4})^{-1}(\mathcal{G})$ can be computed likewise.Denote $\alpha = \{1,\cdots, r \} ,\bar{\alpha} = \{1,\cdots,n_1\}\backslash \alpha , \beta = \{1,\cdots,n_2 \}\backslash (\alpha \cup \bar{\alpha}) ,$ where $r$ is the tubal rank. Since the lower right corners of $\Omega_{\sigma,\sigma}^{\mu}$ and $\Omega_{\sigma,-\sigma}^{\mu}$ are zero,  we can only store and use $(\mathcal{G}_1)_{\alpha\alpha},(\mathcal{G}_1)_{\alpha\bar{\alpha}},(\mathcal{G}_1)_{\bar{\alpha}\alpha}
$,$(\mathcal{G}_1)_{\alpha\beta}$ in computation. Then the low-rank structure can be fully exploited and the total computational cost for each inner iteration is reduced to $O(n_1n_3r^2)$.

\section{Numerical experiments}  \label{sec6}
In this section, we conduct proof-of-concept numerical experiments on ALPDSN. We implement
ALPDSN under MATLAB R2021b. All experiments are performed on a Linux server
with a sixteen-core Intel Xeon Gold 6326 CPU and 256 GB memory.

\subsection{Image restoration with two regularizations} \label{sec:ir}

\subsubsection{CT image restoration}
\begin{figure}[htbp]
\centering
\begin{minipage}{0.48\textwidth}
    \centering
\subfigure[Radon transforms of Phantom]{\includegraphics[width=0.45\textwidth]{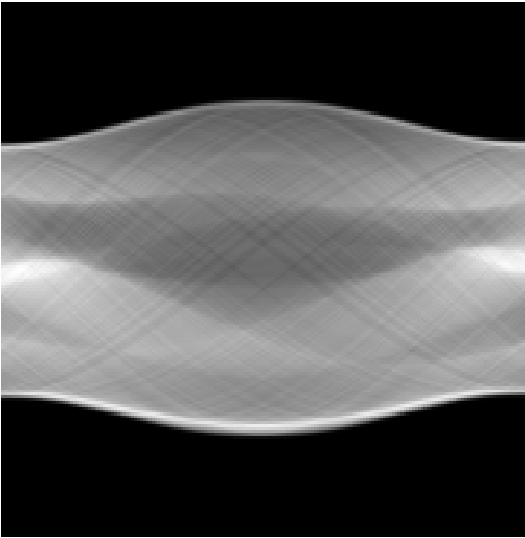}}
    \hspace{0.1in} %
\subfigure[ Sampled
Radon transforms]{\includegraphics[width=0.45\textwidth]{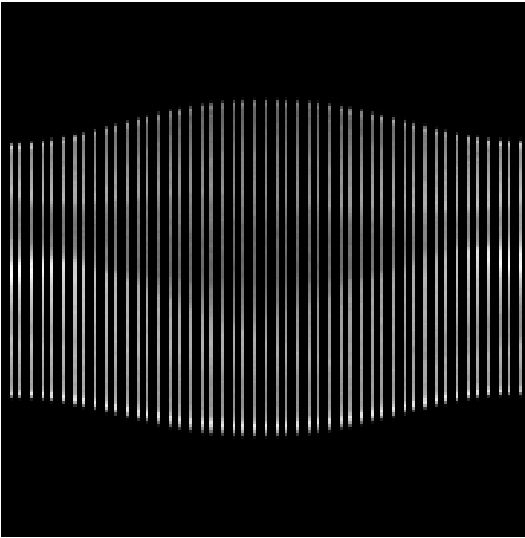}}
\end{minipage}
\hfill
\begin{minipage}{0.48\textwidth}
    \centering
\subfigure[Restored image by FBP ]{\includegraphics[width=0.45\textwidth]{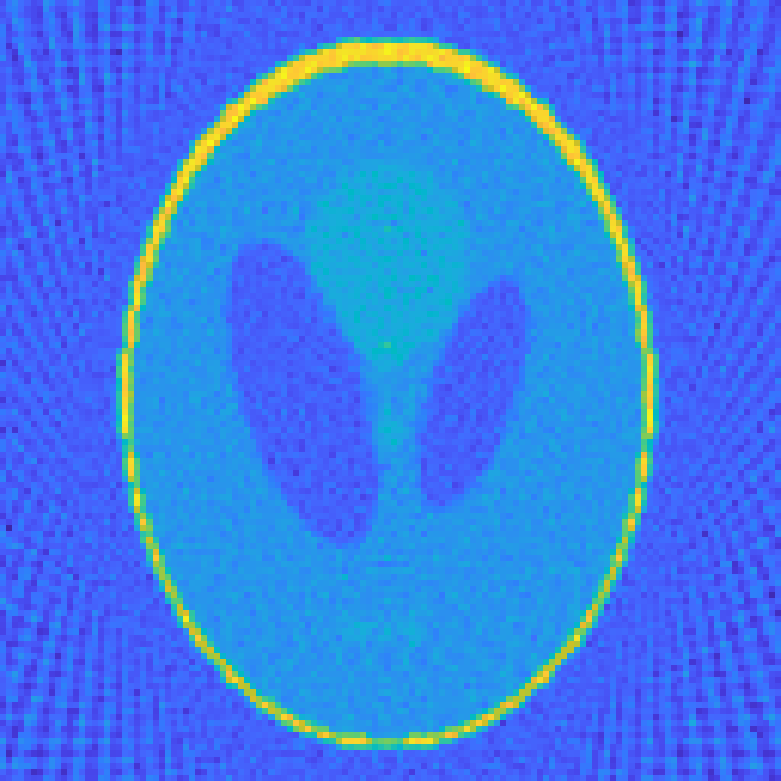}}
    \hspace{0.1in} %
\subfigure[Restored image by ALPDSN]{\includegraphics[width=0.45\textwidth]{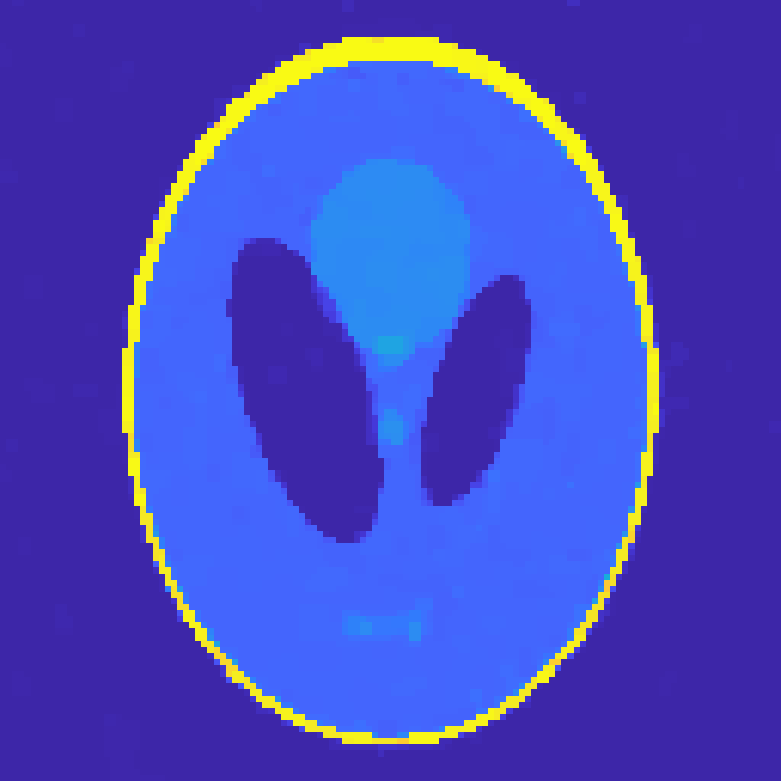}}
\end{minipage}
\caption{Phantom: the PSNR of FBP and  ALPDSN are $19.26$, $35.37$, respectively.}
\label{fig:FBP1}
\end{figure}

\begin{figure}[htbp]
\centering
\begin{minipage}{0.48\textwidth}
    \centering
\subfigure[Radon transforms of Carcinoma]{\includegraphics[width=0.45\textwidth]{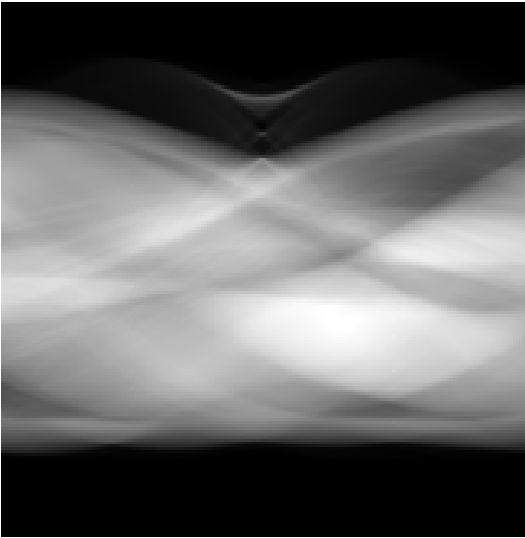}}
    \hspace{0.1in} %
\subfigure[Sampled
Radon transforms]{\includegraphics[width=0.45\textwidth]{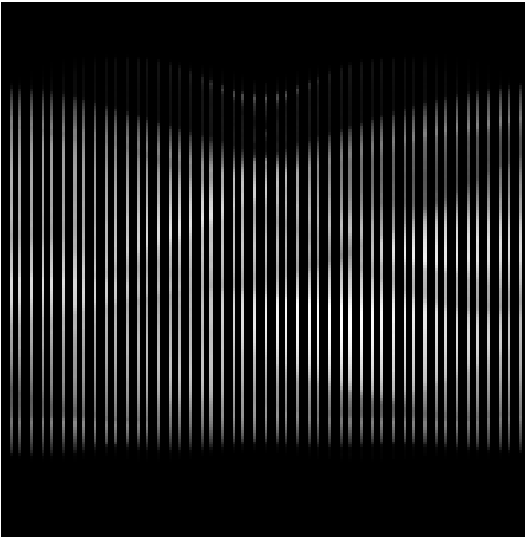}}
\end{minipage}
\hfill
\begin{minipage}{0.48\textwidth}
    \centering
\subfigure[Restored image by FBP  ]{\includegraphics[width=0.45\textwidth]{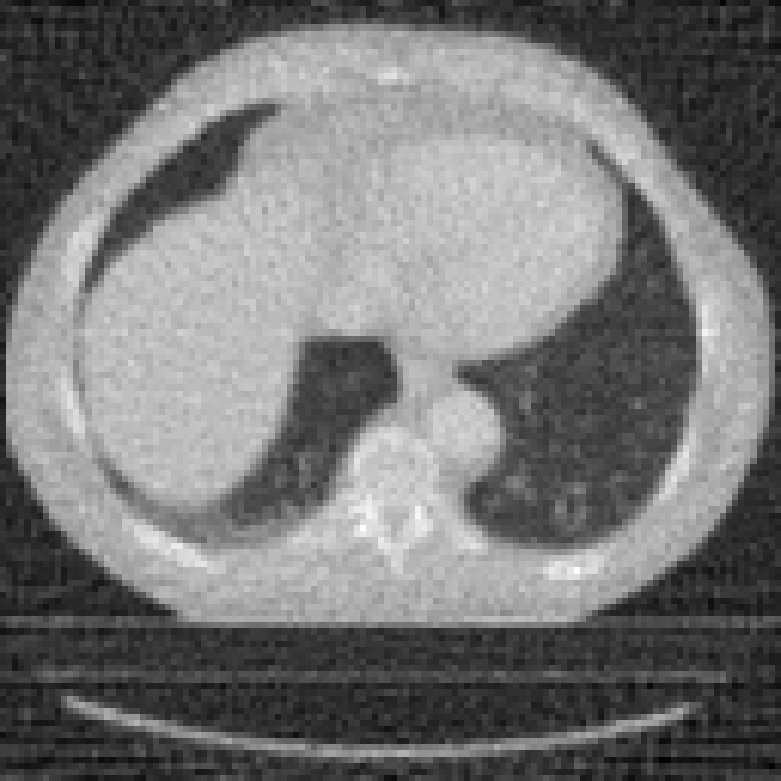}}
    \hspace{0.1in} %
\subfigure[Restored image by ALPDSN ]{\includegraphics[width=0.45\textwidth]{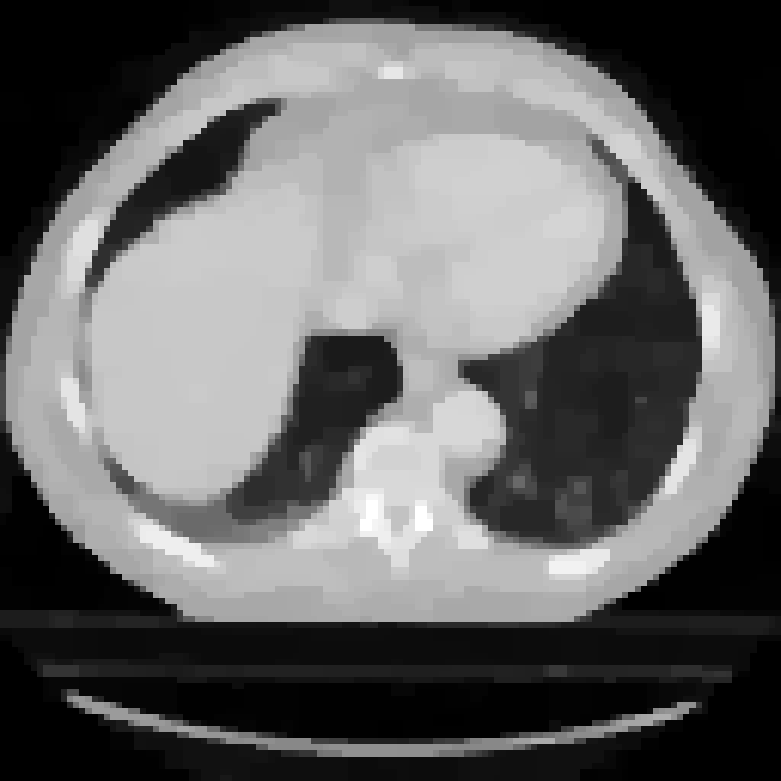}}
\end{minipage}
\caption{Carcinoma:
 the PSNR of FBP and  ALPDSN are $24.53$, $31.57$, respectively.} \label{fig:FBP2}
\end{figure}

In this example, we test the CT image restoration problem. We  generate the data from 180 angles at
$1^{\circ}$ increments from  $1^{\circ}$  to  $180^{\circ}$ with equidistant parallel X-ray beams and then choose 50 angles randomly.
The sizes of tested classic images Phantom and Carcinoma\footnote{\href{https://www.kaggle.com/datasets/kmader/siim-medical-images?resource=download}{ https://www.kaggle.com/datasets/kmader/siim-medical-images}} are 128 $\times$ 128.
Then the corresponding operator matrix
$\bm{A} \in \mathbb{R}^{9250 \times 16384}$ represents a sampled radon transform, while $\bm{D}: \mathbb{R}^{128 \times 128} \rightarrow \mathbb{R}^{256 \times 128}$ corresponds to the anisotropic total variation regularization. The sampling process and the images recovered by the filter-backprojection (FBP) method are visualized in Figures \ref{fig:FBP1} and \ref{fig:FBP2} which demonstrate the effectiveness of our model.  In numerical experiments, we choose $\lambda = 0.1,0.05$.   We
define the relative function error $\mbox{ferror}:= \frac{f-f^*}{f}$  between the current objective function value $f$ and the optimal function value $f^*$, where $f^*$ is obtained by stopping PD3O for 10000 iterations. We compare our algorithm with PD3O \cite{yan2018new},   PDFP \cite{chen2016primal},  and Condat-Vu  \cite{condat2013primal}.
Based on the convergence requirements therein and observed empirical performance, we set the step sizes $\gamma$ of these algorithms to 1.8, 1.5, and 1, respectively.

In consistency with \cite{yan2018new} and to fully demonstrate the efficiency and comparison of these algorithms, we stop the algorithms when $\mbox{ferror} < 10^{-7}$ or the maximal number of iterations, $8000$, is reached. The numerical results
are shown in Table \ref{tab:image1}, where  ``iter'' denotes the number of iterations, ``num'' represents the number of matrix-vector multiplications of $\bm{Au}$, and ``time'' is the total wall-clock time. We present the iterative trajectories of all algorithms in
Figure \ref{fig:image1}, where ``uerror'' refers to $\frac{\|\bm{u}-\bm{u}^*\|_F}{\|\bm{u}^*\|_F}$  and $\bm{u}^*$ is the optimal solution. We observe that ALPDSN requires significantly less time than the other three algorithms to achieve a satisfactory PSNR. Additionally, in terms of ferror and uerror, ALPDSN converges the fastest.Moreover, ALPDSN exhibits the locally superlinear convergence on the ``uerror''.

{\renewcommand{\arraystretch}{0.9}{
\begin{table}[t]
\caption{
Summary on the CT image restoration problem}\label{tab:image1}
\begin{tabular}{|c|c|c|c|c|c|}
\cline{1-6}
Image Name &  &ALPDSN & PD3O & PDFP & Condat-Vu \\ \cline{1-6}

\multirow{4}{*}{Phantom($\lambda = $0.1)}& iter & 346 &2689&3196 &4678 \\ \cline{2-6}
 & num &2506 & 2689 & 3185 & 4686 \\\cline{2-6}
 & ferror & 9.45e-8 & 9.98e-8 & 9.97e-8 & 9.99e-07  \\ \cline{2-6}
&time &\textbf{14.42}&19.80 & 22.85 & 35.19\\        \cline{1-6}

\multirow{4}{*}{Phantom($\lambda = $0.05)}& iter & 404 &2992&3571 &5312 \\ \cline{2-6}
 & num &2125 & 2992 & 3571 & 5312 \\\cline{2-6}
 & ferror & 7.39e-8 & 9.98e-8 & 9.97e-8 & 9.98e-08  \\ \cline{2-6}
&time &\textbf{14.05}&22.23 & 26.35 & 37.58\\        \cline{1-6}
\multirow{4}{*}{Carcinoma($\lambda = $0.1)}& iter & 380 &3593&4269 &6260\\ \cline{2-6}
 & num&2975 & 3595 & 4269 & 6260\\\cline{2-6}
 & ferror & 7.03e-8 & 9.99e-8 & 9.99e-8 & 9.99e-08  \\ \cline{2-6}
&time &\textbf{18.44} & 29.43 &33.65 & 51.31\\       \cline{1-6}

\multirow{4}{*}{Carcinoma($\lambda = $0.05)}& iter & 361 &4830&5766 &8000 \\ \cline{2-6}
 & num&3971 & 4830& 5766 & 8000 \\\cline{2-6}
 & ferror & 7.42e-8 & 9.99e-8 & 9.99e-8 & 1.71e-07  \\ \cline{2-6}
&time &\textbf{23.80} & 40.30 &46.46 & 67.21\\        \cline{1-6}
\end{tabular}
\end{table}
}

\begin{figure}[htbp]
    \centering
      \subfigure[ferror vs time]{\includegraphics[scale = 0.2]{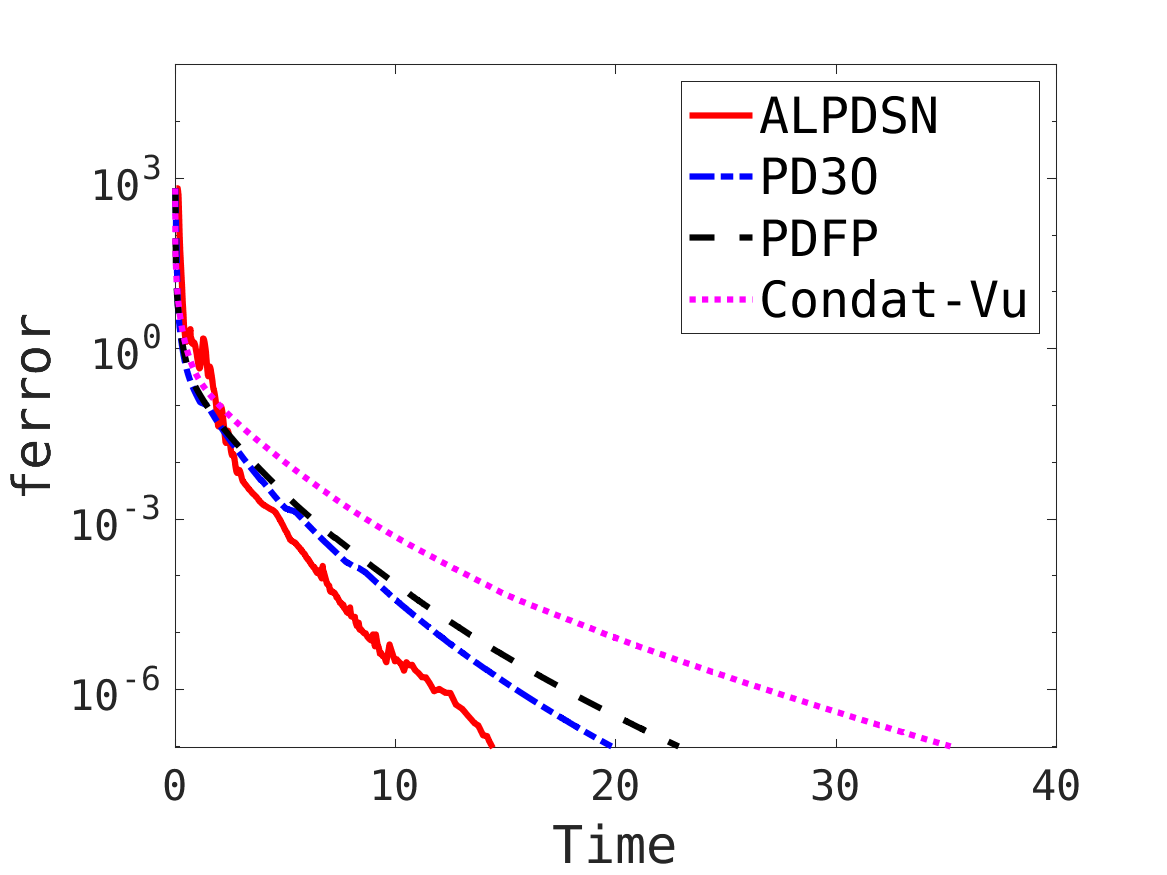}}
      \subfigure[PSNR vs time]{\includegraphics[scale = 0.2]{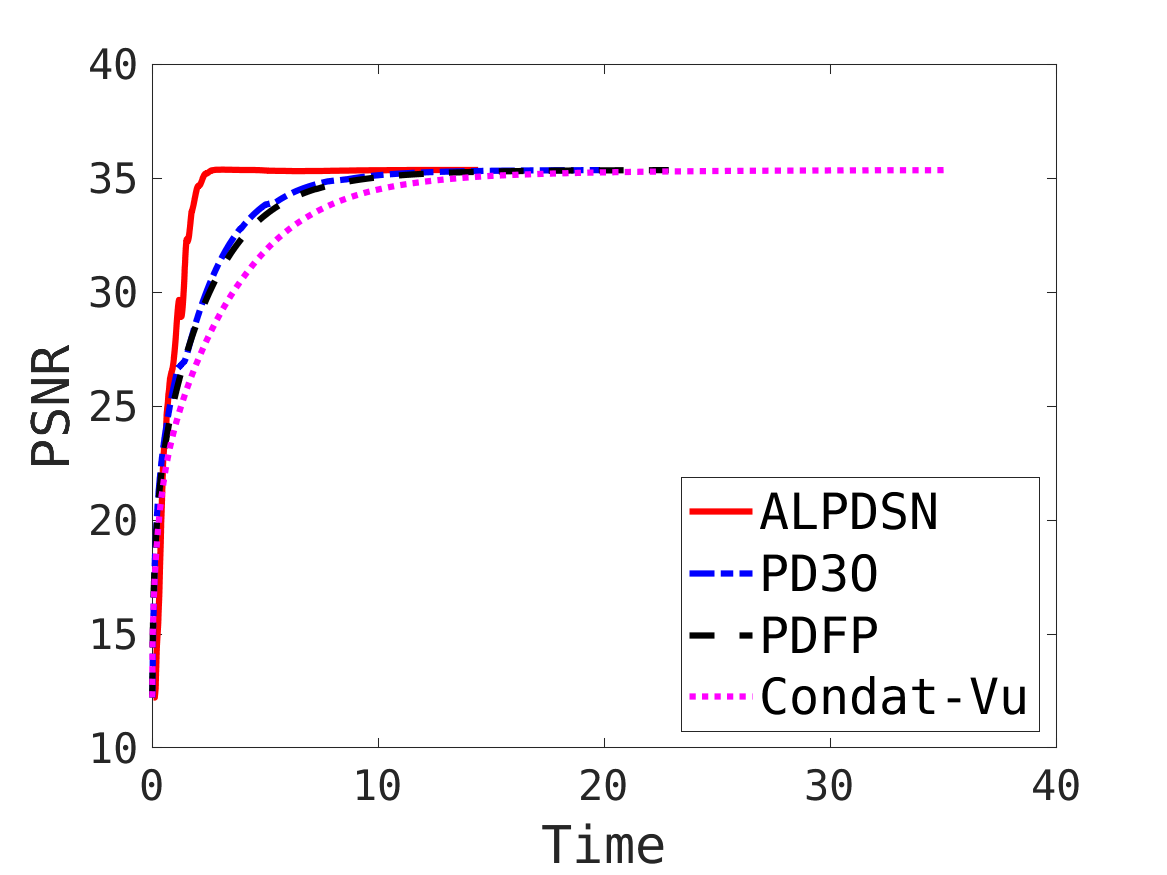}}
      \subfigure[uerror vs time ]{\includegraphics[scale = 0.2]{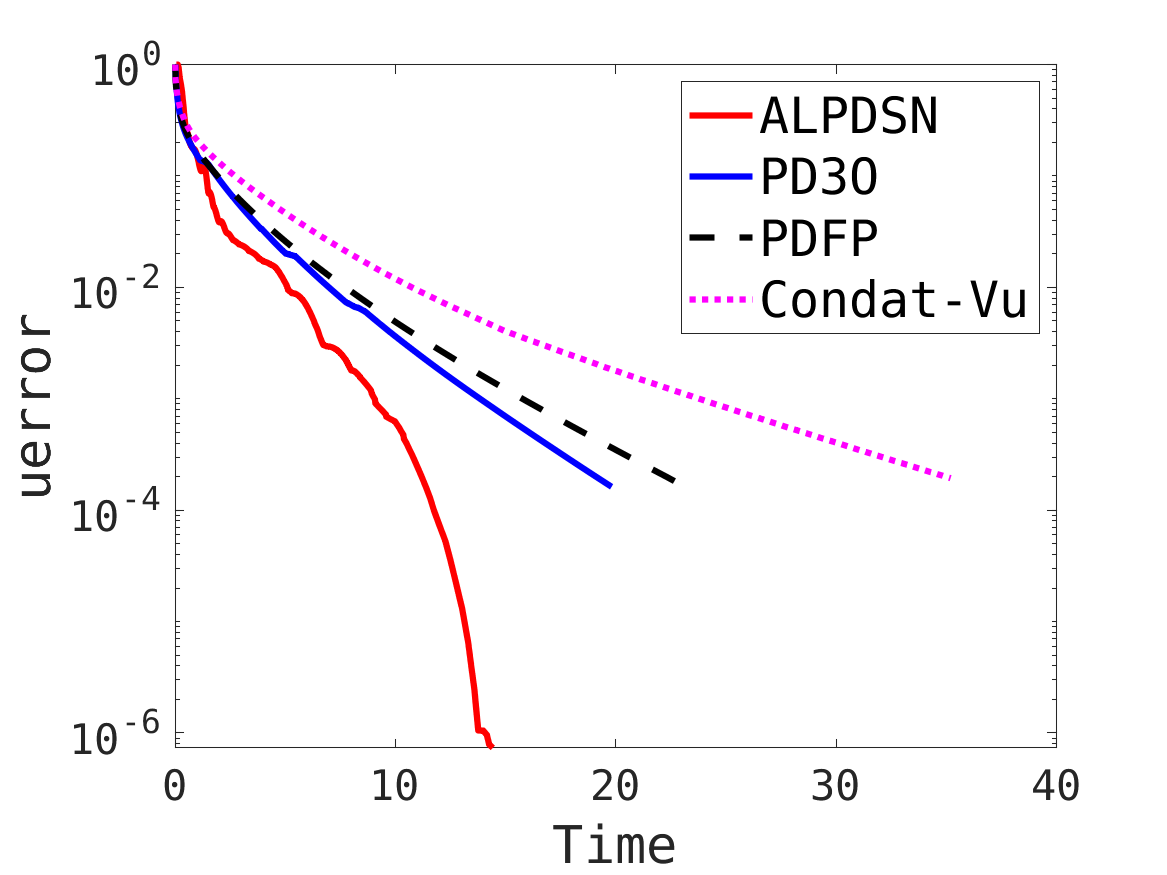}} \\
      \subfigure[ferror vs time]{\includegraphics[scale = 0.2]{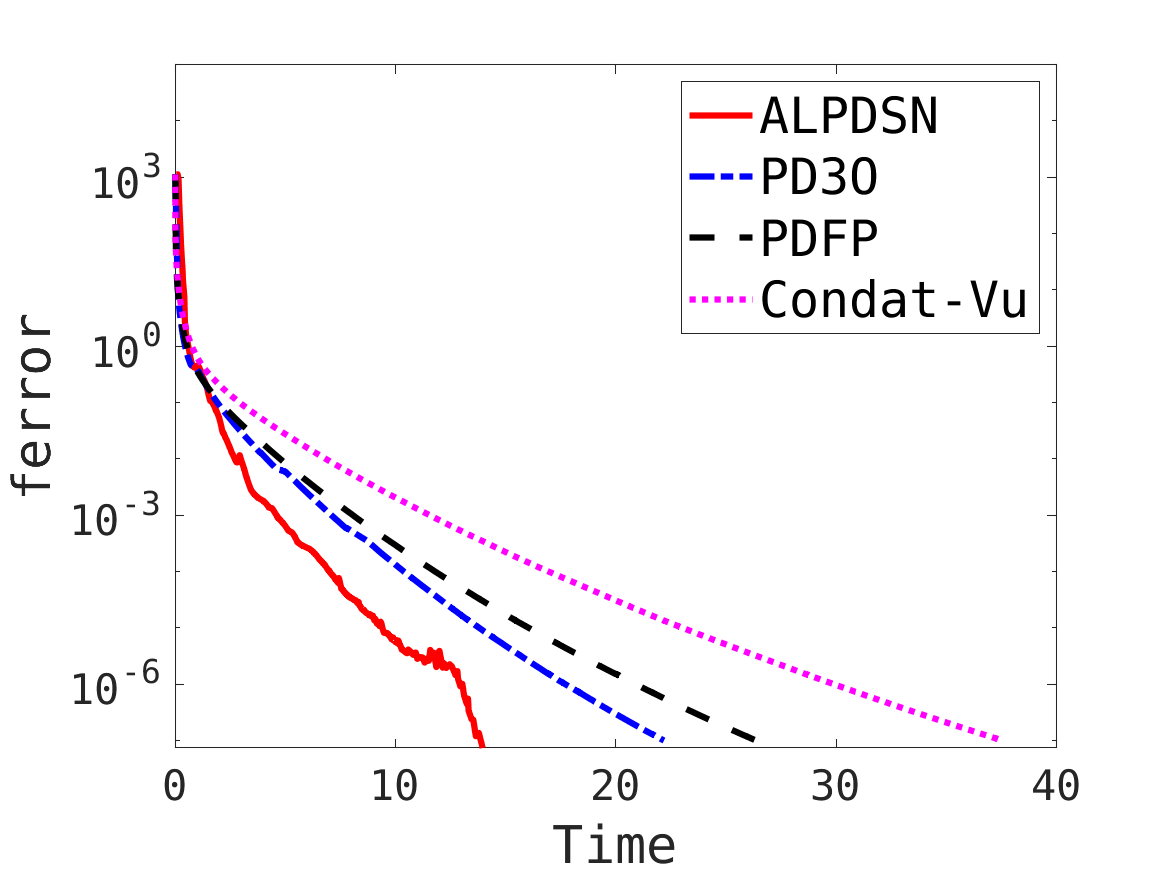}}
      \subfigure[PSNR vs time]{\includegraphics[scale = 0.2]{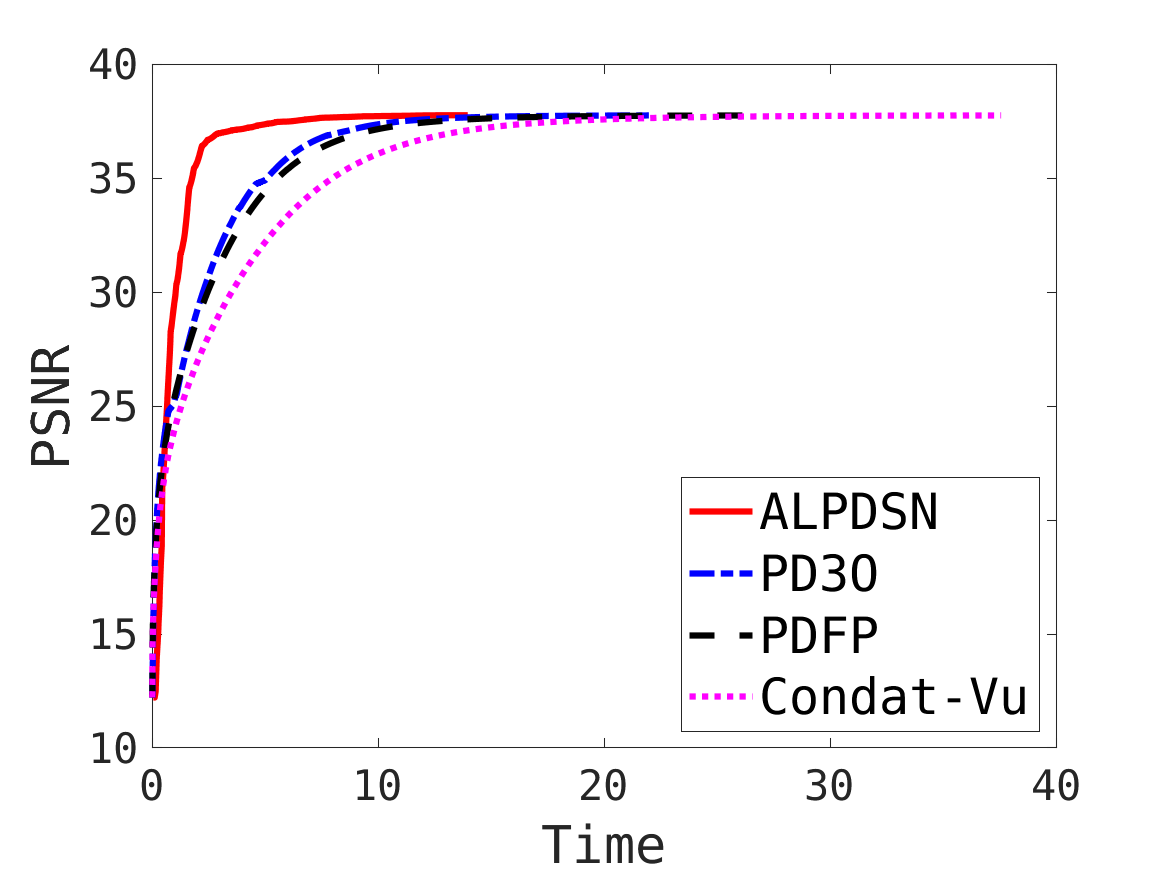}}
      \subfigure[uerror vs time ]{\includegraphics[scale = 0.2]{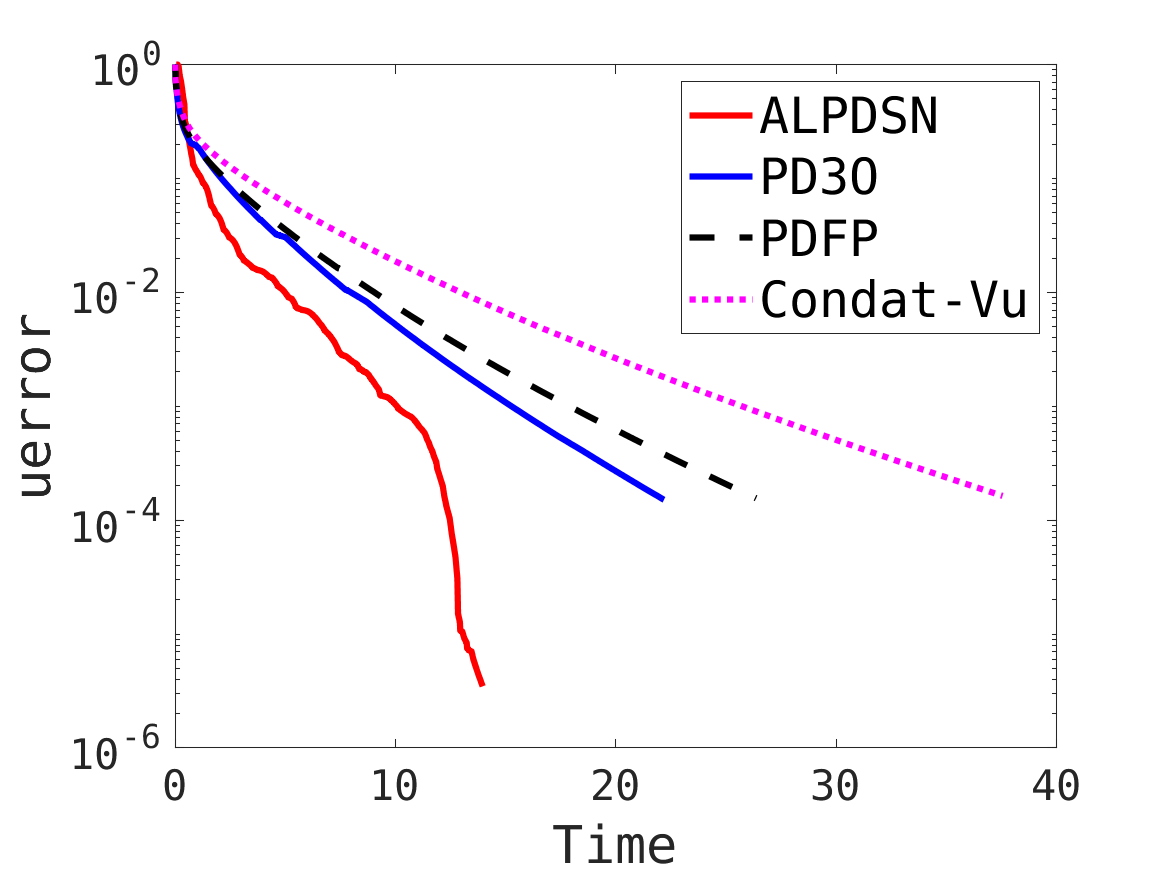}}
    \caption{Comparison between ALPDSN and other first-order methods on Phantom. First row: $\lambda$ = 0.1, second row: $\lambda$ = 0.05.}
      \label{fig:image1}
\end{figure}

\subsubsection{Image deblurring}
 We then test the image deblurring problem, where $\bm{A}$ is the  Gaussian blurring operator. To fully demonstrate the robustness of our algorithm, we test both the gray images and color images. The tested images are Peppers (256 $\times$ 256), Cameraman (256 $\times$ 256), Animals (321 $\times$ 481 $\times$ 3) and Flowers\footnote{\href{https://www2.eecs.berkeley.edu/Research/Projects/CS/vision/bsds/}{https://www2.eecs.berkeley.edu/Research/Projects/CS/vision/bsds/}} (321 $\times$ 481 $\times$ 3). Note that, even though color images are tensors, we treat them as a 963 $\times$ 481 matrix to be consistent with \eqref{image}. The pixels of the tested images are normalized to the range $[0, 255].$
 The data $\bm{b}$ is generated by $\bm{b} = \bm{Au}+ \beta\bm{\varepsilon}$, where $\bm{\varepsilon}$ is the standard Gaussian distribution and $\beta = 1$ is the noise level. The step sizes $\gamma$ of these algorithms are also set to 1.8, 1.5, and 1,  respectively. For the gray images, we stop the iteration when $\mbox{ferror}< 10^{-5}$. We use the stopping criterion, $\mbox{ferror}< 10^{-4}$, for the color images due to the slow tail convergence of first-order algorithms.
The numerical results are shown in Table \ref{tab:image3} and the detailed iterative trajectories are presented in Figure \ref{fig:image4}. In Figure  \ref{fig:image4}(c), we stop ALPDSN  when ferror $< 10^{-7}$ and the superlinear convergence of the ALPDSN algorithm near $\bm{u}^*$ is observed. Similar to the above example, we see that ALPDSN converges the fastest among all algorithms and across both high and low accuracies.
In Table \ref{tab:image3}, ALPDSN always returns a point with lower ferror using a much lower number of matrix-vector multiplications and reduced total wall-clock time.
For gray images, PDFP and Condat-Vu fail to reach the desired accuracy for $\lambda = 0.01$.

\begin{figure}[htbp]
\centering
\begin{minipage}{0.4\textwidth}
    \centering
\subfigure[Peppers]{\includegraphics[width=0.45\textwidth]{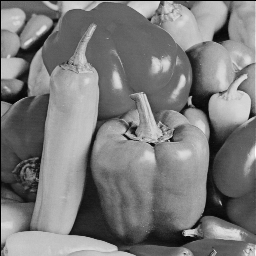}}
    \hspace{0.1in} %
\subfigure[Cameraman]{\includegraphics[width=0.45\textwidth]{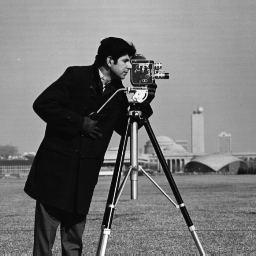}}
\end{minipage}
\hfill
\begin{minipage}{0.58\textwidth}
    \centering
\subfigure[Animals]{\includegraphics[width=0.45\textwidth]{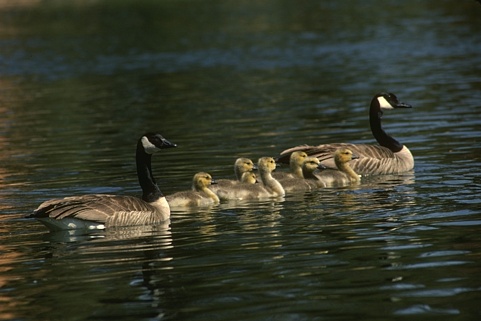}}
    \hspace{0.1in} %
\subfigure[Flowers]{\includegraphics[width=0.45\textwidth]{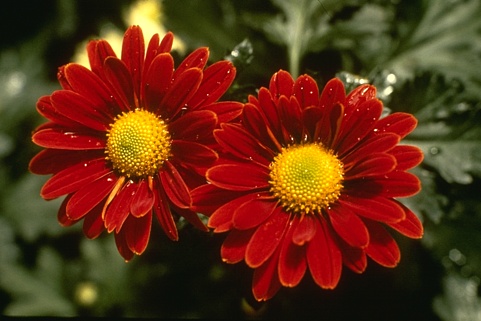}}
\end{minipage}
\caption{Tested Images in image deblurring problem.}
\end{figure}

{{
\begin{table}[t]
\caption{
Summary on the gray image deblurring problem.}\label{tab:image3}
\begin{tabular}{|c|c|c|c|c|c|}
\cline{1-6}
Image Name &  &ALPDSN & PD3O & PDFP  & Condat-Vu\\ \cline{1-6}

\multirow{4}{*}{Peppers($\lambda = $0.025)}&  iter & 318 &4340&5208& 7811 \\ \cline{2-6}
 & num &1300 & 4340 & 5208& 7811 \\\cline{2-6}
 &  ferror & 9.81e-6 & 9.99e-6 & 7.42e-5 & 1.74e-5  \\ \cline{2-6}
&time &\textbf{12.90} & 27.68 &35.40 & 50.20 \\        \cline{1-6}
\multirow{4}{*}{Peppers($\lambda = $ 0.01)}& iter & 375 &6348&7607 & 8000 \\ \cline{2-6}
 & num &1555 & 6348 & 7607 & 8000 \\\cline{2-6}
 & ferror & 9.88e-6 & 9.99e-6 & 1.34e-5 & 5.85e-5  \\ \cline{2-6}
&time &\textbf{16.65} & 45.23 &54.21& 48.08 \\        \cline{1-6}

\multirow{4}{*}{Cameraman($\lambda = $0.025)}& iter & 314 &2827&3440 &5159\\ \cline{2-6}
 & num &1184 & 2827& 3430 &5159 \\\cline{2-6}
 & ferror & 9.77e-6 & 9.99e-6 & 9.99e-6 & 9.99e-5  \\ \cline{2-6}
&time &\textbf{15.74} & 23.57 &26.97& 46.73 \\        \cline{1-6}

\multirow{4}{*}{ Cameraman($\lambda = $ 0.01)}& iter & 393 &5593&6793 &8000 \\ \cline{2-6}
 & num &1583 & 5593 & 6793 & 8000 \\\cline{2-6}
 & ferror & 9.69-6& 9.99-6 & 9.99e-6 & 2.82e-5  \\ \cline{2-6}
&time & \textbf{17.05} & 38.32 &45.70 & 61.49 \\        \cline{1-6}
\end{tabular}
\end{table}
}

{{
\begin{table}[t]
\caption{
Summary on the color image deblurring problem}\label{tab:image4}
\begin{tabular}{|c|c|c|c|c|c|}
\cline{1-6}
Image Name &  &ALPDSN & PD3O & PDFP  & Condat-Vu\\ \cline{1-6}

\multirow{4}{*}{Animals($\lambda = $0.01)}&  iter & 162 &2751&3305 & 4956 \\ \cline{2-6}
 & num &724 & 2751 & 3305 & 4956 \\\cline{2-6}
 &  ferror & 9.70e-5 & 9.99e-5 & 9.99e-5 & 9.99e-5  \\ \cline{2-6}
&time &\textbf{43.36} & 114.4 &155.9 & 199.1 \\        \cline{1-6}
\multirow{4}{*}{Animals($\lambda = $ 0.005)}& iter & 184 &4325&5199 & 7792 \\ \cline{2-6}
 & num &982 & 4325 &5199  & 7792\\\cline{2-6}
 & ferror & 9.88e-6 & 9.99e-6 & 1.34e-5 & 5.85e-5  \\ \cline{2-6}
&time &\textbf{52.96} & 165.8&238.0& 306.8\\        \cline{1-6}

\multirow{4}{*}{Flowers($\lambda = $0.01)}& iter & 192 &3369&4034  &6052\\ \cline{2-6}
 & num &911 & 3369 & 4034 &6052 \\\cline{2-6}
 & ferror & 9.91e-5 & 9.99e-5 & 9.99e-5 & 9.99e-5  \\ \cline{2-6}
&time &\textbf{71.30} & 169.2 &204.7 & 224.1 \\        \cline{1-6}

\multirow{4}{*}{ Flowers($\lambda = $ 0.005)}& iter & 263 &5197&6225 &8000 \\ \cline{2-6}
 & num &1286 & 5197 & 6225 & 8000 \\\cline{2-6}
 & ferror & 9.96-6& 9.99-6 & 9.99e-6 & 1.61e-5  \\ \cline{2-6}
&time & \textbf{105.3} &264.4 &324.1 & 347.1 \\        \cline{1-6}
\end{tabular}
\end{table}
}
\begin{figure}[htbp]
    \centering
      \subfigure[ferror vs time]{\includegraphics[scale = 0.2]{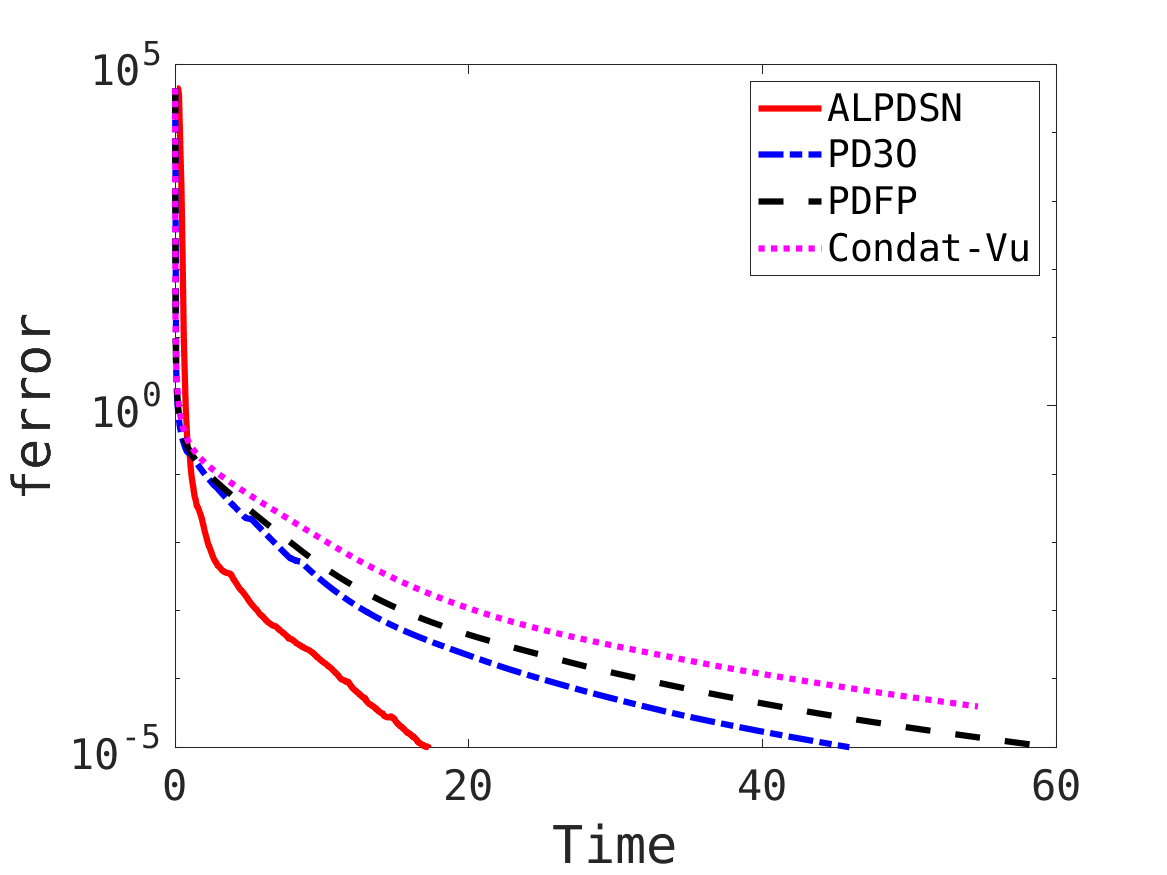}}
      \subfigure[PSNR vs time]{\includegraphics[scale = 0.2]{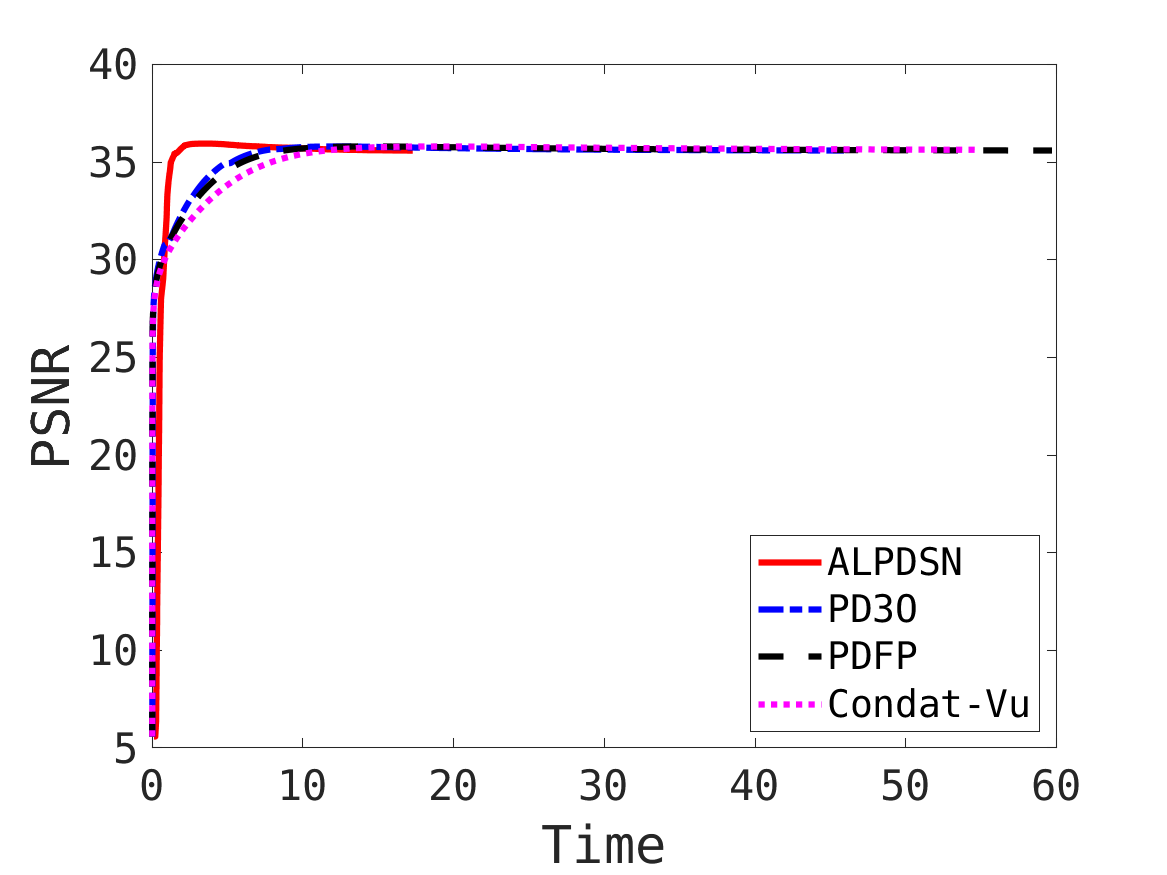}}
      \subfigure[uerror vs time ]{\includegraphics[scale = 0.2]{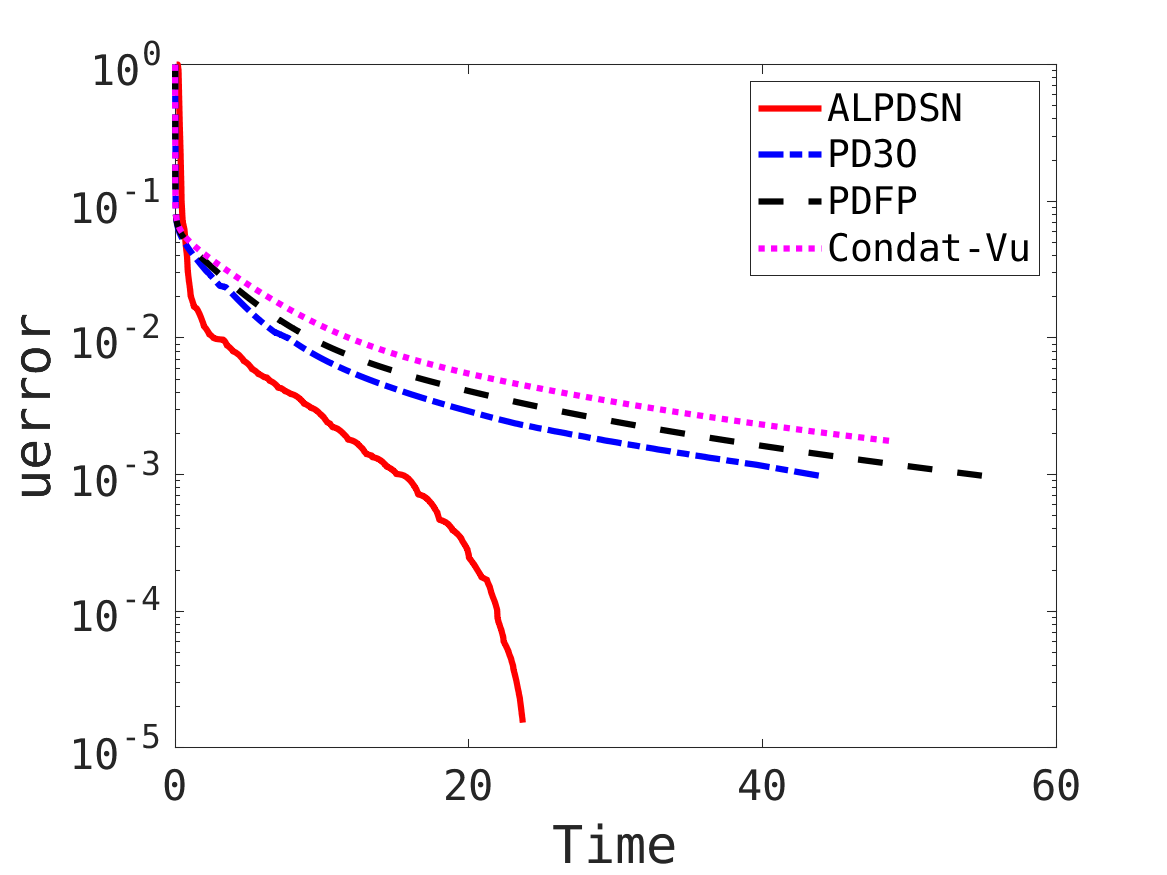}} \\
      \subfigure[ferror vs time]{\includegraphics[scale = 0.2]{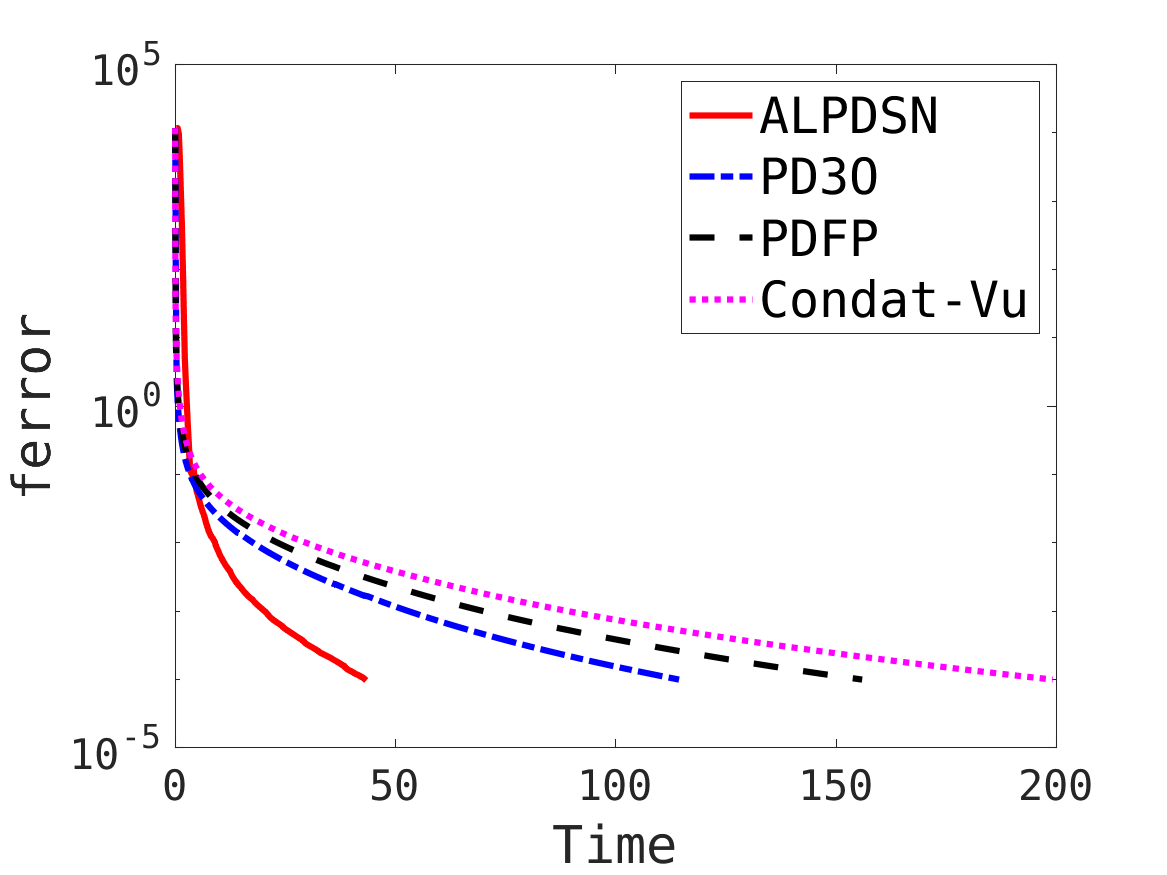}}
      \subfigure[PSNR vs time]{\includegraphics[scale = 0.2]{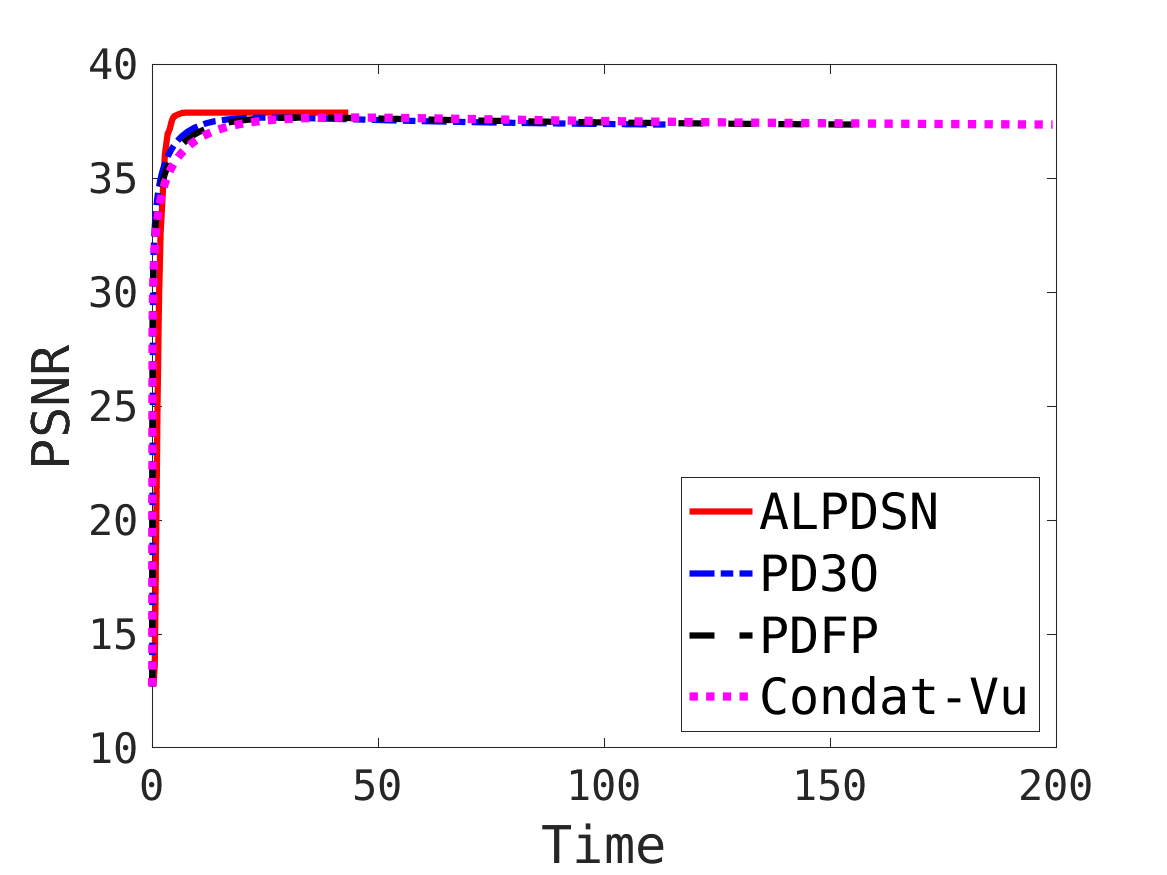}}
      \subfigure[uerror vs time ]{\includegraphics[scale = 0.2]{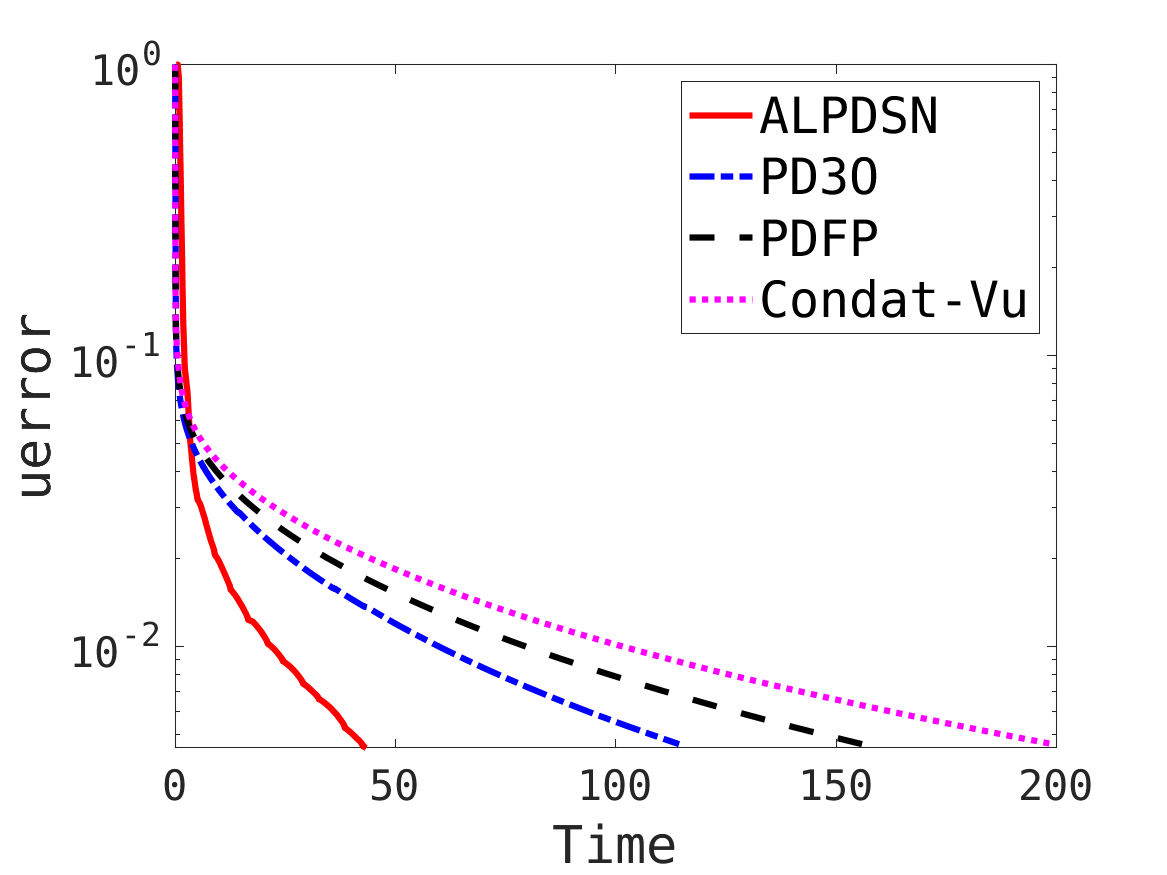}}
    \caption{Comparison between ALPDSN and other first-order methods on image deblurring. First row: Peppers with $\lambda$ = 0.01, second row: Animals with $\lambda$ = 0.01.}
      \label{fig:image4}
\end{figure}

\subsection{The CTNN model}
\subsubsection{Synthetic data}
In this subsection, we consider the CTNN model \eqref{model:dual:CTNN}.  We generate the ground-truth tensor by  $\mathcal{X}_0 = \mathcal{A} * \mathcal{B} \in \mathbb{R}^{n_1 \times n_2 \times n_3}$ with tubal rank $r$, where $\mathcal{A}$ and $\mathcal{B}$ are generated by MATLAB commands
randn($n_1$, $r$, $n_3$) and randn($r$, $n_2$, $n_3$), respectively. $\mathcal{X}_m$ is generated by the TNN model and $F$ is set the same as in \cite{zhang2019corrected}, $c$ is set to $\|\mathcal{X}_0\|_{\infty}$. The observed  vector with noise is generated by $\bm{y} = \mathcal{D}_{\Omega}(\mathcal{X}_0 + \beta \mathcal{Y}),$ where $\mathcal{Y}$ is generated by the MATLAB command randn$(n_1,n_2,n_3), $ the noise level $\beta$ is set to 0.01, $\sigma$ and $m$ are set to 1.  The sampling
ratio is 0.3.
To test the robustness of the algorithms, we choose two groups of tensors with different tubal ranks and tensor sizes. The noise level $\beta$ is set to 0.01. We compare our algorithm with sGs-ADMM \cite{li2016schur}.
In order to better compare the efficiency of the two algorithms, we choose appropriate parameters $\sigma$ based on different images for sGs-ADMM, and fix $\sigma$ to $1$ for ALPDSN. In accordance with \cite{zhang2019corrected}, we define the
 regularization parameter $\mu = \rho_1 \sigma \frac{\|\bm{v}\|}{\sqrt{m}} \sqrt{\frac{2 \log(n_1n_3+n_2n_3)}{m \tilde{n} n_3}}$ with given $\rho_1 > 0$, $\tilde{n} = \min\{n_1,n_2\} $ and $\pi_r := \max \{\pi_d, \pi_u, \pi_w ,\pi_z \},$
\begin{equation}
\begin{aligned}
\pi_d &= \frac{\|\mathcal{Y} - \mu F(\mathcal{X}_m) - \mathcal{D}^*_{\Omega}(\bm{v}) - \mathcal{W} \|_F}{1+\|\mu F(\mathcal{X}_m)\|_F},\quad \pi_u =  \frac{\|m \bm{v} -\bm{y} +\mathcal{D}_{\Omega}(\mathcal{X})  \|_F}{1+\|\bm{y}\|}, \\
\pi_w &= \frac{\|\mathcal{X}-\prox_{\delta_{\mathfrak{U}} }(\mathcal{X} - \mathcal{W}) \|_F}{1+\|\mathcal{X}\|_F +\|\mathcal{W}\|_F }, \quad \pi_y = \frac{\|\mathcal{Y} - \prox_{\delta_{\mathfrak{X}}}(\mathcal{Y} +\mathcal{X})   \|_F}{1+ \|\mathcal{Y}\|_F + \|\mathcal{X}\|_F },
\end{aligned}
\end{equation}
 and $\mathcal{W} = - \frac{1}{\sigma}\mathcal{X} -\mu F(\mathcal{X}_m) - \mathcal{D}^*_{\Omega}(\bm{u})  + \mathcal{Y}   + \frac{1}{\sigma}\prox_{\delta_{\mathfrak{U}}} \left( \mathcal{X} +  \sigma(\mu F(\mathcal{X}_m) + \mathcal{D}_{\Omega}^*(\bm{u}) - \mathcal{Y} )\right) $.
 We stop the algorithms when $\pi_r < 10^{-7}$. The results are listed in Table \ref{tab:tensor2}, where sAD1 means sGs-ADMM with $\sigma$ = 1 and similarly for others, ``time'', ``iter'', and ``error'' represent the total wall-clock time, the number of iterations and $\pi_r$, respectively. We observe that for all values of $\sigma$, ALPDSN achieves the same accuracy as sGs-ADMM in significantly less time. The detailed iterative trajectories are presented in Figure \ref{fig:tensor}, where ``RXerror'' denotes $\frac{\|\mathcal{X}-\mathcal{X}^*\|_F}{\|\mathcal{X}^*\|_F}$ and $\mathcal{X}^*$ is obtained by stopping sGs-ADMM when $\pi_r < 10^{-9}$.
As shown in Figure \ref{fig:tensor}, the convergence speed of ALPDSN is faster than sGs-ADMM at both high and low accuracy  in terms of RXerror.

\subsubsection{Real image data}
 Since CTNN is a model for noisy low-rank tensor completion, we choose two low-rank images texture1 ($194 \times 194 \times 3$ ) and texture2\footnote{\href{https://www2.eecs.berkeley.edu/Research/Projects/CS/vision/bsds/}{https://www2.eecs.berkeley.edu/Research/Projects/CS/vision/bsds/}} ($350 \times 350 \times 3$), and two natural images zebras ($321 \times 481 \times 3$) and butterfly\footnote{\href{https://www.kaggle.com/datasets/kmader/siim-medical-images?resource=download}{ https://www.kaggle.com/datasets/kmader/siim-medical-images}} ($321 \times 481 \times 3$ ), see Figures \ref{fig:tensor-texture-image}(a)-(d). The pixels of the tested images are normalized to $[0, 255]$ , $\sigma$ and $\beta$ are set to 1 and 0.01,  respectively. For sGs-ADMM, we choose the best parameter  $\sigma$ to ensure its performances for different images. To fully demonstrate the  comparison of sGs-ADMM, we stop the iteration when $\pi_r <  10^{-6}.$ The sampled and restored images for CTNN model can be seen from Figures \ref{fig:tensor-texture-image}(e)-(h).  We can see from  Table \ref{tab:tensor-image}
that ALPDSN is faster than sGs-ADMM for all different values of $\sigma$.
 The detailed iteration process of  synthetic data, Texture1 and Zebra can be seen in Figure \ref{fig:tensor}.

\begin{figure}[htbp]
\centering
\begin{minipage}{0.95\textwidth}
    \centering
\subfigure[ Texture1]{\includegraphics[width=0.22\textwidth]{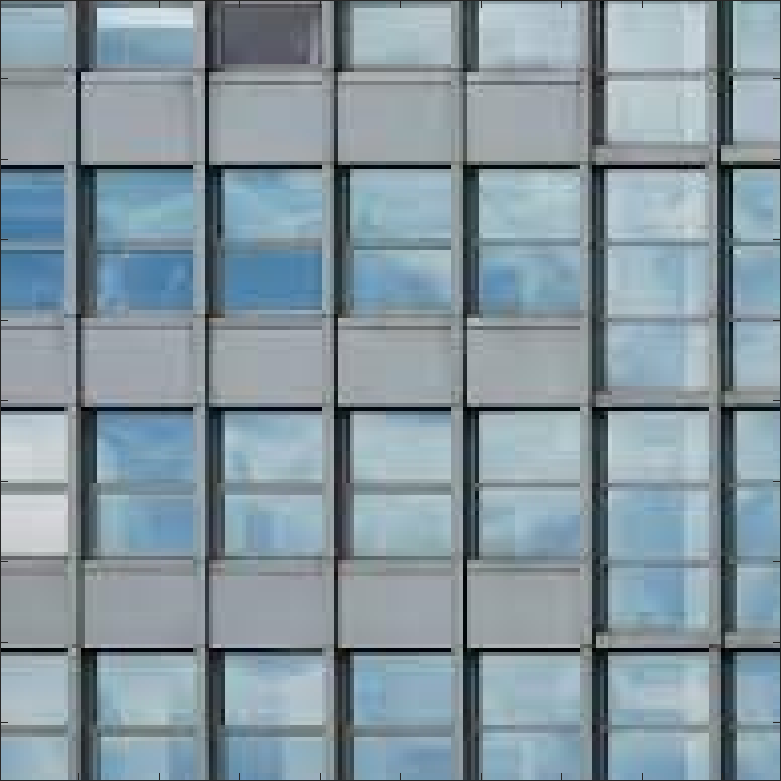}}
    \hspace{0.05in} %
\subfigure[ Texture2]{\includegraphics[width=0.22\textwidth]{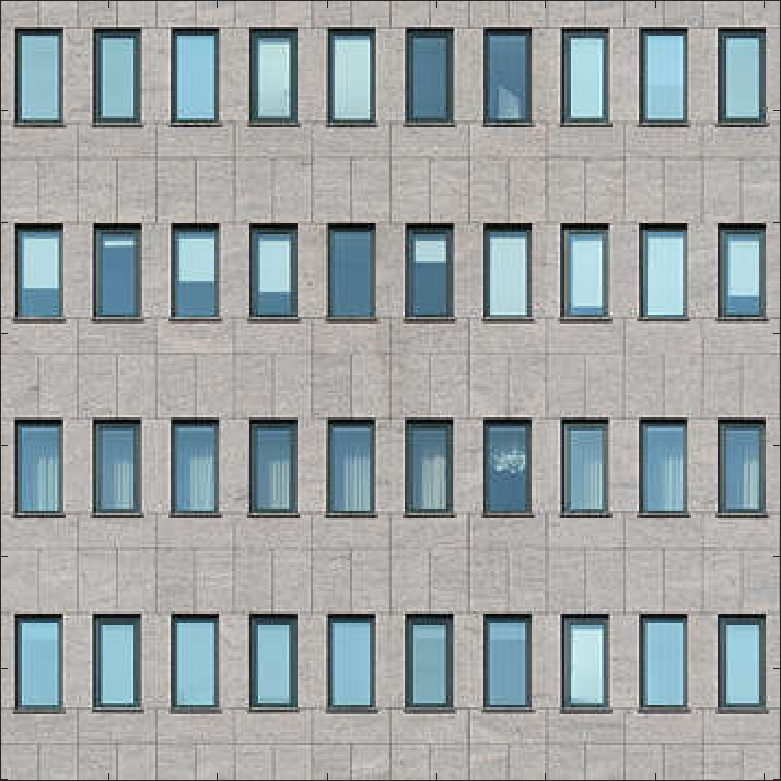}}
\hspace{0.05in}
\subfigure[Zebra]{\includegraphics[width=0.22\textwidth,height=0.22\textwidth]{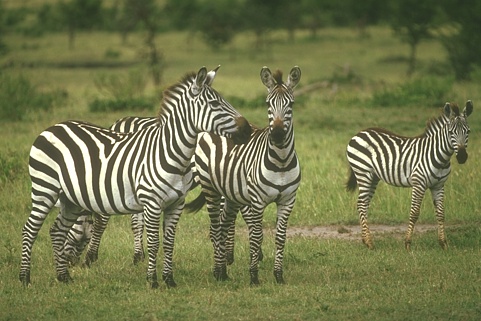}}
\hspace{0.05in}
\subfigure[Butterfly]{\includegraphics[width=0.22\textwidth,height=0.22\textwidth]{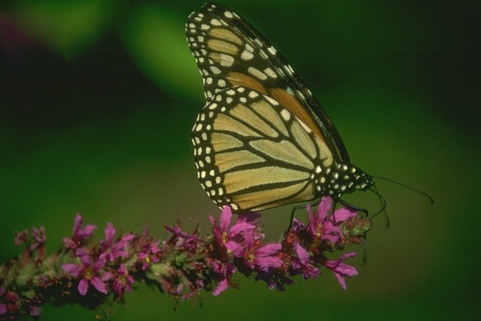}}
\hspace{0.05in} \\
\subfigure[Sampled Texture1]{\includegraphics[width=0.22\textwidth]{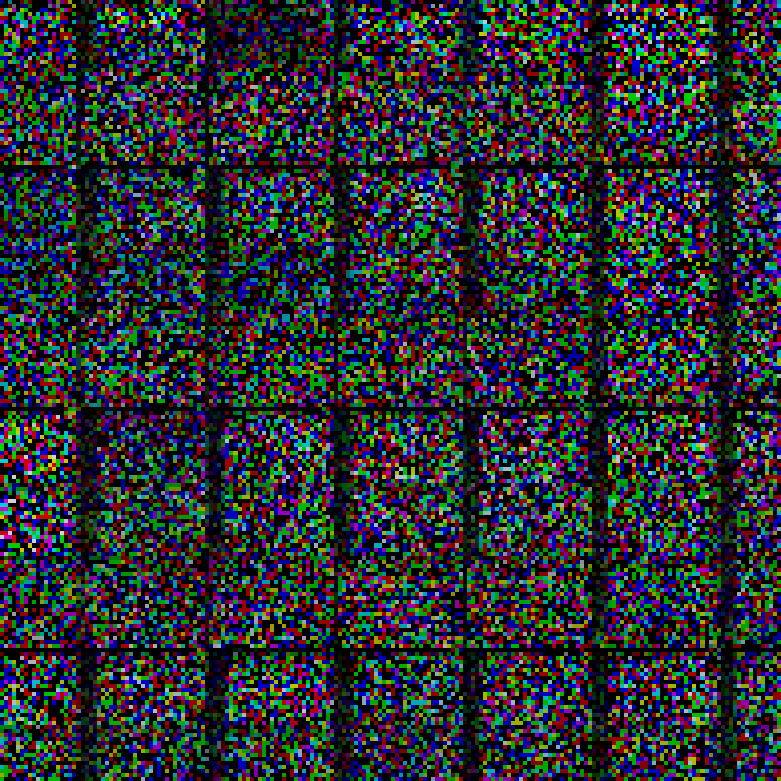}}
    \hspace{0.05in} %
\subfigure[Restored Texture1]{\includegraphics[width=0.22\textwidth]{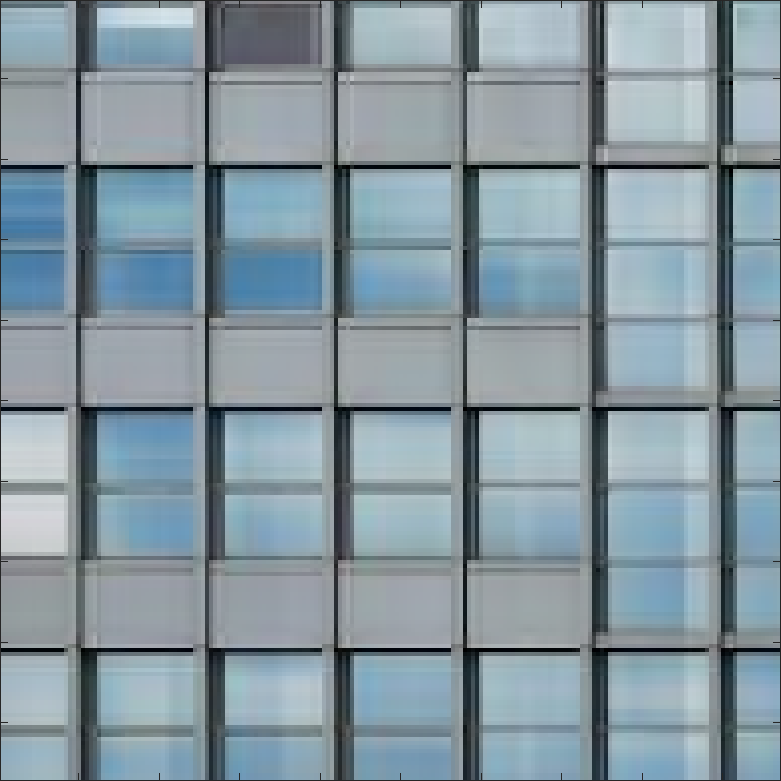}}
    \hspace{0.05in} %
\subfigure[Sampled Zebra]{\includegraphics[width=0.22\textwidth,height=0.22\textwidth]{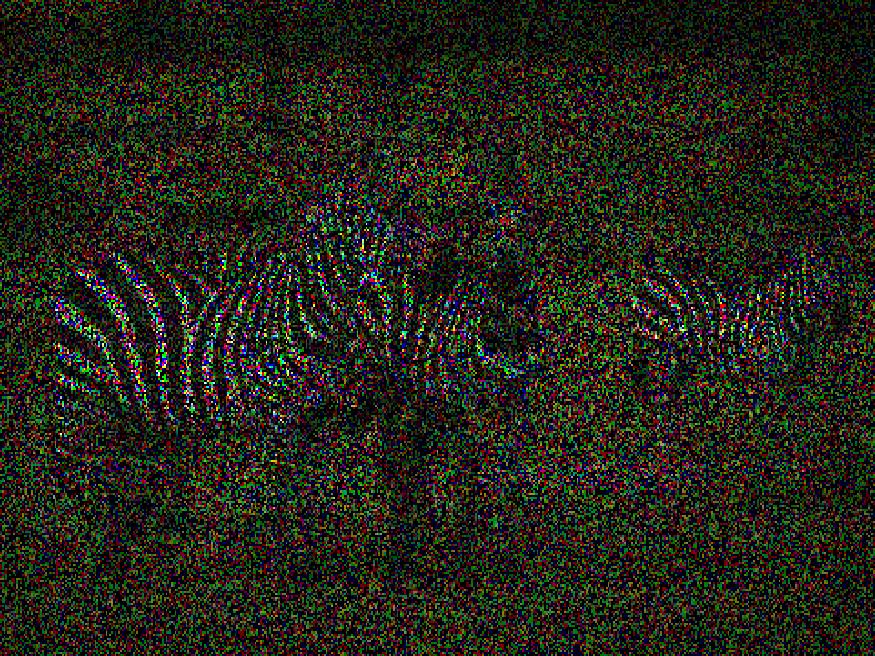}}
\hspace{0.05in}
\subfigure[Restored Zebra]{\includegraphics[width=0.22\textwidth,height=0.22\textwidth]{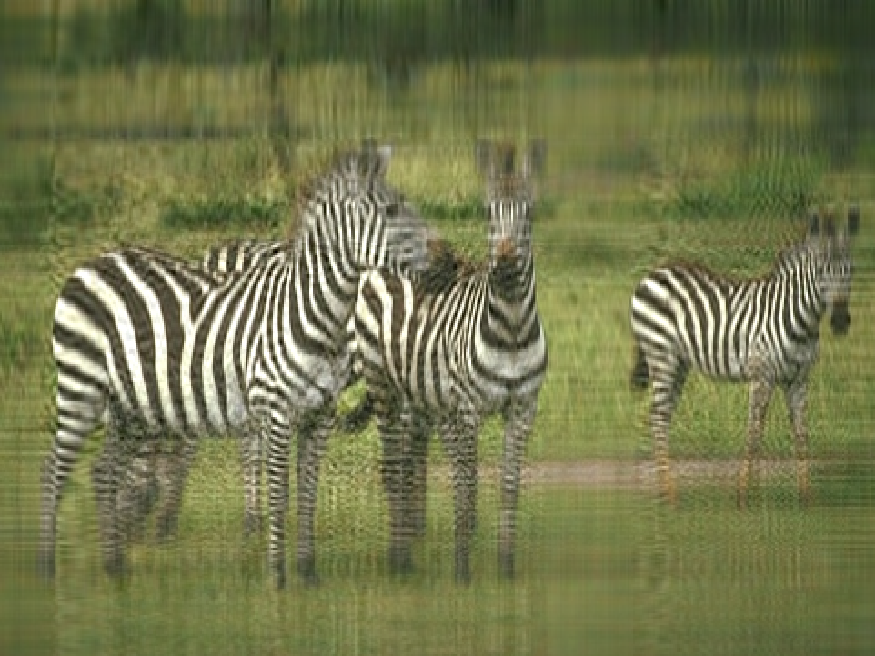}}
\hspace{0.05in}
\end{minipage}
\label{fig:tensor-texture-image}
\caption{Results on CTNN model. First row: tested images, second row: sampled and restored images of Texture1 and Zerbra.}
\end{figure}

{
\begin{table}[t]
\centering
\caption{
Numerical results of CTNN on synthetic data }\label{tab:tensor2}
\begin{tabular}{|c|c|c|c|c|c|c|}
\cline{1-7}
Tensor & rank &&ALPDSN & sAD1 & sAD5& sAD10\\ \cline{1-7}
 \multirow{6}{*}{200 $\times$ 250 $\times$ 10}&
  \multirow{3}{*}{5} & time&\textbf{20.44} & 49.77 &  44.10 & 66.74\\ \cline{3-7}
   & & iter& 15 & 281 &252& 384\\ \cline{3-7}
 & & error& 4.84e-8 &9.42e-8 &9.32e-8 & 9.90e-8\\ \cline{2-7}
 & \multirow{3}{*}{10} & time&\textbf{18.97} &  78.07 & 32.48& 32.37  \\  \cline{3-7}
  & & iter& 19 & 451 & 180 & 181 \\ \cline{3-7}
 & & error& 3.86e-8 & 9.73e-8 &9.32e-8 & 9.19e-8\\  \cline{1-7}
\multirow{6}{*}{300 $\times$ 300 $\times$ 5}&
  \multirow{3}{*}{10} & time&\textbf{17.88} &47.20& 41.81 & 64.27 \\ \cline{3-7}
  & & iter& 16 & 242&215 & 322\\ \cline{3-7}
 & & error& 1.76e-8& 8.81e-8 & 9.95e-8 & 9.96e-8\\ \cline{2-7}
 & \multirow{3}{*}{15} & time&\textbf{15.10} & 63.57 & 28.99& 30.68 \\  \cline{3-7}
  & & iter& 18 & 325 &149 & 157\\ \cline{3-7}
 & & error& 8.75e-8& 9.67e-8 & 9.23e-8 & 9.63e-8\\
 \cline{1-7}
  \end{tabular}
\end{table}

\begin{figure}[htbp]
    \centering
      \subfigure[Synthetic data]{\includegraphics[scale = 0.2]{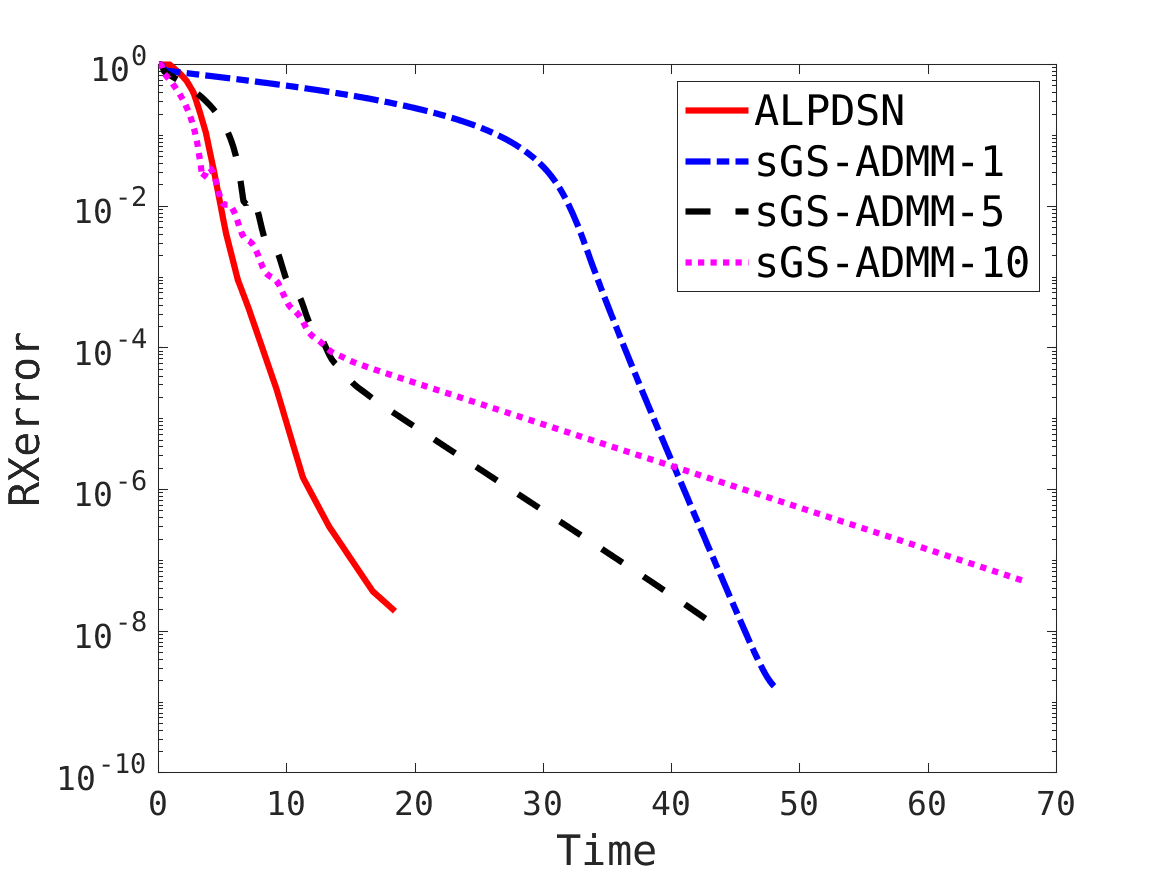}}
      \subfigure[Texture1]{\includegraphics[scale = 0.2]{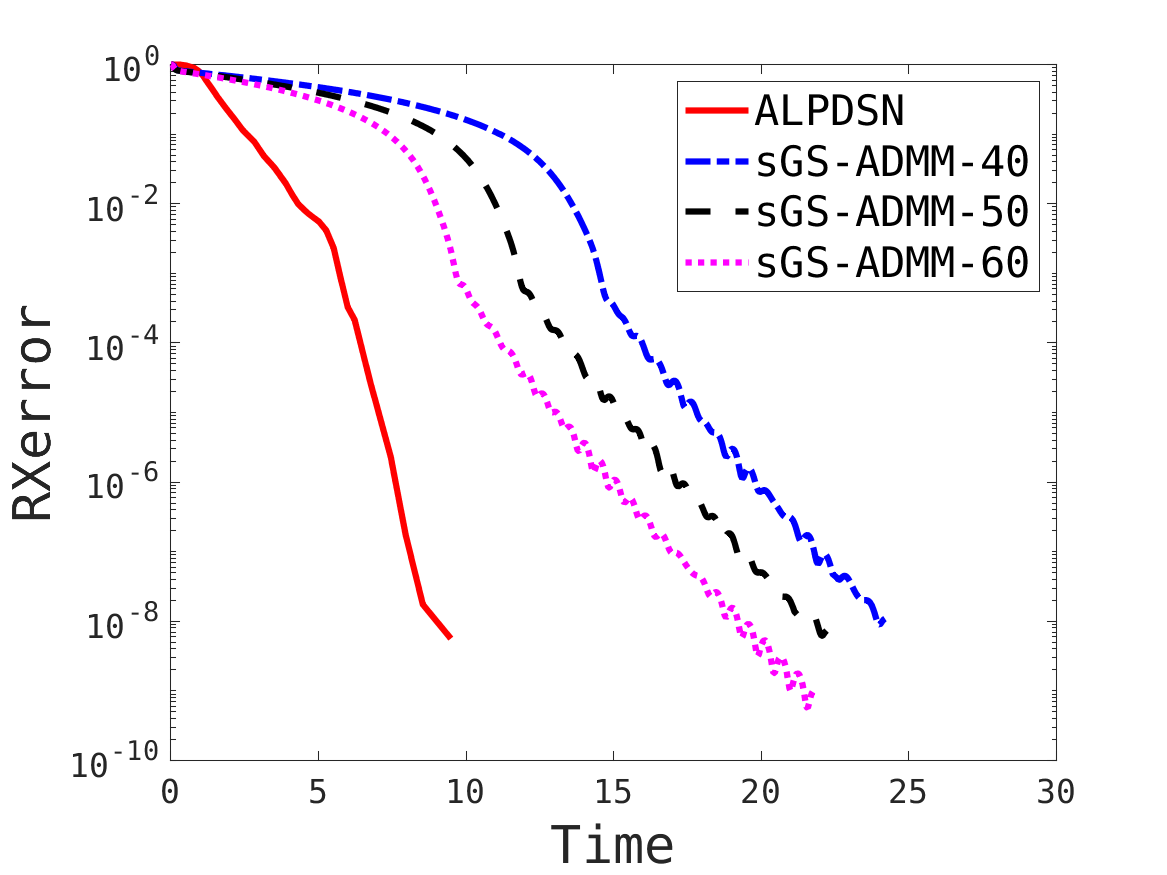}}
            \subfigure[Zebra]{\includegraphics[scale = 0.2]{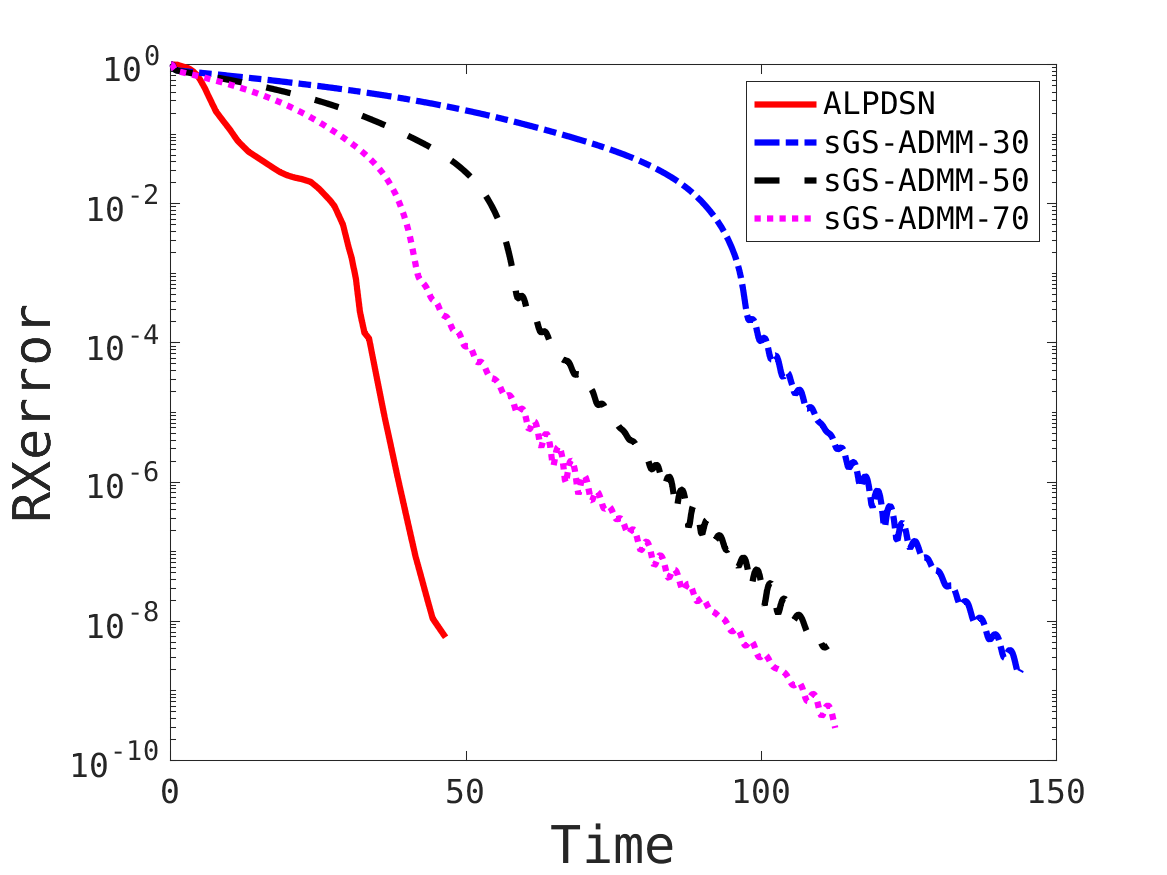}}
      \caption{Comparison of ALPDSN and sGs-ADMM on the CTNN model.}
      \label{fig:tensor}
\end{figure}

\section{Conclusions} \label{sec7}

This paper introduces the ALPDSN, an AL primal dual semismooth Newton method for efficiently solving multi-block convex composite optimization problems. To navigate the challenges of inherent nonsmoothness and multi-block structures, we leverage the concept of the Moreau envelope, resulting in an equivalent smooth saddle point problem. Subsequently, a semismooth system of nonlinear equations is presented to characterize the optimality condition of the original problem. Integrating this technique with a novel nonmonotone line search strategy, our proposed ALPDSN performs purely regularized semismooth Newton steps. The global convergence and locally superlinear convergence are established by carefully elaborating the regularity conditions of the problem and the converged solution, e.g., twice differentiability, convexity, partial smoothness, piecewise linear-quadratic, and the strict complementarity. Numerical comparisons with state-of-the-art methods validate the high efficiency and robustness of our ALPDSN on the image restoration with two regularization terms and the CTNN model.

Although the tested applications are limited to problems with three-block variables due to the restriction on the length of this paper, they possess large dimensionality and complex structures, especially in the tensor-based CTNN problem. Our approach is highly adaptable and can be extended to address more complex multi-block scenarios. It is worth mentioning that our algorithm also yields favorable results for doubly nonnegative semidefinite programming, which will be presented comprehensively in a companion paper.

\begin{table}[t]
\centering
\caption{
Numerical results on low-rank textures and natural images .}\label{tab:tensor-image}
\begin{tabular}{|c|c|c|c|c|c|}
\cline{1-6}
Image name  &&ALPDSN & sAD 40 & sAD 50& sAD 60\\ \cline{1-6}
\multirow{3}{*}{Texture 1}&
 time&\textbf{9.53} & 24.02 & 21.93&  21.45 \\ \cline{2-6}
   &  iter& 33&  538 & 487 & 490\\ \cline{2-6}
 &  error& 3.98e-7 & 9.88e-7& 9.90e-7 & 9.31e-7\\  \cline{1-6}
 Image name  &&ALPDSN & sAD 50 & sAD 80& sAD 100\\ \cline{1-6}
 \multirow{3}{*}{Texture 2}&
   time&\textbf{29.67} &105.7& 104.6& 108.7\\ \cline{2-6}
   &  iter&32 & 614 &   598& 633\\ \cline{2-6}
 &  error&4.38e-8 &9.75e-7 & 9.75e-7& 9.58e-7\\ \cline{1-6}
 Image name  &&ALPDSN & sAD 30 & sAD 50& sAD 70\\ \cline{1-6}
\multirow{3}{*}{Zebra}&
   time&\textbf{46.51} &144.0& 111.5& 112.5\\ \cline{2-6}
   &  iter&38 & 685 &   530& 536\\ \cline{2-6}
 &  error&2.81e-7 &8.80e-7 & 9.45e-7& 9.58e-7\\ \cline{1-6}
 Image name  &&ALPDSN & sAD 30 & sAD 50& sAD 70\\ \cline{1-6}
 \multirow{3}{*}{Butterfly}&
   time&\textbf{38.80} &112.9& 93.75& 108.2\\ \cline{2-6}
   &  iter&29 & 538 &   449& 516\\ \cline{2-6}
 &  error&3.63e-07 &9.59e-7 & 9.84e-7& 9.96e-7\\ \cline{1-6}
  \end{tabular}
\end{table}

\section*{Acknowledgments}
We would like to thank Xiongjun Zhang for sharing their codes on the CTNN problem.

\bibliographystyle{siamplain}
\bibliography{ref}

\begin{thebibliography}{10}

\bibitem{beck2017first}
{\sc A.~Beck}, {\em First-order Methods in Optimization}, SIAM, 2017.

\bibitem{boyd2011distributed}
{\sc S.~Boyd, N.~Parikh, E.~Chu, B.~Peleato, J.~Eckstein, et~al.}, {\em
  Distributed optimization and statistical learning via the alternating
  direction method of multipliers}, Foundations and Trends{\textregistered} in
  Machine learning, 3 (2011), pp.~1--122.

\bibitem{chambolle2011first}
{\sc A.~Chambolle and T.~Pock}, {\em A first-order primal-dual algorithm for
  convex problems with applications to imaging}, Journal of Mathematical
  Imaging and Vision, 40 (2010), pp.~120--145.

\bibitem{chan2008constraint}
{\sc Z.~X. Chan and D.~Sun}, {\em Constraint nondegeneracy, strong regularity,
  and nonsingularity in semidefinite programming}, SIAM Journal on
  optimization, 19 (2008), pp.~370--396.

\bibitem{chen2017efficient}
{\sc L.~Chen, D.~Sun, and K.-C. Toh}, {\em An efficient inexact symmetric
  {G}auss--{S}eidel based majorized {ADMM} for high-dimensional convex
  composite conic programming}, Mathematical Programming, 161 (2017),
  pp.~237--270.

\bibitem{chen2016primal}
{\sc P.~Chen, J.~Huang, and X.~Zhang}, {\em A primal-dual fixed point algorithm
  for minimization of the sum of three convex separable functions}, Fixed Point
  Theory and Applications, 2016 (2016), pp.~1--18.

\bibitem{chen2020proximal}
{\sc S.~Chen, S.~Ma, A.~M.-C. So, and T.~Zhang}, {\em Proximal gradient method
  for nonsmooth optimization over the {S}tiefel manifold}, SIAM Journal on
  Optimization,  (2020).

\bibitem{chen2020alternating}
{\sc S.~Chen, S.~Ma, L.~Xue, and H.~Zou}, {\em An alternating manifold proximal
  gradient method for sparse principal component analysis and sparse canonical
  correlation analysis}, INFORMS Journal on Optimization,  (2020).

\bibitem{condat2013primal}
{\sc L.~Condat}, {\em A primal--dual splitting method for convex optimization
  involving {L}ipschitzian, proximable and linear composite terms}, Journal of
  optimization theory and applications, 158 (2013), pp.~460--479.

\bibitem{davis2016convergence}
{\sc D.~Davis and W.~Yin}, {\em Convergence rate analysis of several splitting
  schemes}, Splitting Methods in Communication, Imaging, Science, and
  Engineering,  (2017), p.~115.

\bibitem{drori2015simple}
{\sc Y.~Drori, S.~Sabach, and M.~Teboulle}, {\em A simple algorithm for a class
  of nonsmooth convex--concave saddle-point problems}, Operations Research
  Letters, 43 (2015), pp.~209--214.

\bibitem{drusvyatskiy2014optimality}
{\sc D.~Drusvyatskiy and A.~S. Lewis}, {\em Optimality, identifiability, and
  sensitivity}, Mathematical Programming, 147 (2014), pp.~467--498.

\bibitem{facchinei2003finite}
{\sc F.~Facchinei and J.-S. Pang}, {\em Finite-dimensional variational
  inequalities and complementarity problems}, Springer, 2003.

\bibitem{fan2004inexact}
{\sc J.~Fan and J.~Pan}, {\em Inexact {L}evenberg-{M}arquardt method for
  nonlinear equations}, Discrete and Continuous Dynamical Systems-B, 4 (2004),
  pp.~1223--1232.

\bibitem{hintermuller2010semismooth}
{\sc M.~Hinterm{\"u}ller}, {\em Semismooth {N}ewton methods and applications},
  Department of Mathematics, Humboldt-University of Berlin,  (2010).

\bibitem{hu2022local}
{\sc J.~Hu, T.~Tian, S.~Pan, and Z.~Wen}, {\em On the local convergence of the
  semismooth {N}ewton method for composite optimization}, arXiv preprint
  arXiv:2211.01127,  (2022).

\bibitem{jiang1995local}
{\sc H.~Y. Jiang and L.~Qi}, {\em Local uniqueness and convergence of iterative
  methods for nonsmooth variational inequalities}, Journal of Mathematical
  Analysis and Applications, 196 (1995), pp.~314--331.

\bibitem{karush1939minima}
{\sc W.~Karush}, {\em Minima of functions of several variables with
  inequalities as side constraints}, M. Sc. Dissertation. Dept. of Mathematics,
  Univ. of Chicago,  (1939).

\bibitem{kilmer2011factorization}
{\sc M.~E. Kilmer and C.~D. Martin}, {\em Factorization strategies for
  third-order tensors}, Linear Algebra and its Applications, 435 (2011),
  pp.~641--658.

\bibitem{komodakis2015playing}
{\sc N.~Komodakis and J.-C. Pesquet}, {\em Playing with duality: An overview of
  recent primal-dual approaches for solving large-scale optimization problems},
  IEEE Signal Processing Magazine, 32 (2015), pp.~31--54.

\bibitem{kuhn2014nonlinear}
{\sc H.~W. Kuhn and A.~W. Tucker}, {\em Nonlinear programming}, in Traces and
  emergence of nonlinear programming, Springer, 2014, pp.~247--258.

\bibitem{lam2021semi}
{\sc X.-Y. Lam, D.~Sun, and K.-C. Toh}, {\em Semi-proximal augmented
  {L}agrangian-based decomposition methods for primal block-angular convex
  composite quadratic conic programming problems}, Informs Journal on
  Optimization, 3 (2021), pp.~254--277.

\bibitem{lewis2002active}
{\sc A.~S. Lewis}, {\em Active sets, nonsmoothness, and sensitivity}, SIAM
  Journal on Optimization, 13 (2002), pp.~702--725.

\bibitem{li2016schur}
{\sc X.~Li, D.~Sun, and K.-C. Toh}, {\em A {S}chur complement based
  semi-proximal {ADMM} for convex quadratic conic programming and extensions},
  Mathematical Programming, 155 (2016), pp.~333--373.

\bibitem{li2018highly}
{\sc X.~Li, D.~Sun, and K.-C. Toh}, {\em A highly efficient semismooth {N}ewton
  augmented {L}agrangian method for solving {L}asso problems}, SIAM Journal on
  Optimization, 28 (2018), pp.~433--458.

\bibitem{li2018efficiently}
{\sc X.~Li, D.~Sun, and K.-C. Toh}, {\em On efficiently solving the subproblems
  of a level-set method for fused {L}asso problems}, SIAM Journal on
  Optimization, 28 (2018), pp.~1842--1866.

\bibitem{li2018qsdpnal}
{\sc X.~Li, D.~Sun, and K.-C. Toh}, {\em Q{S}{D}{P}{N}{A}{L}: A two-phase
  augmented {L}agrangian method for convex quadratic semidefinite programming},
  Mathematical Programming Computation, 10 (2018), pp.~703--743.

\bibitem{li2020asymptotically}
{\sc X.~Li, D.~Sun, and K.-C. Toh}, {\em An asymptotically superlinearly
  convergent semismooth {N}ewton augmented {L}agrangian method for linear
  programming}, SIAM Journal on Optimization, 30 (2020), pp.~2410--2440.

\bibitem{li2018semismooth}
{\sc Y.~Li, Z.~Wen, C.~Yang, and Y.-x. Yuan}, {\em A semismooth {N}ewton method
  for semidefinite programs and its applications in electronic structure
  calculations}, SIAM Journal on Scientific Computing, 40 (2018),
  pp.~A4131--A4157.

\bibitem{lu2019tensor}
{\sc C.~Lu, J.~Feng, Y.~Chen, W.~Liu, Z.~Lin, and S.~Yan}, {\em Tensor robust
  principal component analysis with a new tensor nuclear norm}, IEEE
  Transactions on Pattern Analysis and Machine Intelligence, 42 (2019),
  pp.~925--938.

\bibitem{mifflin1977semismooth}
{\sc R.~Mifflin}, {\em Semismooth and semiconvex functions in constrained
  optimization}, SIAM Journal on Control and Optimization, 15 (1977),
  pp.~959--972.

\bibitem{milzarek2014semismooth}
{\sc A.~Milzarek and M.~Ulbrich}, {\em A semismooth {N}ewton method with
  multidimensional filter globalization for $\ell_1$-optimization}, SIAM
  Journal on Optimization, 24 (2014), pp.~298--333.

\bibitem{milzarek2019stochastic}
{\sc A.~Milzarek, X.~Xiao, S.~Cen, Z.~Wen, and M.~Ulbrich}, {\em A stochastic
  semismooth {N}ewton method for nonsmooth nonconvex optimization}, SIAM
  Journal on Optimization, 29 (2019), pp.~2916--2948.

\bibitem{mordukhovich2017geometric}
{\sc B.~S. Mordukhovich and N.~M. Nam}, {\em Geometric approach to convex
  subdifferential calculus}, Optimization, 66 (2017), pp.~839--873.

\bibitem{nemirovski2004prox}
{\sc A.~Nemirovski}, {\em Prox-method with rate of convergence
  $\mathcal{O}(1/t)$ for variational inequalities with {L}ipschitz continuous
  monotone operators and smooth convex-concave saddle point problems}, SIAM
  Journal on Optimization, 15 (2004), pp.~229--251.

\bibitem{nocedal1999numerical}
{\sc J.~Nocedal and S.~J. Wright}, {\em Numerical Optimization}, Springer
  Series in Operations Research and Financial Engineering, Springer, New York,
  second~ed., 2006.

\bibitem{robinson1981some}
{\sc S.~M. Robinson}, {\em Some continuity properties of polyhedral
  multifunctions}, Springer, 1981.

\bibitem{rockafellar2009variational}
{\sc R.~T. Rockafellar and R.~J.-B. Wets}, {\em Variational Analysis},
  vol.~317, Springer Science \& Business Media, 2009.

\bibitem{sun2002semismooth}
{\sc D.~Sun and J.~Sun}, {\em Semismooth matrix-valued functions}, Mathematics
  of Operations Research, 27 (2002), pp.~150--169.

\bibitem{toh2004solving}
{\sc K.-C. Toh}, {\em Solving large scale semidefinite programs via an
  iterative solver on the augmented systems}, SIAM Journal on Optimization, 14
  (2004), pp.~670--698.

\bibitem{xiao2018regularized}
{\sc X.~Xiao, Y.~Li, Z.~Wen, and L.~Zhang}, {\em A regularized semi-smooth
  {N}ewton method with projection steps for composite convex programs}, Journal
  of Scientific Computing, 76 (2018), pp.~364--389.

\bibitem{xu2015augmented}
{\sc J.~Xu, S.~Zhu, Y.~C. Soh, and L.~Xie}, {\em Augmented distributed gradient
  methods for multi-agent optimization under uncoordinated constant stepsizes},
  in 54th IEEE Conference on Decision and Control (CDC), IEEE, 2015,
  pp.~2055--2060.

\bibitem{yan2018new}
{\sc M.~Yan}, {\em A new primal--dual algorithm for minimizing the sum of three
  functions with a linear operator}, Journal of Scientific Computing, 76
  (2018), pp.~1698--1717.

\bibitem{yue2019family}
{\sc M.-C. Yue, Z.~Zhou, and A.~M.-C. So}, {\em A family of inexact {SQA}
  methods for non-smooth convex minimization with provable convergence
  guarantees based on the {Luo--Tseng} error bound property}, Mathematical
  Programming, 174 (2019), pp.~327--358.

\bibitem{zhang2019corrected}
{\sc X.~Zhang and M.~K. Ng}, {\em A corrected tensor nuclear norm minimization
  method for noisy low-rank tensor completion}, SIAM Journal on Imaging
  Sciences, 12 (2019), pp.~1231--1273.

\bibitem{zhao2022robust}
{\sc X.~Zhao, M.~Bai, D.~Sun, and L.~Zheng}, {\em Robust tensor completion:
  Equivalent surrogates, error bounds, and algorithms}, SIAM Journal on Imaging
  Sciences, 15 (2022), pp.~625--669.

\bibitem{zhao2010newton}
{\sc X.-Y. Zhao, D.~Sun, and K.-C. Toh}, {\em A {N}ewton-{C}{G} augmented
  {L}agrangian method for semidefinite programming}, SIAM Journal on
  Optimization, 20 (2010), pp.~1737--1765.

\bibitem{zhou2005superlinear}
{\sc G.~Zhou and K.-C. Toh}, {\em Superlinear convergence of a {N}ewton-type
  algorithm for monotone equations}, Journal of optimization theory and
  applications, 125 (2005), pp.~205--221.

\bibitem{zhu2008efficient}
{\sc M.~Zhu and T.~Chan}, {\em An efficient primal-dual hybrid gradient
  algorithm for total variation image restoration}, UCLA CAM Report, 34 (2008),
  pp.~8--34.

\end{thebibliography}
\end{document}